\UseRawInputEncoding
\documentclass[12pt,a4paper]{article}
\usepackage{amsmath,amssymb,amsthm,url,amsfonts,enumerate}
\usepackage{graphicx}
\usepackage{hyperref}
\hypersetup{
   colorlinks   = true, 
   urlcolor     = blue, 
   linkcolor    = blue, 
   citecolor   = red 
}
\graphicspath{{talk/}}
\input{xy}\xyoption{all}
\usepackage[a4paper,left=2.5cm,right=1.7cm,top=2cm,bottom=2cm]{geometry}
\usepackage[usenames,dvipsnames]{color}

\long\def\comment#1\endcomment{}

\def\Z{{\mathbb Z}} \def\R{{\mathbb R}}  

\def\lk{\mathop{\fam0 lk}}
\def\Int{\mathop{\fam0 Int}}

\DeclareMathOperator{\diag}{diag}
\newcommand{\id}{\mathop{\mathrm{id}}}
\def\forg{\mathop{\fam0 forg}}
\def\sgn{\mathop{\fam0 sgn}}
\def\t{\widetilde}

\DeclareMathOperator{\supp}{supp}
\newcommand{\iprod}{\raisebox{-0.5ex}{\scalebox{1.8}{$\cdot$}}\,}

\theoremstyle{plain}
\newtheorem{theorem}{Theorem}[section]
\newtheorem{lemma}[theorem]{Lemma}
\newtheorem{proposition}[theorem]{Proposition}

\theoremstyle{definition}
\newtheorem{problem}[theorem]{Problem}
\newtheorem{remark}[theorem]{Remark}

\begin{document}

\title{A user's guide to the topological Tverberg conjecture
\footnote{
Research supported by the Russian Foundation for Basic Research Grant No. 15-01-06302,
by Simons-IUM Fellowship and by the D. Zimin's Dynasty Foundation Grant.
\newline
Subsection \ref{s:oz} is written jointly with R. Karasev.
I am grateful to S. Avvakumov, P. Blagojevi\'c, V. Buchstaber, G. Kalai, R. Karasev, I. Mabillard, S. Melikhov, A. Ryabichev, M. Tancer, T. Tolozova and U. Wagner for useful remarks, to H. McFadden for improving English, and to I. Mabillard and U. Wagner for allowing me to use some figures.}}

\author{A. Skopenkov
\footnote{Moscow Institute of Physics and Technology,
and Independent University of Moscow.
Email: \texttt{skopenko@mccme.ru}.
\texttt{http://www.mccme.ru/\~{}skopenko}
}
}

\date{}

\maketitle

\begin{abstract}
The \emph{topological Tverberg conjecture} was considered a central unsolved problem of topological combinatorics.
The conjecture asserts that {\it for every integers $r,d$ and any continuous map $f\colon\Delta\to \R^d$ of the $(d+1)(r-1)$-dimensional simplex there are pairwise disjoint faces $\sigma_1,\ldots,\sigma_r\subset\Delta$ such that $f(\sigma_1)\cap \ldots \cap f(\sigma_r)\ne\emptyset$. }
The conjecture was proved for a prime power $r$.
Recently counterexamples for other $r$ were found.
Analogously, the \emph{$r$-fold van Kampen-Flores conjecture} holds for a prime power $r$ but does not hold for other $r$.
The arguments form a beautiful and fruitful interplay between combinatorics, algebra and topology.
We present a simplified exposition accessible to non-specialists in the area.
We also mention some recent developments and open problems.
\end{abstract}


\tableofcontents

\section{Introduction}\label{s:intr}

\subsection{Description of main results}\label{s:intrm}

\begin{remark}[Linear Radon, Tverberg and van Kampen-Flores theorems]\label{r:rt}
We start with classical results motivating the topological Tverberg conjecture.
For history, more motivation and related problems see \cite{Zi11, BBZ, Sk18}, \cite[\S1-\S3]{BZ16}, \cite[\S6]{Sk}.
This remark is formally not used in what follows.

A subset of $\R^d$ is {\it convex},
if for any two points from this subset the segment joining these two points is in this subset.
The {\it convex hull} of a subset $X\subset\R^d$ is the minimal convex set that contains $X$.

(a) The well-known \emph{Radon theorem} asserts that
{\it for every integer $d>0$ any $d+2$ points in $\R^d$ can be decomposed into two groups such that the convex hulls of the groups intersect.}
See inductive geometric proofs in \cite{Ko18, Pe72, RRS} and standard algebraic proof e.g. in \cite{BZ16, RRS}.

(b) The linear version of well-known \emph{van Kampen-Flores theorem} asserts that
{\it for every integer $k>0$ from any $2k+3$ points in $\R^{2k}$ one can choose two disjoint $(k+1)$-tuples whose convex hulls intersect.}
For $k=1$ this implies the linear non-planarity of the complete graph $K_5$ on 5 vertices.
For $k>1$ this implies the linear non-realizability in $\R^{2k}$ of the complete $(k+1)$-homogeneous hypergraph on $2k+3$ vertices.

For an odd-dimensional version of the van Kampen-Flores theorem (Conway-Gordon-Sachs-Lovasz-Schrijver-Taniyama)  see $(VKF_d)$ of \S\ref{s:cgst} for $d$ odd.
An inductive geometric proof (involving the odd-dimensional version) is analogous to the proof for $k=2$ \cite{Sk14}.
For algebraic proof see \cite{BM15}.

\begin{figure}[h]
\centerline{\includegraphics[width=8cm]{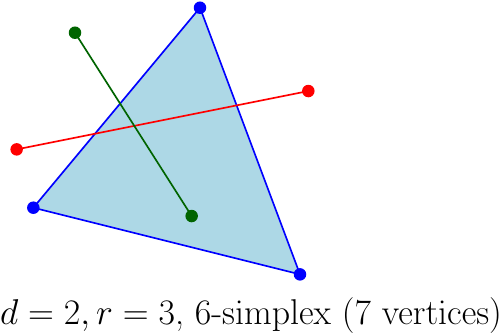}\qquad \includegraphics[width=6cm]{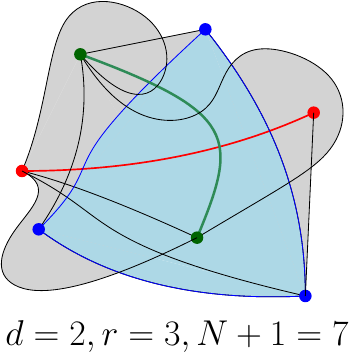}}
\caption{Linear and topological Tverberg Theorems for $r=3$}\label{f:lt}
\end{figure}

(c) The Radon Theorem is generalized by \emph{Tverberg theorem} stating that
{\it for every integers $d,r>0$ any $(d+1)(r-1)+1$ points in $\R^d$ can be decomposed into $r$ groups
such that all the $r$ convex hulls of the groups have a common point.} See fig. \ref{f:lt}.

For a motivated exposition of a well-known proof see \cite{RRS}.


Clearly, every $(d+1)(r-1)$ points in general position in $\R^d$ (or vertices of a $d$-dimensional simplex taken with multiplicity $r-1$) do not satisfy this property.
So if one is bothered by remembering the number $(d+1)(r-1)+1$ in the Tverberg theorem, one can remember that
this is the minimal number such that general position and calculation of the dimension of the intersection do not produce a counterexample.

(d) The linear version of \emph{$r$-fold van Kampen-Flores conjecture} asserts that
{\it for every integers $k,r>0$ from any $(r-1)(kr+2)+1$ points in $\R^{kr}$ one can choose $r$ pairwise disjoint $(k(r-1)+1)$-tuples all the $r$ convex hulls of the tuples have a common point.}

This is true for a prime $r$ \cite{Sa91g} and even for a prime power $r$ \cite{Vo96v}; this is an open problem for other $r$ \cite[beginning of \S2]{Fr17}.

Take in $\R^{kr}$ the vertices of a $kr$-dimensional simplex and its center.
Either take every of these $kr+2$ points with multiplicity $r-1$ or for every point take close $r-1$ points in general position.
We obtain $(r-1)(kr+2)$ points in $\R^{kr}$ such that for any $r$ pairwise disjoint $(k(r-1)+1)$-tuples all the $r$ convex hulls of the tuples do not have a common point.\footnote{For $r=3k=3$ cf. \cite[Example 6.7.4]{Ma03}: `It is not known whether such triangles can always be found for 9 points in $\R^3$'.}
So if one is bothered by remembering the number $(r-1)(kr+2)+1$ in the conjecture, one can remember that
this is the minimal number such that the above trivial construction does not produce a counterexample.
\end{remark}

Denote by $\Delta_N$ the $N$-dimensional simplex.

\smallskip
{\bf ($TR_d$) Topological Radon Theorem.} \cite{BB} For any continuous map $\Delta_{d+1}\to\R^d$ there are two disjoint faces whose images intersect.

\smallskip
{\bf ($VKF_{2k}$) Van Kampen-Flores Theorem.} For any continuous map
$\Delta_{2k+2}\to\R^{2k}$ there are two disjoint $k$-dimensional faces whose images intersect.

\smallskip
These results generalize the theorems of Remark \ref{r:rt}.ab above
(see also $(VKF_d)$ of \S\ref{s:cgst} for $d$ odd).
They are nice in themselves, and are also interesting because they are corollaries of the celebrated Borsuk-Ulam Theorem \ref{vf-bu}, of which the topological Radon Theorem is also a simplicial version.
They can also be proved inductively (Remark \ref{r:rela}) or using {\it van Kampen (Radon) number} \cite{Sk18}.

The particular case $r=2$ of the Constraint Lemma \ref{p:redu} below asserts that
$(TR_{2k+1})\Rightarrow (VKF_{2k})$.
For similar implications see Remark \ref{r:rela}.


The PL (piecewise-linear) versions of assertions $(TR_d)$ and $(VKF_{2k})$ are as interesting and non-trivial as the stated topological versions (this is similar to Remark \ref{r:ae}.c below).

The well-known \emph{topological Tverberg conjecture} was raised by E. Bajmoczy and I. B\'ar\'any~\cite{BB} and
H. Tverberg~\cite[Problem~84]{GS}.
It was considered a central unsolved problem of topological combinatorics.
The conjecture asserts that
{\it for every integers $r,d$ and any continuous map $f\colon\Delta_{(d+1)(r-1)}\to \R^d$ there are pairwise disjoint faces $\sigma_1,\ldots,\sigma_r\subset\Delta_{(d+1)(r-1)}$ such that
$f(\sigma_1)\cap \ldots \cap f(\sigma_r)\ne\emptyset$.}

This conjecture generalizes both the Tverberg theorem and the topological Radon theorem above.
In my opinion, one of the main reasons why this conjecture was considered important is that its investigation was one of the first examples of a problem in topological combinatorics, for which $\Z_2$-actions are insufficient and one needs actions of more complicated groups.

The conjecture was proved for $r$ a prime by I.~B\'{a}r\'{a}ny, S. Shlosman, A.~Sz{\H{u}}cs \cite{BSS},
and then for $r$ a prime power by M.~\"Ozaydin and A. Volovikov \cite{Oz, vo96} (Theorem \ref{t:tvepo} below).
A stronger result is considered in \cite{FS20}.

Recently and somewhat unexpectedly, it turned out that there are counterexamples for $r$ not a prime power.
For these counterexamples papers \cite{Oz, Gr10, BFZ14, Fr15, MW15} (by M.~\"Ozaydin; M. Gromov; P. Blagojevi\'c, F. Frick, and G. Ziegler; F. Frick; I. Mabillard and U. Wagner) are important (Theorem \ref{t:tvepo} below).\footnote{\label{f:discl} Here we do not discuss  the credits for the counterexample, but
we provide accurate references to each step of the proof (Remark \ref{r:hystor} below and \S\ref{s:plan})
so that the reader could make his/her own opinion.}

Counterexamples were first constructed for $d\ge3r+1$ and then for $d\ge2r+1$ \cite{AMSW} (the case $d=3r$ is considered in \cite{MW15}, but I did not check the argument, cf. \cite[Remark 3.1.a]{AMSW}).
The lowest dimensional counterexample
is an almost $6$-embedding $\Delta_{70}\to\R^{13}$, cf. \cite[\S1.3]{AMSW}.
For stronger counterexamples see \cite[\S5]{BBZ} (with a simpler exposition in Remark \ref{l:lon}.c) and \cite[Theorem 1]{AKS}.
This conjecture is still open for $r$ not a prime power and $d\le2r$, in particular, for $d\le12$ and so for $d=2$.
For $d=2r$ see \cite[Remark 3.1.a]{AMSW}.



Analogously, \emph{$r$-fold van Kampen-Flores conjecture} asserts that
{\it for every integers $r,k>0$ and any continuous map $f\colon\Delta_{(kr+2)(r-1)}\to \R^{kr}$ there are pairwise disjoint $k(r-1)$-dimensional faces $\sigma_1,\ldots,\sigma_r\subset\Delta_{(kr+2)(r-1)}$ such that
$f(\sigma_1)\cap \ldots \cap f(\sigma_r)\ne\emptyset$.}\footnote{Let us call this statement a conjecture to emphasize that it is not always true, and that it is analogous to the topological Tverberg conjecture.
This was not stated as a conjecture before positive result for $r$ a prime power \cite{Vo96v} appeared.}


This is true for a prime $r$ \cite{Sa91g} and even for a prime power $r$ \cite{Vo96v} (Theorem \ref{t:vkfr} below).
This is false for other $r$ \cite{MW15} (Theorem \ref{t:vkfr} below and Remark \ref{r:hystor}.b).
Counterexamples were first constructed for $k\ge3$ \cite{MW15} and then for $k=2$ \cite{AMSW}.
For stronger counterexamples see \cite[Theorem 4]{AKS}.
For $k=1$ and $r$ not a prime power this conjecture is still open.

The profs are based on a beautiful and fruitful interplay among combinatorics, algebra and topology.
Namely, the arguments for positive solutions (\S\ref{s:posi}) use equivariant algebraic topology of configuration spaces.
The arguments for counterexamples (\S\ref{s:coex}) use, besides that, elementary combinatorics
and geometric topology (eliminating multiple intersections by an analogue of Whitney trick).

\smallskip
{\bf Why this survey might be interesting.}
We present a simplified explanation of the arguments (proving and disproving
the topological Tverberg and the $r$-fold van Kampen-Flores conjectures).
Our exposition is accessible to non-specialists in the area.
We present a proof of $r$-fold Whitney trick (the Mabillard-Wagner Theorem \ref{t:mawa}) from \cite{AMSW}, which is simpler than the original one.
Despite being shorter, our exposition of other parts is not an alternative proof but just a different
presentation, making clear the structure and avoiding sophisticated language.
When we use a certain theory, we explicitly state a specific result proved with the help of the theory (\S\ref{s:posi}, \S\ref{s:oz}) but in terms not involving the theory.
This makes the {\it application} of the result  accessible to non-specialists in the area, and it also makes both the proof of the result and the theory more accessible.
So if a statement is clear to a reader but the proof uses unfamiliar notions, the reader can skip the proof and just make use of the statement.
We present accurate description of references on counterexamples to the topological Tverberg conjecture (leaving to a reader distribution of credits), and specific critical remarks on misleading descriptions.

\subsection{Main definitions and statements}\label{s:state}

\begin{figure}[h]
\centerline{\includegraphics[width=7cm]{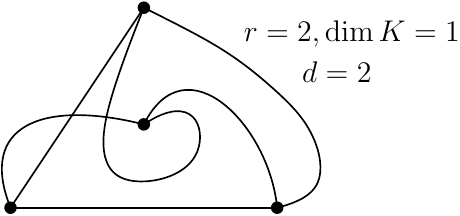}}
\caption{An almost 2-embedding $K_4\to\R^2$}\label{f:ak4}
\end{figure}

A map  $f\colon K\to \R^d$ from a union $K$ of some faces of a simplex is an {\bf almost $r$-embedding} if $f(\sigma_1)\cap \ldots \cap f(\sigma_r)=\emptyset$ whenever $\sigma_1,\ldots,\sigma_r$ are pairwise disjoint faces of $K$. See fig.  \ref{f:ak4}.

We mostly omit `continuous' for maps and group actions.

\begin{remark}\label{r:ae}
(a) In this language the topological Tverberg conjecture and the $r$-fold van Kampen-Flores conjecture state that

$\bullet$ for every integers $r,d$ there are no almost $r$-embeddings $\Delta_{(d+1)(r-1)}\to\R^d$.

$\bullet$ for every integers $r,k$ there are no almost $r$-embeddings of
the union of $k(r-1)$-faces of $\Delta_{(kr+2)(r-1)}$ in $\R^{kr}$.

(b) The notion of an almost $2$-embedding implicitly arose in studies realizability of graphs and complexes/hypergraphs (see ($VKF_{2k}$) of \S\ref{s:intrm} and \cite{We67}).
It was explicitly formulated in the Freedman-Krushkal-Teichner paper \cite{FKT}.

(c) Any sufficiently small perturbation of an almost $r$-embedding is again an almost $r$-embedding.
So the existence of a {\it continuous} almost $r$-embedding is equivalent to the existence of a {\it PL} almost $r$-embedding, and to the existence of a {\it general position PL} almost $r$-embedding.

Topological Tverberg conjecture means that this is also equivalent to the existence of a {\it linear} almost $r$-embedding, because the conjecture for {\it linear} maps is Tverberg theorem of Remark \ref{r:rt}.c.
Thus the conjecture is a higher-dimensional $r$-fold analogue of F\'ary Theorem stating that {\it if a graph is planar, then it is linearly planar}, cf.  \cite{PW}.

(d) The studies of almost $r$-embeddings should be compared to studies of {\it immersions without $r$-tuple points}.
E.g.

$\bullet$ a closed 2-manifold $M$ admits an immersion to $\R^3$ without triple points if and only if $\chi(M)$ is even (folklore);

$\bullet$ any 3-manifold admits an immersion to $\R^4$ with a quadruple point \cite{Fr78}.
\end{remark}

\begin{theorem}[\cite{BSS, Oz, vo96}; proved in \S\ref{s:posi}]\label{t:tvepo}
If $r$ is a prime power and $d>0$ is an integer, then there are no almost $r$-embeddings $\Delta_{(d+1)(r-1)}\to\R^d$.
\end{theorem}

\begin{theorem}[see \cite{Oz, Gr10, BFZ14, Fr15, MW15} and footnote \ref{f:discl}]\label{t:tve}
If $r$ is not a prime power and $d\ge3r+1$, then there is an almost $r$-embedding
$\Delta_{(d+1)(r-1)}\to\R^d$.
\end{theorem}

For $d=3r+1$ Theorem \ref{t:tve} follows from Theorem \ref{t:vkfr-} and the Constraint Lemma \ref{p:redu},
both below.
The case $d>3r+1$ is implied by the case $d=3r+1$ and the following Remark \ref{l:lon}.b.

\begin{remark}[using the join construction]\label{l:lon}
(a) For two maps $f:\Delta_a\to B^p$ and $g:\Delta_b\to B^q$ define the {\it join}
$$
f*g:\Delta_{a+b+1}=\Delta_a*\Delta_b \to B^p * B^p =B^{p+q+1}
\quad\text{by the formula}\quad
(f*g)(\lambda x \oplus \mu y):= \lambda f(x)\oplus\mu f(y).
$$
A join of almost $r$-embeddings is an almost $r$-embedding.

(b) {\it For every integers $r,N,d>0$ if there is an almost $r$-embedding $\Delta_N\to\R^d$,
then there is an almost $r$-embedding $\Delta_{N+r-1}\to \R^{d+1}$} \cite[Proposition~2.5]{Lo}.

Indeed, assume that $f:\Delta_N\to\R^d$ is an almost $r$-embedding.
The map $g$ of $\Delta_{r-2}$ to a point $B^0$ is an almost $r$-embedding.
By (a), the join $f*g:\Delta_{N+r-1}=\Delta_N*\Delta_{r-2}\to B^d*B^0=B^{d+1}$ is an almost $r$-embedding.

(c) {\it For every integers $r,a,d,k>0$ if there is an almost $r$-embedding $\Delta_a\to\R^d$, then there is an almost $r$-embedding $\Delta_{k(a+1)-1}\to\R^{k(d+1)-1}$} \cite[Lemma 5.2]{BFZ}.

Indeed, by (a) the $k$-fold join power of an almost $r$-embedding $\Delta_a\to B^d$ is an almost $r$-embedding $\Delta_{k(a+1)-1}\to B^{k(d+1)-1}$.
\end{remark}

\begin{theorem}[$r$-fold van Kampen-Flores Theorem]\label{t:vkfr}
If $r$ is a prime power and $k>0$ is an integer, then there are no almost $r$-embeddings of
the union of $k(r-1)$-faces of $\Delta_{(kr+2)(r-1)}$ in $\R^{kr}$.
\end{theorem}

This result was proved in \cite{Vo96v} analogously to Theorem \ref{t:tvepo};
this result also follows from Theorem \ref{t:tvepo} by the Constraint Lemma \ref{p:redu} below.


\begin{theorem}\label{t:vkfr-}
If $r$ is not a prime power and $k\ge3$ is any integer, then there is an almost $r$-embedding of
the union of $k(r-1)$-faces of $\Delta_{(kr+2)(r-1)}$ in $\R^{kr}$.
\end{theorem}

This result follows from Theorems \ref{t:mawa} and \ref{t:ozmawa} of Mabillard-Wagner and \"Ozaydin
\cite{MW15, Oz}.
Its analogue for {\it codimension $k=2$} is correct \cite[Theorem 1.2.a]{AMSW}.

\begin{lemma}[Gromov-Blagojevi\'c-Frick-Ziegler Constraint Lemma; see Remark \ref{r:hystor}.ac]\label{p:redu}
For every integers $r,k>0$ if there is an almost $r$-embedding of the union of $k(r-1)$-faces of
$\Delta_{(kr+2)(r-1)}$ in $\R^{kr}$, then there is an almost $r$-embedding $\Delta_{(kr+2)(r-1)}\to\R^{kr+1}$.
\end{lemma}

\begin{proof}
(In order to make this argument more accessible we consider the case $r=6$ and $k=3$.
The general case is analogous.)

Take an almost 6-embedding $\Delta_{100}^{15}\to\R^{18}$ of the union of $15$-dimensional faces of $\Delta_{100}$.
Extend it arbitrarily to a map $f:\Delta_{100}\to \R^{18}$.
Denote by $\rho(x)$ the distance from $x\in \Delta_{100}$ to $\Delta_{100}^{15}$.
It suffices to prove that $f\times\rho:\Delta_{100}\to \R^{19}$ is an almost 6-embedding.

Suppose to the contrary that 6 points $x_1,\ldots, x_6\in\Delta_{100}$ lie in pairwise disjoint faces and are mapped to the same point under $f\times\rho$.
Dimension of one of those faces does not exceed $\frac{101}{6} - 1$, so it is at most 15.
W.l.o.g. this is the first face, hence $\rho(x_1)=0$.
Then $\rho(x_2) = \ldots = \rho(x_6)=\rho(x_1)= 0$, i.e. $x_1,\ldots, x_6\in\Delta_{100}^{15}$.
Now the condition $f(x_1) = \ldots = f(x_6)$ contradicts the fact that $f|_{\Delta_{100}^{15}}$ is an almost $6$-embedding.
\end{proof}

The remaining sections of this paper are independent on each other.
Sections 2, 3 and 4 use important notions introduced in \S\ref{s:posiex}.

\begin{remark}[Historical]\label{r:hystor}
(a) {\it `The topological Tverberg theorem, whenever available, implies the van Kampen-Flores theorem'}
\cite[2.9.c, p. 445, lines --1 and --2]{Gr10}.
This is the Constraint Lemma \ref{p:redu}.
What M. Gromov called `theorem', we call `conjecture'.
See (c) for the rediscovery of this lemma.

This lemma was proved in \cite[2.9.c, p.446, 2nd paragraph]{Gr10} by a beautiful combinatorial trick reproduced in the current paper.\footnote{\label{f:gromov} This lemma was by no means among main results of the long and deep paper \cite{Gr10}.
The statement in \cite[2.9.c, p. 446]{Gr10} of a more general result is hard to read.
So note that
\newline
$\bullet$ the number $T_{top}(q,n)$ is the number of topological Tverberg partitions, see \cite[p. 444 above and the third paragraph of 2.9.a]{Gr10};
\newline
$\bullet$ instead of $T_{top}(q,n)$ there should be $T_{top}(q,n+1)$;
\newline
$\bullet$
The Constraint Lemma \ref{p:redu} is obtained by taking $q:=r$, $k:=k(r-1)$, $n:=kr$, $N=N_{qn}=N_{nq}:=(kr+2)(r-1)$, and using the implication `$T_{top}(q,n+1)>0\ \Rightarrow \ m(q,n)>0$' not stronger inequality $[VKF]_q$.
\newline
See Remarks \ref{r:known}.a, \ref{r:bbz}.be, \ref{r:bs}.ad  and \ref{r:sh}.adefg.}



(b) {\it A counterexample to the $r$-fold van Kampen-Flores conjecture for $r$ not a prime power}
(Theorem \ref{t:vkfr-})\footnote{\label{f:gropro} The counterexample is a weaker version of a positive answer to the following open problem \cite[end of 2.9.c, p. 446]{Gr10} (what he denotes by $q$, we denote by $r$):
{\it `It also seems unknown, for any $q=6,10,\ldots$ which \emph{is not} a prime power, whether
every compact $k$-dimensional topological space $X$, where $n=qk/(q-1)$ is an integer, admits a $(q-1)$-to-1 map to $\R^n$.'}}
was implicitly obtained by I. Mabillard and U. Wagner in \cite{MW14, MW15}\footnote{\label{f:mw15} Publications in computer science conferences (like \cite{MW14}) are not peer review publications, see \cite{Fo}.
The paper \cite{Fr15} refers to \cite{MW14} as to `extended abstract'.
The arXiv preprint \cite{MW15} was not published by 2020; as far as I know, it remains unpublished by January, 2022.
A complete proof (slightly different from \cite{MW15}) is published in the survey \cite{Sk18u}
(of course attributing the result to Mabillard-Wagner).}
by mentioning the \"Ozaydin  Theorem \ref{theorem:ozaydin-product}, and proving that the \"Ozaydin  Theorem \ref{theorem:ozaydin-product} implies Theorem \ref{t:vkfr-}.
Indeed,

$\bullet$ the \"Ozaydin  Theorem \ref{theorem:ozaydin-product} states that some equivariant map exists for $r$ not a prime power;\footnote{\label{f:ozaydin} `{\it On the other hand, \"Ozaydin [22] also showed that for every $r$ that is not a prime power, there \emph{does} exist an equivariant map $(\Delta_N)^r_\Delta\to_{\Sigma_r}S^{d(r-1)-1}$}'
\cite[p. 173, the paragraph before Theorem 3]{MW14}. The notation used for the equivariant map is introduced in \S\ref{s:posiex}.
What \cite{MW14} denotes by $(\Delta_N)^r_\Delta$, we denote by $\t{\Delta_N}^r$.
What \cite{MW14} calls `a map without $r$-Tverberg point', we call `an almost $r$-embedding'.
See also \cite[p. 174, Motivation \& Future Work, 2nd paragraph]{MW14}.
This \"Ozaydin Theorem was not stated in unpublished paper \cite{Oz} but is attributed to \cite{Oz} in \cite{MW14}, and indeed easily follows from \cite{Oz} as it was shown in \cite[Proof of Corollary 3]{Fr15}.}

$\bullet$ \cite[Theorem 3]{MW14} states that for codimension at least 3 the existence of the above equivariant map implies the existence of an almost $r$-embedding as in Theorem \ref{t:vkfr-}.\footnote{This is Theorem \ref{t:mmw} for $(r-1)d=rs$. See Remark \ref{r:mw} for the history.}

Failure of the $r$-fold van Kampen-Flores conjecture for $r$ not a prime power was first {\it explicitly stated} by F. Frick \cite{Fr15}.\footnote{The above shows that the following phrase of \cite[p. 2]{Fr15},
\cite[p. 1]{Fr15o} is misleading: `{\it we first derive / obtain counterexamples to $r$-fold versions of the van Kampen-Flores theorem}'. Cf. Remark \ref{r:known}.c and \cite[bottom of p. 1]{BZ16}: `{\it In a spectacular recent development, Isaac Mabillard and Uli Wagner [33] [34] have developed an $r$-fold version of the classical `Whitney trick' (cf. [50]), which yields the failure of the generalized Van Kampen-Flores theorem when $r\ge6$ is not a prime power.}'
See Remarks \ref{r:bbz}.e and \ref{r:bs}.b.}


(c) {\it Since $O$ and $O\Rightarrow\overline{VKF}$ and $TTC\Rightarrow VKF$ imply $\overline{TTC}$,
by (a) and (b) papers \cite{Oz, Gr10, MW14} together give counterexample to the topological Tverberg conjecture} (so that there remained the non-trivial task of writing versions of \cite{Oz, MW14} that could be accepted to peer-review journals).
However, neither Mabillard and Wagner nor the topological combinatorics community were aware of (a) before 2016.
This is surprising because the next part \cite[2.9.e]{Gr10} of Gromov's paper was discussed during the problem session at 2012 Oberwolfach Workshop on Triangulations.\footnote{In spite of that, I am sure that the rediscovery of (a) described below is independent.
This situation shows how important it is for efficient research to thoroughly study the relevant literature.}
So this community saw a serious problem with the Mabillard-Wagner approach to a counterexample to the topological Tverberg conjecture: maps from the $(d+1)(r-1)$-simplex to $\R^d$ do not satisfy the codimension $\ge3$ restriction required for \cite[Theorem 3]{MW14};
these maps actually have negative codimension.
F. Frick \cite[proof of Theorem 4]{Fr15} realized that this problem can be overcome by (a).
He did this by rediscovering (a), not by finding the reference \cite[2.9.c]{Gr10}.
In fact, (a) was implicitly rediscovered earlier by Blagojevi\'{c}-Frick-Ziegler
\cite[Lemma 4.1.iii and 4.2]{BFZ14}.\footnote{The purpose of \cite{BFZ14} was not (a) but rather its generalizations, `constraint method', not required for disproof of the topological Tverberg conjecture.
The Constraint Lemma \ref{p:redu} of (a) is not explicitly stated in \cite{BFZ14} but is implicitly proved
in the proof of other results.
The lemma is also not explicitly stated in \cite{Fr15, Fr15o, BZ16}.
Thus the lemma is proved implicitly, separately for $r$ a prime power \cite[Lemma 4.1.iii and 4.2]{BFZ14},
\cite[\S4.1]{BZ16} and for other $r$ \cite[proof of Theorem 4]{Fr15}, \cite[\S5]{BZ16}, although neither case  uses the fact that $r$ is a prime power or not.
See also remarks listed at the end of footnote \ref{f:gromov}.}




(d) An accurate description of references on counterexample to the topological Tverberg conjecture
is given here and in \S\ref{s:plan}.
For other descriptions of such references
see \cite{BBZ},
\cite[\S1]{BFZ}, \cite[\S1 and beginning of \S5]{BZ16}, \cite[\S1.1]{JVZ}, \cite[\S1.1]{Si16},
\cite[\S1.1]{AMSW},
\cite[bottom of p. 1]{Fr17}, \cite{BS17}, \cite[\S3, \S4]{Sh18}.
The specific remarks presented \S\ref{s:app} show that the description in
\cite{BBZ, BFZ, BZ16, Si16, Fr17, BS17, Sh18} is misleading.
\end{remark}

\subsection{Complexes and configuration spaces}\label{s:posiex}

{\bf Definitions of a complex, a polyhedron, $D^n$, $S^{n-1}$ and $m$-connectedness.}
A {\bf complex} is a collection of closed faces (=simplices) of some simplex.\footnote{This is an abbreviation of `an (abstract) finite simplicial complex'.
In combinatorial terms, a complex is a collection of subsets of a finite set such that
if a subset $A$ is in the collection, then each subset of $A$ is in the collection,
A close but different notion widely studied in combinatorics is {\it hypergraph}.
}
A {\it $k$-complex} is a complex containing at most $k$-dimensional simplices.
The {\it body} (or geometric realization) $|K|$ of a complex $K$ is the union of simplices of $K$.
Thus continuous or piecewise-linear (PL) maps $|K|\to\R^d$ and continuous maps $|K|\to S^m$ are defined.
Below we abbreviate $|K|$ to $K$; no confusion should arise.
However, we note that the property of being an almost embedding depends on $K$ and not just on $|K|$.

A {\bf polyhedron} is a subset of $\R^d$ homeomorphic to the body of certain complex.

Denote $D^n:=[0,1]^n$ and $S^{n-1}:=\partial D^n$.

A polyhedron $K$ is called {\bf $m$-connected} if for each $j=0,1,\ldots,m$ every map $f:S^j\to K$ extends over $D^{j+1}$.





\smallskip
{\bf Some results and open problems.}
The following result is deduced in \cite{MW15, AMSW} from Theorem \ref{t:mawa} below.

\begin{theorem}[\cite{MW15, AMSW}]\label{t:algal}
For every fixed $s,d,r$ such that $(r-1)d=rs$ and $d\ge2r\ne4$
there is a polynomial algorithm for checking almost $r$-embeddability of $s$-complexes in $\R^d$.
\end{theorem}

A version of Theorem \ref{t:algal} under the weaker (for $d\ge s+3$) dimension restriction  $rd\ge(r+1)s+3$
is \cite[Corollary 5]{MW16}, \cite[Theorem 1.6]{FV21}.
This version follows from Theorem \ref{t:mmw} below and the result \cite[Theorem 4]{MW16} of \cite{FV21} (see \cite[Remark 1.7]{Sk21d}).

It would be interesting to answer the following open questions.

$\bullet$ Which 2-complexes admit a PL map to $\R^3$ without triple points?

$\bullet$  Which 2-complexes are almost 3-embeddable in $\R^3$?

$\bullet$  Are there algorithms for checking the above properties of 2-complexes?

Same questions for $\R^3$ replaced by $\R^2$.

Any complex that is a triangulation of a $k(r-1)$-manifold admits an almost $r$-embedding to $\R^{kr}$ \cite{MW15}.
We conjecture that any contractible (or collapsible) $k(r-1)$-complex admits an almost $r$-embedding to $\R^{kr}$.

\smallskip
{\bf Toy examples.}
The following two exercises introduce the idea of a configuration space.
They are not formally used in the proof.

$\bullet$ Alice and Bob stand at the vertices $A,B$ of the triangle $ABC$.
They can walk continuously along the edges of the triangle in such a way that at each moment one of them is in a vertex, and the other is on the opposite edge.
Can they exchange their positions?

$\bullet$ Alice, Bob and Claude stand at the vertices $A,B,C$ of the tetrahedron $ABCD$.
They can walk continuously along the edges of the tetrahedron in such a way that at each moment two of them are at two different vertices, and the third is on the opposite edge.
Can they move to the respective vertices $D,B$, and $A$?

\begin{figure}[h]
\centerline{\includegraphics[width=3.5cm]{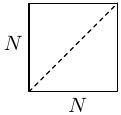}\qquad\qquad  \includegraphics[width=4cm]{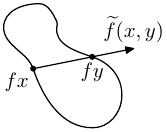}}
\caption{The deleted product and the Gauss map}\label{f:dp}
\end{figure}

\smallskip
{\bf The deleted product and Gauss map.}
The {\it deleted product} $\t K$ (another notation: $\t K^2$) of a set $K$ is the product of $K$ with itself, minus the diagonal:
$$
\t K:=\{(x,y)\in K\times K\ :\ x\ne y\},
$$
see fig. \ref{f:dp}. This is the configuration space of ordered pairs of distinct points in $K$.

Suppose that $f:K\to\R^d$ is an embedding of a subset $K\subset \R^m$.
Then the map $\t f:\t K\to S^{d-1}$ is well-defined by the Gauss formula
$$\t f(x,y)=\frac{f(x)-f(y)}{|f(x)-f(y)|}.$$
We have $\t f(y,x)=-\t f(x,y)$, that is, this map is equivariant with respect to the `exchanging factors' involution $(x,y)\mapsto(y,x)$ on $\t K$ and the antipodal involution on $S^{d-1}$.
Thus the existence of an equivariant map $\t K\to S^{d-1}$ is a necessary condition for the embeddability
of $K$ in $\R^d$.

Now assume that $K$ is a complex.
Then the {\it simplicial deleted product} of $K$ is
$$\cup \{ \sigma\times\tau\ : \sigma, \tau \textrm{ are simplices of }K,\ \sigma\cap\tau = \emptyset\}.$$
The existence of an equivariant map from this space to $S^{d-1}$ is necessary even for the almost 2-embeddability of $K$ in $\R^d$.
For results on the sufficiency of this condition see survey \cite[\S5, \S8]{Sk06}.
Those results are not used in this paper, but they are `predecessors' of the Mabillard-Wagner Theorems \ref{t:mawa} and \ref{t:mmw}.

\smallskip
{\bf The $r$-fold deleted product.}
The {\it $r$-fold deleted product} of a set $K$ is
$$
\t K^r:=\{(x_1,\ldots,x_r)\in K^r\ |\ x_i\ne x_j \mbox{ for every }i \neq j\}.
$$
This is the configuration space of ordered $r$-tuples of pairwise distinct points in $K$.
A nice notation for $\t K^r$ would be $K^{\underline r}$, because $n^{\underline r}:=n(n-1)\ldots(n-r+1)$
is the number of ordered $r$-tuples of pairwise different elements of an $n$-element set.

In this paper we shall use the notation $\t K^r$ for the following slightly different notion.
The {\it simplicial $r$-fold deleted product} of a complex $K$ is
\[
\t K^r:= \cup \{ \sigma_1 \times \cdots \times \sigma_r
\ : \sigma_i \textrm{ a simplex of }K,\ \sigma_i \cap \sigma_j = \emptyset \mbox{ for every }i \neq j \}.
\]
This is the union of products $\sigma_1\times\ldots\times\sigma_r$ formed by pairwise disjoint simplices of $K$.
Another notation: $K^{\times r}_{\Delta}$.


The set $\t K^r\subset K^r$ has no natural structure of a complex.
However, $\t K^r$ is the union of products of simplices with `nice' intersections, and therefore $\t K^r$ is a polyhedron.

\smallskip
{\bf Examples.}
Denote by $\t{\Delta_N}^3$ the set of ordered triples $(x,y,z)$ of points of $\Delta_N$ lying in pairwise disjoint faces.
Then:

(1) $\t{\Delta_1}^3=\emptyset$.

(2) The set $\t{\Delta_2}^3$ has 6 elements.

(3) The set $\t{\Delta_3}^3$ is a union of arcs forming a closed polygonal line.

(4) The set $\t{\Delta_4}^3$ is a union of triangles and squares forming a 2-dimensional polyhedron, and
one can prove that it is homeomorphic to the 2-sphere.


(5) The set $\t{\Delta_N}^3$ is a union of products of simplices, each product of dimension at most $N-2$.

\smallskip
{\bf Definitions of $\Sigma_r$ and $S^{d(r-1)-1}_{\Sigma_r}$.}
Denote by $\Sigma_r$ the permutation group of $r$ elements.
The group $\Sigma_r$ has a natural action on the set $\t K^r$, permuting the points in an $r$-tuple
$(p_1,\ldots, p_r)$.
This action is clearly free and PL, that is, compatible with some structure of a complex on $\t K^r$.

The group $\Sigma_r$ acts on the set of real $d\times r$-matrices by permuting the columns.
Denote by $S^{d(r-1)-1}_{\Sigma_r}$ the set formed by all those of such matrices, in which the sum of the elements in each row zero, and the sum of the squares of all the matrix elements is 1.
This set is homeomorphic to the sphere of dimension $d(r-1)-1$, and is invariant under the action of $\Sigma_r$.

The existence of an equivariant map $\t K^r\to S^{d(r-1)-1}_{\Sigma_r}$ is a necessary condition for the existence
of an almost $r$-embedding $K\to\R^d$.
The proof of this is similar to the above argument for $r=2$, and to the Configuration Space Lemma \ref{l:conf}.
A sufficiency result for this condition is Theorem \ref{t:mmw}.

\comment

Is it correct that a closed 2k-manifold  K  admits an immersion to  R^{3k}
without triple points if and only if  \chi(K)  is even?

Do you know any results and references on the following problem (for r>2):

Which n-manifolds  N  admits an immersion to  R^m  without r-tuple points?

    The answer to your first question is: No.
    1) The simplest contrexample is the connected sum of two copies of
    CP^2. This can not even be immersed to R^6.

    (Proof: If such an immersen exited , then the square of the normal
    Euler class would be the first normal Pontrjagin class, which is + or
    - 6 times the generator. But there is no such a class in the second
    cohomology group.)

    2) Actually thee following is true: If an oriented closed 4-manifold
    can be immersed into R^6 without triple points , then it must be
    null-cobordant (in oriented sense).

    (Proof: This follows from the Herbert formula. This formula says, that
    the triple points set represents the homology class dual to the square
    of the normal Euler class, i.e. the normal Pontrjagin class, So this
    must be zero.)

    Concerning your second question I do not know anything but the Herbert
    formula, that gives restriction for such manifolds

a counterexample is more simple (but less interesting).
If you take an immersion $f$ of an oriented manifold M^4 into R^6 with triple points, then the algebraic sum
T(f) is well defined as an integer. It is well know that there exists an immersion $f$ with A(f)\ne 0.
Then F=f \sharp f: M^4 \sharp M^4 \looparrowright R^6 is well-defined and A(F)=2A(f). Obviously, \chi(M \sharp M) is even. But there is no immersion of M \sharp M into R^6 without triple points: the algebraic number of triple points of an immersion depends no of the choice of  immersions. Petya.


 Denote $N:=(r-1)(3r+2)$.
Take an almost $r$-embedding $\Delta_N^{(3(r-1))} \to \R^{3r}$ of the $3(r-1)$-skeleton of the $N$-simplex.
Extend it arbitrarily to a map $f:\Delta_N\to \R^{3r}$.
Denote by $\rho(x)$ the distance from $x\in \Delta_N$ to the $3(r-1)$-skeleton.
It suffices to prove that $f\times\rho:\Delta_N\to \R^{3r+1}$ is an almost $r$-embedding.

Suppose to the contrary that $r$ points $x_1,\ldots, x_r\in\Delta_N$ lie in pairwise disjoint simplices
and are mapped to the same point under $f\times\rho$.
Dimension of one of those faces does not exceed $\frac{N+1}{r} - 1$, so it is at most $3(r-1)$.
W.l.o.g. this is the first face, hence $\rho(x_1)=0$.
Then $\rho(x_2) = \ldots = \rho(x_r)=\rho(x_1)= 0$, i.e. $x_1,\ldots, x_r\in\Delta_N^{(3(r-1))}$.
Therefore $f(x_1) = \ldots = f(x_r)$ contradicts the fact that $f|_{\Delta_N^{(3(r-1))}}$ is an almost $r$-embedding.





\endcomment

\section{Proof of Theorem \ref{t:tvepo}}\label{s:posi}

\subsection{The Borsuk-Ulam Theorem and its applications}\label{s:bu}

The arguments for Theorem \ref{t:tvepo} extend Borsuk-Ulam type results \cite{Ma03}.

\begin{theorem}[Borsuk-Ulam]\label{vf-bu}
For every $n$ and map $f:S^n\to\R^n$ there exists $x\in S^n$ such that $f(x)=f(-x)$.
\end{theorem}


For $n=2$ this means that at each moment of time there are two antipodal points on the Earth with the same  temperature and the same pressure.


For elementary proofs see \cite[p. 153-154]{Ma03}, \cite[\S8.8]{Pr06} and references therein.

This theorem has many equivalent formulations \cite{Ma03}.
The Borsuk-Ulam Theorem and its versions are famous due to their applications to combinatorics \cite{Ma03}
and to the proof of the existence of equilibria in mathematical economics (game theory)
\cite{AGL86, Ya99, SS03, SST95, SST02}.
The following well-known result gives some feeling for the applications.

\begin{theorem}[Ham Sandwich]\label{t:hasa} For any three convex polyhedra in 3-dimensional space there is a plane that splits their volumes into halves.\footnote{Usually this result is formulated as follows: in 3-dimensional space a sandwich, which is a slice of ham and two clices of bread (that is, any three measures) can all be simultaneously bisected with a single cut (that is, a plane).
We give an elementary formulation not involving measures.}
\end{theorem}

\begin{proof}
For each $x\in S^2$ there is a unique plane $\alpha(x)$ orthogonal to the vector~$x$ and~splitting the first polyhedron into two parts of equal volume.
Clearly, $\alpha(x)=\alpha(-x)$.
For $k=2,3$ this plane splits the $k$-th polyhedron into two (possibly empty) parts.
Denote by$\varphi_k(x)$ the difference between the volumes of the part in the half-space at which $x$ points, and the opposite part.
Clearly, $\varphi_2(-x)=-\varphi_2(x)$ and $\varphi_3(-x)=-\varphi_3(x)$.
By the Borsuk-Ulam Theorem~\ref{vf-bu} there is $x_0\in S^2$ such that $\varphi_2(x_0)=\varphi_2(-x_0)$ and
$\varphi_3(x_0)=\varphi_3(-x_0)$.
Then $\varphi_2(x_0)=\varphi_3(x_0)=0$.
So the plane $\alpha(x_0)$ is as required.
\end{proof}


\subsection{Proof of Theorem \ref{t:tvepo} for $r$ a prime}\label{s:posipr}

Our exposition is based on the well-structured and clearly written paper \cite{BSS}.
For a simpler alternative proof see \cite{VZ93}, \cite[\S6.4]{Ma03}, \cite[\S2.3.4]{Sk18}
(this proof uses deleted joins not deleted products, and so significantly simplifies the technical step of proving the Connectivity Lemma \ref{l:bss}).

Recall the definitions in \S\ref{s:posiex}.

\begin{lemma}[Configuration Space]\label{l:conf}
For each $r,N,d$ if there is an almost $r$-embedding $f\colon\Delta_N\to \R^d$,
then there is a $\Sigma_r$-equivariant map $\t{\Delta_N}^r\to S^{(r-1)d-1}_{\Sigma_r}$.
\end{lemma}

\begin{proof} (We present the argument for $r=3$. The case of an arbitrary $r$, not necessarily a prime, is proved similarly.)
The group $\Sigma_3$ acts freely on $(\R^d)^3-\diag$, where $\diag:=\{(x,x,x)\in(\R^d)^3\ :\ x\in\R^d\}$.
The required map is constructed as a composition
$$\t{\Delta_N}^3 \overset{f^3}\to(\R^d)^3-\diag \overset{\pi}\to S^{2d-1}_{\Sigma_3}.$$
Here:

$\bullet$ the equivariant map $f^3$ is well-defined because of the condition $f(\sigma_1)\cap f(\sigma_2) \cap f(\sigma_3)=\emptyset$.

$\bullet$ informally, $\pi$ is the projection to the subspace orthogonal to the diagonal; formally, for
$x_1,x_2,x_3\in\R^d$ not all equal we define
$$
\pi'(x_1,x_2,x_3)\ :=\ \left(\ 2x_1-x_2-x_3,\ 2x_2-x_1-x_3,\ 2x_3-x_1-x_2\ \right)
\quad\text{and}\quad \pi:=\frac{\pi'}{|\pi'|}.
$$
\end{proof}

\begin{theorem}[$r$-fold Borsuk-Ulam Theorem]\label{t:bu3}
Let $X$ be a $k$-connected $(k+1)$-complex, $r$ a prime number, and $\omega_X:X\to X$ and
$\omega_S:S^k\to S^k$ simplicial maps without fixed points such that $\omega_X^r=\id X$ and $\omega^r_S=\id S^k$.
Then there are no maps $f:X\to S^k$ such that $f\circ\omega_X=\omega_S\circ f$.
\end{theorem}

\begin{proof}[Comments on the proof]  Since $X$ is $k$-connected and $\omega_X$ is simplicial,
there is a map $g:S^k\to X$ such that $g\circ\omega_S=\omega_X\circ g$.
Then the composition $f\circ g:S^k\to S^k$ extend over $D^{k+1}$ and commutes with $\omega_S$.
The non-existence of such map is in turn an analogue of the Borsuk-Ulam Theorem.
It is proved analogously to the Borsuk-Ulam Theorem \cite[Lemma 2]{BSS}.
In \cite{Ma03} an alternative proof (using {\it Lefschetz trace formula}) and a
further generalization ({\it Dold Theorem}) are presented.
\end{proof}

\begin{lemma}[Connectivity; {\cite[Lemma 1]{BSS}}]\label{l:bss}
For each $r,N$ the polyhedron $\t{\Delta_N}^r$ is $(N-r)$-connected.
\end{lemma}

\begin{proof}[Comments on the proof]
{\it Hurewicz theorem.}
Let $Y$ be a simply connected complex such that $H_j(Y;\Z)=0$ for each $j\in\{2,3,\ldots,k\}$.

(a) Then any map $f:S^k\to Y$ extends over $D^{k+1}$.

(b) If there are fixed point free maps $\omega:Y\to Y$ and $\omega:S^k\to S^k$ such that $\omega^3=\id Y$ and
$\omega^3=\id S^k$, then there is a map $f:S^k\to Y$ such that $f\circ\omega=\omega\circ f$.

(Part (b) is not a part of the usual Hurewicz Theorem but is easily deduced from (a).)

By the Hurewicz Theorem it suffices to prove that $\t{\Delta_N}^r$ is
1-connected and $H_j(\t{\Delta_N}^r;\Z)=0$ for each $j\in\{2,3,\ldots,N-r\}$.
This is proved in \cite[Proof of Lemma 1]{BSS}, cf. Examples in \S\ref{s:posiex}.

See alternative proof in \cite{BZ16}.
\end{proof}


\begin{proof}[Proof of Theorem \ref{t:tvepo} for a prime $r$]
Let $N:=(r-1)(d+1)$.
Then $N-r=(r-1)d-1$.
Suppose to the contrary that there is an almost $r$-embedding $\Delta_N\to \R^d$.
Take a $\Sigma_r$-equivariant map $\t{\Delta_N}^r\to S^{(r-1)d-1}_{\Sigma_r}$ given by the Configuration Space Lemma \ref{l:conf}.
The cyclic-shift-of-$r$-tuples-by-one self-map $\omega_{\t{\Delta_N}^r}$ of $\t{\Delta_N}^r$ has no fixed points.
Since $r$ is a prime, the cyclic-shift-of-columns-by-one self-map $\omega_S$ of $S^{(r-1)d-1}_{\Sigma_r}$ has no  fixed points.
Therefore by the $r$-fold Borsuk-Ulam Theorem \ref{t:bu3} for $k=N-r$ and $X=\t{\Delta_N}^r$ we obtain a contradiction to the Connectivity Lemma \ref{l:bss}, because $X$ is clearly $(N-r+1)$-dimensional.
\end{proof}

\subsection{Proof of Theorem \ref{t:tvepo} for $r$ a prime power}\label{s:posipow}

Our exposition is based on \cite[\S2]{vo96} and \cite[\S2,\S3]{Oz};
see alternative proofs in \cite{Vo96v, Zi98, Sa00}.

An action is said to be {\it fixed point free} if the space has no points fixed by {\it each} element of the group.

\begin{theorem}[$r$-fold Borsuk-Ulam Theorem; {\cite[Lemma]{vo96}, \cite[Lemma 3.3]{Oz}}]\label{t:bu4}
Let $\alpha>0$ be an integer and $X$ a $k$-connected $(k+1)$-complex with a simplicial action of $\Z_p^\alpha$.
For each fixed point free action of $\Z_p^\alpha$ on $S^k$ there are no $\Z_p^\alpha$-equivariant maps $X\to S^k$.
\end{theorem}

This is deduced below from Lemma \ref{l:sp} and Localization Theorem \ref{t:loc}.

\begin{proof}[Proof of Theorem \ref{t:tvepo}]
Denote by $G$ the subgroup of $\Sigma_r$ formed by all permutations preserving $\{kp^s+1,kp^s+2,\ldots,(k+1)p^s\}$
for each $s=1,2,\ldots,\alpha-1$ and $k=0,1,\ldots,p^{\alpha-s}-1$.
Clearly, $G\cong\Z_p^\alpha$.
Let $N:=(d+1)(r-1)$.
The theorem now follows (as in the case of prime $r$) by the Configuration Space Lemma \ref{l:conf},
the $r$-fold Borsuk-Ulam Theorem \ref{t:bu4} for $k=N-r$ and $X=\t{\Delta_N}^r$, together with the Connectivity Lemma \ref{l:bss}, because $X$ is clearly $(N-r+1)$-dimensional.
\end{proof}

The following Lemma holds for any coefficients, but is used for the coefficients $\Z_p$.

\begin{lemma}\label{l:sp} (a) If $X$ is a $k$-connected $(k+1)$-complex and $X\to E\to B$ is a bundle,
then $p^*:H^{k+1}(B)\to H^{k+1}(E)$ is injective.

(b) If $S^k\to E\to B$ is a bundle over a connected base $B$, \ $\pi_1(B)$ acts trivially on $H^k(S^k)$,
and  $p^*_{k+1}:H^{k+1}(B)\to H^{k+1}(E)$ is injective, then $p^*_j:H^j(B)\to H^j(E)$ is injective for each $j$.
\end{lemma}

\begin{proof} The Lemma is a standard exercise in spectral or Gysin sequences.
The reader can omit the proof, and just use the statement of the lemma.

In order to prove (a) consider the Serre spectral sequence \cite[\S21 in the Russian version]{FF89}
for which $E_2^{i,j}\cong H^i(B,H^j(X)_T)$.
Since $H^j(X)=0$ unless $j\in\{0,k+1\}$, there are no non-trivial differentials
$d_s^{i,j}:E_s^{i,j}\to E_s^{i+s,j-s+1}$ to or from
$$H^{k+1}(B)\cong E_2^{k+1,0}=E_3^{k+1,0}=\ldots=E_\infty^{k+1,0}.$$
Thus, $p^*$ is the composition of these isomorphisms and the inclusion $E_\infty^{k+1,0}\to H^{k+1}(E)$.
Hence $p^*$ is injective.

In order to prove (b) look at the following segment of the Gysin exact sequence \cite[\S22.4.A in the Russian version]{FF89}:
$$H^{j-k-1}(B)\overset{d_j}\to H^j(B)\overset{p^*_j}\to H^j(E),$$
where $d_j(x):=x\cup e$ for the characteristic class $e\in H^{k+1}(B)$ of the bundle.
Since $p^*_{k+1}$ is injective, $d_{k+1}=0$, so $e=0$.
Hence, for each $j$ we have that $d_j=0$ and $p^*_j$ is injective.
\end{proof}

For complexes $A,X$ with simplicial actions of a group $G$ define an action of $G$ on $A\times X$ by $g(a,x):=(g(a),g(x))$.
Define
$$\overline X_A:=(A\times X)/G\quad\text{and define}\quad
\pi_X=\pi_{X,A}:\overline X_A\to A/G \quad\text{by}\quad \pi_X[a,x]:=[a].$$

\begin{theorem}\label{t:loc} \cite[Proposition III.3.14]{Di}
For each prime $p$, \ $\alpha\ge1$ and $G:=\Z_p^\alpha$ there is an infinite
complex $E=E_G$ with a free action of $G$ such that

$\bullet$ for each complex $X$ with an action of $G$ the map $\pi_X=\pi_{X,E}$ is a bundle;

$\bullet$ for each fixed point free action of $G$ on $S^k$ and some $j>0$ the map
$\pi^*_{S^k}=\pi^*_{S^k,E}:H^j(E/G;\Z_p)\to H^j(\overline{S^k}_E;\Z_p)$ is not injective.
\end{theorem}

\begin{proof}[Proof of Theorem \ref{t:bu4}]
Let $G:=\Z_p^\alpha$.
Assume to the contrary that a $G$-equivariant map $f:X\to S^k$ exists.
Take a complex $E$ as given by Theorem \ref{t:loc}.
Consider the diagram
$$\xymatrix{ \overline X_E \ar[rr]^{\overline f=(\id E\times f)/G} \ar[dr]_{\pi_X} & & \overline{S^k}_E \ar[dl]^{\pi_{S^k}} \\
 & E/G
}$$
Here $\overline f$ is well-defined by $\overline f[a,x]:=[a,f(x)]$.
Apply Lemma \ref{l:sp}.a to the bundle $\pi_X$.
We obtain that $\pi_X$ induces an injection on $H^{k+1}(\cdot;\Z_p)$.
By the commutativity the map $\pi_{S^k}$ also induces an injection on $H^{k+1}(\cdot;\Z_p)$.
A self-homeomorphism of $S^k$ of order $p$ induces an isomorphism of $H^k(S^k;\Z_p)$.
Hence $\pi_1(E/G)$ acts trivially on $H^k(S^k;\Z_p)$.
Apply Lemma \ref{l:sp}.b to the bundle $\pi_{S^k}$.
We obtain a contradiction to Theorem \ref{t:loc}.
\end{proof}



\section{Proof of Theorem \ref{t:vkfr-}}\label{s:coex}

\subsection{Plan of the proof of Theorem \ref{t:vkfr-}}\label{s:plan}

We recall from \S\ref{s:state} that

$\bullet$ Theorem \ref{t:tve} for $d=3r+1$ (counterexample to the topological Tverberg Conjecture)
follows from Theorem \ref{t:vkfr-} (counterexample to the $r$-fold van Kampen-Flores Conjecture) and the Constraint Lemma \ref{p:redu} (which states that the topological Tverberg Conjecture implies the $r$-fold van Kampen-Flores Conjecture);

$\bullet$ the Constraint Lemma \ref{p:redu} is proved by Gromov and later rediscovered by Blagoevic-Frick-Ziegler;

$\bullet$ Theorem \ref{t:vkfr-} follows from Theorems \ref{t:mawa} and \ref{t:ozmawa}, which are due to Mabillard-Wagner and Ozaydin, see references below in this section
(observe our proof of Theorem \ref{t:vkfr-} is different from \cite{Oz, MW15}).

See more in Remark \ref{r:hystor}.

\begin{figure}[h]
\centerline{\includegraphics[width=5cm]{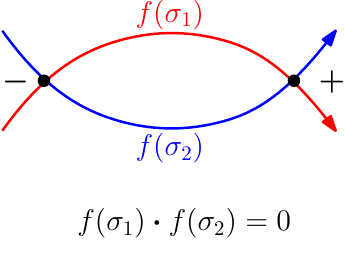}\qquad \qquad \includegraphics[width=8cm]{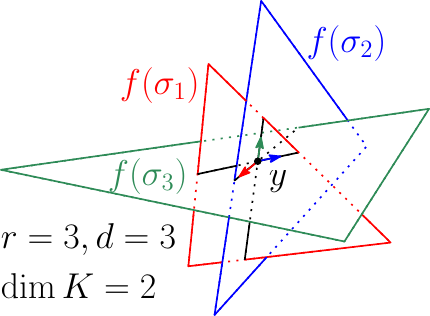}}
\caption{$r$-fold points and their $r$-intersection signs for $r=2,3$}\label{f:gl2}
\end{figure}

For statement and proof of Theorems \ref{t:mawa} and \ref{t:ozmawa} we need the following definitions.


In this section we shorten $f(x)$ to $fx$.

Let $K$ be a $k(r-1)$-complex for some $k\ge1$, $r\ge2$, and $f\colon K\to\R^{kr}$ a PL map in general position.

Then preimages $y_1,\ldots,y_r\in K$ of any $r$-fold point $y\in\R^{kr}$ (that is, of a point having
$r$ preimages) lie in the interiors of $k(r-1)$-dimensional simplices of $K$.
Choose arbitrarily an orientation for each of the $k(r-1)$-simplices.
By general position, $f$ is affine on some neighborhood $U_j$ of $y_j$ for each $j=1,\ldots,r$.
Take a positive basis of $k$ vectors in the oriented normal space to oriented $fU_j$.
{\bf The $r$-intersection sign} of $y$ is the sign $\pm 1$ of the basis in $\R^{kr}$ formed by $r$ such $k$-bases.
(This construction is classical for $r=2$ \cite{BE82} and is analogous for $r\ge3$, cf. \cite[\S~2.2]{MW15}.)

We call the map $f$ a {\bf $\Z$-almost $r$-embedding} if $f\sigma_1\iprod\ldots\iprod f\sigma_r=0$ whenever $\sigma_1,\ldots,\sigma_r$ are pairwise disjoint simplices of $K$.
Here the {\it algebraic $r$-intersection number}  $f\sigma_1\iprod\ldots\iprod f\sigma_r\in\Z$ is defined as the sum of the $r$-intersection signs of all $r$-fold points $y\in f\sigma_1\cap\ldots\cap f\sigma_r$.
The sign of the number $f\sigma_1\iprod\ldots\iprod f\sigma_r$ depends on an arbitrary choice of orientations for each $\sigma_i$, and on the order of $\sigma_1,\ldots,\sigma_r$, but the condition $f\sigma_1\iprod\ldots\iprod f\sigma_r=0$ does not.
See fig. \ref{f:gl2}.

Clearly, an almost $r$-embedding is a $\Z$-almost $r$-embedding.

\begin{theorem}[Mabillard-Wagner]\label{t:mawa}
For every integers $k\ge3$ and $r>0$ if there is a $\Z$-almost $r$-embedding of a $k(r-1)$-complex $K$ in $\R^{kr}$, then there is an almost $r$-embedding of $K$ in $\R^{kr}$.
\end{theorem}

This result follows from \cite[Theorem 7]{MW15}; in \S\ref{s:mw} we outline a simpler proof \cite[\S2]{AMSW}
which can be generalized to a proof of the codimension $k=2$ version \cite{AMSW}.
See Remark \ref{r:mw} for background.

\begin{theorem}\label{t:ozmawa} If $r$ is not a prime power and $k>0$ is an integer,
then there is a $\Z$-almost $r$-embedding of any  $k(r-1)$-complex in $\R^{kr}$.
\end{theorem}

This result is implied by the following \"Ozaydin Theorem \ref{theorem:ozaydin-product}
for $d=kr$, and Proposition \ref{cor:equiv-alm} of Mabillard-Wagner.
Recall the definitions in \S\ref{s:posiex}.
For a $k(r-1)$-complex $K$ we have $\dim\t K^r\le r\dim K=kr(r-1)$.

\begin{theorem}[\"Ozaydin]\label{theorem:ozaydin-product}
If $r$ is not a prime power, $d>0$ is an integer and $X$ is a $d(r-1)$-complex with a free PL action of $\Sigma_r$, then there is a $\Sigma_r$-equivariant map $X\to S^{d(r-1)-1}_{\Sigma_r}$.
\end{theorem}

In \S\ref{s:oz} we present a clearer proof than the original (see footnote \ref{f:ozaydin}).

\begin{proposition}\label{cor:equiv-alm}
Let $k,r\ge2$ be integers and $K$ a $k(r-1)$-complex.
There is a $\Z$-almost $r$-embedding of $K$ in $\R^{kr}$ if and only if there is a $\Sigma_r$-equivariant map
$\t K^r\to S^{kr(r-1)-1}_{\Sigma_r}$.
\end{proposition}

This is a generalization of known results, see \cite[\S~4.1 and Corollary 44]{MW15} and \S\ref{s:fm}.

The following extension of the Mabillard-Wagner Theorem \ref{t:mawa} is not used in the proofs of Theorems \ref{t:tve} and \ref{t:vkfr-}.
However, it is both interesting in itself, and is useful for stronger counterexamples (\S\ref{s:intrm}, \cite{AKS}).


\begin{theorem}[\cite{MW15, AMSW, MW16, MW16', Sk17, Sk17o}]\label{t:mmw}
Suppose that $K$ is a $s$-complex and either $rd\ge(r+1)s+3$ or $d=2r=s+2\ne4$.
There exists an almost $r$-embedding $f:K\to\R^d$ if and only if there exists a $\Sigma_r$-equivariant map $\widetilde K^r \to S^{d(r-1)-1}_{\Sigma_r}$.
\end{theorem}

The case $(r-1)d=rs$ of Theorem \ref{t:mmw} was proved in \cite{MW15, AMSW}; this case is covered by
the Mabillard-Wagner Theorem \ref{t:mawa}, its generalization from \cite{AMSW} and Proposition \ref{cor:equiv-alm}.


The remaining subsections of this section are independent on each other.

\begin{remark}[On the Mabillard-Wagner Theorem \ref{t:mawa}]\label{r:mw}
The Mabillard-Wagner Theorem \ref{t:mawa} is the most non-trivial part of the disproof of the topological Tverberg conjecture.
Their idea is similar to, but different from, `Haefliger's $h$-principle for embeddings'
\cite[2.1.1, (E), p. 50-51]{Gr86}, \cite[\S5]{Sk06} (for `$h$-principle' itself see \cite[p. 3]{Gr86})
and Whitney trick \cite[Whitney Lemma 5.12]{RS72}.
The {\it $r$-fold} analogues of Haefliger's $h$-principle for embeddings and of Whitney trick were `in the air' since 1960s \cite[\S5.6 `The Generalized Haefliger-Wu invariant']{Sk06}.
Cf. Remark \ref{r:ae}.d.
`Positive results' were available for links, and an argument involving triple Whitney trick was sketched
by S. Melikhov in 2007 (but published only in \cite{Me17, Me18}; the authors of \cite{MW15} received early private version of \cite{Me17} before arxiv submission, but do not refer to it in \cite{MW15}).
The open problem cited in footnote \ref{f:gropro} suggests that Gromov was aware that Theorem \ref{t:vkfr-} might hold, to the extent of asking the right question.
However, some counterexamples were known \cite[\S5.6]{Sk06}.
So $r$-fold analogue of Haefliger's $h$-principle for {\it almost} $r$-embeddings and of the Whitney trick is an important contribution of Mabillard and Wagner.

Their $r$-fold Whitney trick involves an analogue of increasing the connectivity of the intersection set by surgery \cite{Ha63}, \cite[Theorem 4.5 and appendix A]{HK98}, \cite[Theorem 4.7 and appendix]{CRS}.
In other words, this is first attaching an embedded 1-handle along an arc (`piping') and then attaching a canceling  embedded 2-handle along a disk
(`unpiping') \cite[\S3]{Ha62}, \cite[proof of Theorem 1.1 in p. 7]{Me17}.
(We apply \cite[Proposition 3.3]{Ha62} for $r=0$ and $r=1$; both times we pass from embedding into $B\times0$ to  embedding into $B\times1$.)
Application of these constructions is non-trivial and is an important achievement of Mabillard and Wagner.
\end{remark}

\subsection{Proof of the \"Ozaydin Theorem \ref{theorem:ozaydin-product}}\label{s:oz}


The idea is to deduce the \"Ozaydin Theorem \ref{theorem:ozaydin-product} from its local versions `away from a prime $p$' for each $p$.
The property that $r$ is not a prime power will be used by using that any group containing a prime power
$r$ of elements cannot act transitively on $\{1,2,\ldots,r\}$.

We use the following well-known results of the group theory and the equivariant obstruction theory (Lemmas \ref{lemma:m} and \ref{lemma:o}).
For the reader's convenience we sketch their proofs below;
however, the reader can just use the lemmas without proof.

Denote the order of $p$ in the factorization of $r! = |\Sigma_r|$ by
$$
\alpha_p=\alpha_{p,r} = \sum_{k=1}^\infty \left\lfloor \frac{r}{p^k}\right\rfloor.
$$

\begin{lemma}[a particular case of the Sylow theorem]\label{lemma:m}
For any integer $r$ and prime $p$ there is a subgroup $G$ of $\Sigma_r$ having $p^{\alpha_p}$ elements.
\end{lemma}

{\bf Remark.}
The simplest counterexample to the topological Tverberg conjecture is obtained for $r=6$, when
the subgroups in Lemma \ref{lemma:m} have a very simple description.
A group $G$ acting on a set $X$ \emph{preserves the splitting $X=X_1\sqcup \ldots \sqcup X_k$} if for any $g\in G$ and any $i=1,\ldots,k$ we have $g(X_i) = X_j$ for some $j$.

$\bullet$ $p=5$, \ $\Z_5\cong G < \Sigma_5 < \Sigma_6$, where $G$ is the subgroup preserving the splitting
$\{1,2,3,4,5,6\} = \{1,2,3,4,5\}\sqcup \{6\}$ and acting in a cyclic way on $\{1,2,3,4,5\}$;

$\bullet$ $p=3$, \ $\Z_3\times \Z_3\cong G < \Sigma_3 \times \Sigma_3 < \Sigma_6$, where $G$ is the subgroup preserving the splitting $\{1,2,3,4,5,6\} = \{1,2,3\} \sqcup \{4,5,6\}$ and acting in a cyclic way on both $\{1,2,3\}$ and $\{4,5,6\}$;

$\bullet$ $p=2$, \ $G$ is the subgroup preserving the splittings
$\{1,2,3,4,5,6\} = \{1,2,3,4\} \sqcup \{5,6\}$  and  $\{1,2,3,4\} = \{1,2\} \sqcup \{3,4\}$.

\begin{lemma}[Obstruction]\label{lemma:o}
Let $r>0$ be an integer and $X$ a $d(r-1)$-complex with a free PL action of $\Sigma_r$.
Then to each subgroup $G$ of $\Sigma_r$ there corresponds
an Abelian group $H_G$ and an element $v(G)\in H_G$ so that:

(*) $v(G)=0$ $\Leftrightarrow$ there is a $G$-equivariant map $X\to S^{d(r-1)-1}_{\Sigma_r}$;

(**) there is a homomorphism $\tau:H_G\to H_{\Sigma_r}$ with $\tau v(G)=[\Sigma_r:G]v(\Sigma_r)$.
\end{lemma}

\begin{proof}[Proof of the \"Ozaydin Theorem~\ref{theorem:ozaydin-product}]
If $\dim X<d(r-1)$, then the existence of an equivariant map $X\to S^{d(r-1)-1}_{\Sigma_r}$ follows because $S^{d(r-1)-1}_{\Sigma_r}$ is $(d(r-1)-2)$-connected.

Now assume that $\dim X=d(r-1)$.

Take any prime $p$.
Take a subgroup $G$ of $\Sigma_r$ as given by Lemma \ref{lemma:m} (a $p$-Sylow subgroup).
Since $r$ is not a power of $p$, the group $G$ cannot act transitively on $\{1,2,\ldots,r\}$,
otherwise this would be a coset of the $p$-group $G$, and would have an order equal to a power of $p$.
Hence we can assume that $G$ preserves the sets $\{1,\ldots,k\}$ and $\{k+1,\ldots, r\}$ for some $k$.

This gives a $G$-invariant $(d\times r)$-matrix $M$ of $d$ rows
$$
(\underbrace{k-r,\ldots,k-r}_{k},\underbrace{k,\ldots,k}_{r-k}).
$$
Then the point $M/|M|\in S^{d(r-1)-1}_{\Sigma_r}$ is $G$-invariant.
Thus a $G$-equivariant map $X\to S^{d(r-1)-1}_{\Sigma_r}$ is defined by mapping $X$ to the $G$-invariant point.
Then by the Obstruction Lemma \ref{lemma:o}, the `$\Leftarrow$' direction of (*) and (**),
$\dfrac{r!}{p^{\alpha_p}}v(\Sigma_r)=\tau v(G)=0$.

The numbers $r!/p^{\alpha_p}$ for all primes $p<r$ have no common multiple.
Therefore $v(\Sigma_r)=0$.
The theorem now follows by the Obstruction Lemma \ref{lemma:o}, the `$\Rightarrow$' direction of (*).
\end{proof}

{\bf Remark.} For each even $r$ there is an elementary proof that $v(\Sigma_r)=-v(\Sigma_r)$, i.e. $2v(\Sigma_r)=0$.
If $r=p^\alpha$, the above proof gives the relation $p^{\alpha(\alpha+1)/2}v(\Sigma_r)=0$.
For an elementary construction of $v(\Sigma_r)$, which could potentially give an elementary proof of
the \"Ozaydin Theorem~\ref{theorem:ozaydin-product}, see \cite{Sk18}.


\begin{proof}[Proof of Lemma \ref{lemma:m}]
{\it Construction of a tree $T$ whose leaf vertices are numbered by $0,1,\ldots,r-1$.}
Take an integer $\ell$ such that $p^\ell< r < p^{\ell+1}$.
Denote by  $T_0$ the graded (i.e. levelled) tree whose vertices are words in the alphabet $0,1,\ldots,p-1$ having at most $\ell+1$ letters, and the children of a vertex $w$ are those words that can be obtained from $w$ by adding a letter to the right.
The level of a vertex is the length of the corresponding word.
The word $a_1a_2a_3\ldots a_{\ell+1}$ is base-$p$ expansion of the number
$$\overline{a_{\ell+1}\ldots a_3a_2a_1}=a_1+a_2p+a_3p^2+\ldots+a_{\ell+1} p^\ell\quad\text{(possibly $a_{\ell+1}=0$)}.$$
Denote by $T$ the tree obtained from $T_0$ by deleting all the vertices that have no descendant whose number
is less than $r$ (in particular, deleting all the vertices whose numbers are greater or equal to $r$).
This gives a strictly smaller tree because $r < p^{\ell+1}$.
The leaf vertices of $T$ have level $\ell+1$, that is, they are all $(\ell+1)$-letter words.

{\it Construction of $G$.}
If there are $p$ children of a given vertex $v\in T$, then order them cyclically as residues modulo $p$.
If there are less than $p$ children of $v$, then order them linearly.
Both orderings are defined in terms of the last letter.
Let $G$ be the group of automorphisms of $T$ that preserve the level and the described (cyclic or linear) order of children of every vertex.

{\it Proof that $|G|=p^{\alpha_p}$.}
First, we rotate independently collections of vertices of level $\ell+1$ (that is, leaf vertices)
that have a cyclic order.
In total, we obtain $p^{\lfloor r/p\rfloor}$ permutations.
Second, we rotate independently collections of vertices of level $\ell$ that have a cyclic order.
In total, we obtain $p^{\lfloor r/p^2\rfloor}$ permutations, and so on.
Every element of $G$ can be build this way by rotating the children of the vertices, starting from the bottom and going to the top of the tree.
This shows that $|G|= p^{\alpha_p}$.
\end{proof}

\begin{proof}[Sketch of a proof of the Obstruction Lemma \ref{lemma:o}]
Let $H_G:=H_G^{(r-1)d}(X; \Z)$ be the {\it $G$-equivariant (=$G$-symmetric) cohomology group}
w.r.t. the action of $G$ on $X$ and the multiplication-by-the-sign (of the permutation) action of $G$ on $\Z$.
This group and the element $v(G)\in H_G$ are defined in a natural way by trying to construct a $G$-equivariant map $X\to S^{d(r-1)-1}_{\Sigma_r}$ by the skeleta of some $G$-invariant subdivision of $X$.
For the details see~\cite{hu59, Sk15, Sk} in the non-equivariant setting, and \cite{Di} in the equivariant setting.\footnote{Alternatively, $H_G=H^{(r-1)d}((X\times E_G)/G;\Z)$, where $E_G$ is the complex in the Localization Theorem \ref{t:loc}.}
Since $\dim X=(r-1)d$, we obtain (*).

Let $\forg:H_{\Sigma_r}\to H_G$ be the `forgetting symmetry' homomorphism.
Clearly, $v(G)=\forg v(\Sigma_r)$.

The property (**) follows because there is a homomorphism $\tau:H_G\to H_{\Sigma_r}$ such that
$\tau\circ\forg$ is the multiplication by $[\Sigma_r:G]$.
In order to define such a map $\tau$ (a {\it transfer} homomorphism) we define a map
$t:C_G^{(r-1)d}(X; \Z)\to C_{\Sigma_r}^{(r-1)d}(X; \Z)$ of the {\it simplicial cochain groups}.
Let $s:=[\Sigma_r:G]$ and take $f_1,\ldots,f_s\in \Sigma_r$ such that
$\Sigma_r=f_1G\sqcup\ldots\sqcup f_sG$.
For a simplex $\sigma$ of the $G$-invariant subdivision of $X$ we define
$$t(x)(\sigma)=x(f_1\sigma)\sgn f_1+\ldots+x(f_s\sigma)\sgn f_s.$$
We extend $t(x)$ to a cochain by linearity.
It is easy to check that $t$ defines the required map $\tau$ in cohomology, see~\cite{brown1982}, \cite[\S5.2]{BLZ} for the details.
\end{proof}

\subsection{Proof of Proposition \ref{cor:equiv-alm}}\label{s:fm}


Let $d:=kr$ and take any:

$\bullet$ orientations on $k(r-1)$-simplices of $K$;

$\bullet$ general position PL map $f\colon K\to \R^d$;

$\bullet$ collection $\sigma_1,\ldots,\sigma_r$ of pairwise disjoint $k(r-1)$-simplices of $K$.


Compare the following with proof of the Configuration Space Lemma \ref{l:conf} in \S\ref{s:posi}.
For $x_1,\ldots,x_r\in\R^d$ which are not all equal we define
$$
S:=x_1+\ldots+x_r,\quad \pi':=\left(x_1-\frac Sr,\ldots,x_r-\frac Sr\right)
\quad\text{and}\quad \pi:=\frac{\pi'}{|\pi'|}.
$$
This defines a map
$$
\pi:(\R^d)^r-\diag\to S^{d(r-1)-1}_{\Sigma_r},\quad\text{where}\quad
\diag:=\{(x,x,\ldots,x)\in(\R^d)^r\ |\ x\in\R^d\}.
$$
By general position $f^r\partial(\sigma_1\times \ldots \times \sigma_r)\subset (\R^d)^r-\diag$, so we obtain the map
$$
\pi\circ f^r:\partial(\sigma_1\times \ldots \times \sigma_r)\to S^{d(r-1)-1}_{\Sigma_r}.
$$
By \cite[Lemma 41.b]{MW15}
$f\sigma_1\iprod\ldots\iprod f\sigma_r=\pm\deg(\pi\circ f^r|_{\partial(\sigma_1\times \ldots \times \sigma_r)})$.
Thus the map $\pi\circ f^r$ extends to $\sigma_1\times \ldots \times \sigma_r$ if and only if  $f\sigma_1\iprod\ldots\iprod f\sigma_r=0$.


\smallskip
{\it Proof of the `only if' part of Proposition \ref{cor:equiv-alm}.}
(This part is not required for Theorem \ref{t:ozmawa}.)
Since $f\sigma_1\iprod\ldots\iprod f\sigma_r=0$ for any
pairwise disjoint $k(r-1)$-simplices $\sigma_1,\ldots,\sigma_r$ of $K$,
we obtain that the $\Sigma_r$-equivariant map  $\pi\circ f^r$ defined on the codimension 1 skeleton of $\t K^r$ extends to a $\Sigma_r$-equivariant map $\t K^r\to S^{d(r-1)-1}_{\Sigma_r}$.

\begin{figure}[h]
\centerline{\includegraphics[width=3.5cm]{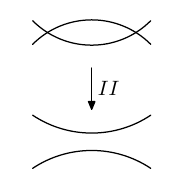}\qquad\includegraphics[width=2.5cm]{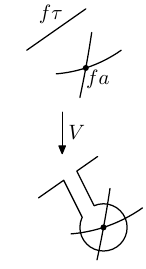} }
\caption{For map of a graph to the plane: II --- Whitney trick, V ---Van Kampen finger move}
\label{f:reidall}
\end{figure}

\smallskip
{\it Sketch of the proof of the `if' part of Proposition \ref{cor:equiv-alm}.}
This is the most non-trivial part of Proposition \ref{cor:equiv-alm},
so it could be omitted for the first reading.

The existence of a $\Sigma_r$-equivariant map $\t K^r\to S^{d(r-1)-1}_{\Sigma_r}$ implies that the $r$-dimensional `partial matrix' with entries $f\sigma_1\iprod\ldots\iprod f\sigma_r$ for pairwise disjoint $k(r-1)$-simplices $\sigma_1,\ldots,\sigma_r$ of $K$, is a sum of {\it elementary coboundaries} \cite{hu59, Sk15, Sk, Di}.
Now we can obtain from $f$ a $\Z$-almost $r$-embedding $K\to\R^d$ using a higher-multilplicty generalization \cite[Corollary 44]{MW15} of {\it van Kampen finger moves} corresponding to  elementary coboundaries \cite{Fo04},  \cite[\S1.5.3]{Sk18}, \cite[\S5.10, \S5.11]{Sk} (cf. fig. \ref{f:reidall}.V).

\subsection{Proof of the Mabillard-Wagner Theorem \ref{t:mawa}}\label{s:mw}

In this subsection we assume that $k\ge3$ and $r\ge2$ are integers.
In this paper we work in the PL category, in particular, all disks, balls and maps are PL.
Let $B^d:=[0,1]^d$ denote the standard ball and $S^{d-1}=\partial B^d$ the standard sphere.
We need to speak about balls of different dimensions and we will use the word `disk' for lower-dimensional
objects and `ball' for higher-dimensional ones in order to clarify the distinction (even though, formally, the disk $D^d$ is the same as the ball $B^d$).
We denote by $\partial M$, respectively $\Int M$, the boundary, respectively the interior, of a manifold $M$.
A map $f\colon M \to B^d$ from a manifold with boundary to a ball is called {\bf proper},
if $f^{-1}S^{d-1}=\partial M$.

The following result is interesting in itself, and in simple terms illuminates the core of the proof of Mabillard-Wagner Theorem~\ref{t:mawa}.

\begin{theorem}[Local Disjunction]\label{l:ld+3}
Let $D=D_1 \sqcup \ldots \sqcup D_r$ be the disjoint union of $k(r-1)$-dimensional  disks and $f:D\to D^{kr}$ a proper general position PL map to the $kr$-ball such that $fD_1\iprod\ldots\iprod fD_r=0$.
(Here $fD_1\iprod\ldots\iprod fD_r$ is the sum of the $r$-intersection signs of all
$r$-fold points
$y \in fD_1\cap\ldots\cap fD_r$.)

Then there exists a proper general position PL map $g:D\to D^{kr}$ such that $g=f$ on $\partial D$ and $gD_1\cap\ldots\cap gD_r=\emptyset$.
\end{theorem}

\begin{figure}[h]
\centerline{\includegraphics[width=7cm]{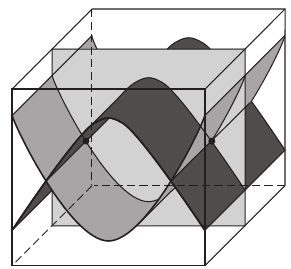} \qquad \qquad \includegraphics[width=7cm]{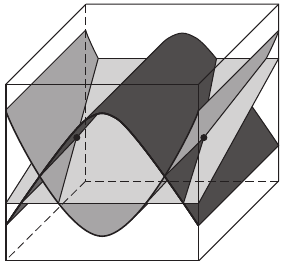}}
\caption{Whitney trick for $r=3$}\label{f:wh3}
\end{figure}

The case $r=2$ is the well known {\it double} Whitney trick \cite{RS72} (cf. fig. \ref{f:reidall}.II).
The main idea of proof for $r>2$ \cite{MW15}, \cite[\S2]{AMSW} is to invent an $r$-tuple analogue of the Whitney trick, that is, an analogue for higher-multiplicity intersections.
To see why this analogue is non-trivial, take $r=3$.
If two triple points of opposite signs in $fD_1\cap fD_2\cap fD_3$ are contained in one connected component of $fD_1\cap fD_2$, then we can `cancel' them by the double Whitney trick applied to $fD_1\cap fD_2$ and $fD_3$
(see fig. \ref{f:wh3} left, where $fD_1$ is the square section and $fD_2,fD_3$ are curvilinear sections).
If not (fig. \ref{f:wh3} right),
then we need to first achieve this property by an analogue of the double Whitney trick applied to $fD_1$ and $fD_2$.
This is analogous to `surgery of the intersection' $fD_1\cap fD_2$, see details in Remark \ref{r:mw}, in Lemmas \ref{l:surg}, \ref{t:gdi} below and in \cite[Remark 1.16.a]{AMSW}.

More precisely, the Mabillard-Wagner Theorem \ref{t:mawa} is easily deduced \cite[\S1.3]{AMSW} from the following theorem by induction.
The Local Disjunction Theorem \ref{l:ld+3} is proved similarly.

\begin{figure}[h]
\centerline{\includegraphics[width=7cm]{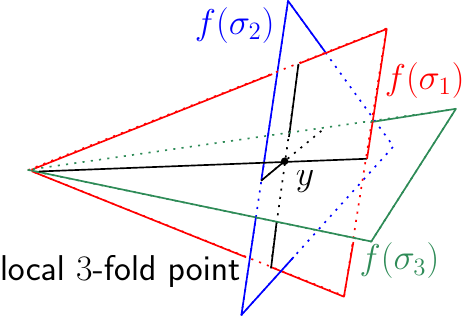}}
\caption{A non-global 3-fold point}\label{f:gl3}
\end{figure}

Let $K$ be a complex and $f\colon K\to \R^d$ a map.
A point $y\in\R^d$ is a {\bf global $r$-fold point} of $f$ if $y$ has $r$ preimages lying in pairwise disjoint simplices of $K$, that is, $y\in f(\sigma_1)\cap\ldots\cap f(\sigma_r)$ and $\sigma_i\cap\sigma_j=\emptyset$ for $i\ne j$.
See fig. \ref{f:gl3}.
Thus, $f$ is an almost $r$-embedding if and only if it has no global $r$-fold points.

\begin{theorem}[Global Disjunction]\label{t:elim3}
Let $f:K\to\R^{kr}$ be a general position PS map of a $k(r-1)$-complex $K$, let $\sigma_1,\ldots,\sigma_r$ be pairwise disjoint simplices of $K$, and $x,y\in f\sigma_1\cap\ldots\cap f\sigma_r$ two global $r$-fold points of opposite $r$-intersections signs.

Then there is a general position PS map $g\colon K\to\R^{kr}$ such that $f=g$ on $K-\Int(\sigma_1\sqcup\dots\sqcup\sigma_r)$, and $g$ has the same global $r$-fold points with the same signs as $f$ except that $x,y$ are not global $r$-fold points of $g$.
\end{theorem}


For a proof of the Global Disjunction Theorem \ref{t:elim3} we need the following lemmas, the first of which is well-known (see the references in the paragraph after Theorem \ref{l:ld+3}).

Throughout this section we fix orientations on the disks $D^m$.

\begin{lemma}[surgery of the intersection]\label{l:surg}
Assume that $d-3\ge p,q$ and that $f:D^p\to D^d$, $g:D^q\to D^d$ are proper general position embeddings such
that $fD^p\cap gD^q$ is a proper submanifold (possibly disconnected) of $D^d$ containing points $x,y$.

(a)  If $p+q>d$, then there is a proper general position map $f':D^p\to D^d$ such that:

$\bullet$ $f'=f$ on $\partial D^p$ and on a neighborhood of the set $\{f^{-1}x,f^{-1}y\}$;

$\bullet$ $x$ and $y$ lie in the interior of an embedded $(p+q-d)$-disk  contained in $f'D^p\cap gD^q$.

(b) If $p+q=d$, \ $fD^p\cap gD^q=\{x,y\}$ and $x,y$ have opposite $2$-intersection signs, then there is a general position map $f':D^p\to D^d$ such that $f'=f$ on $\partial D^p$ and $f'D^p\cap gD^q=\emptyset$.
\end{lemma}

\begin{lemma}\label{t:gdi}
Under the assumptions of the Global Disjunction Theorem \ref{t:elim3} for each $n=1,\ldots,r-1$ there is a general position PL map $f'\colon K\to\R^{kr}$ such that

$\bullet$ $f=f'$ on $K-\Int(\sigma_1\sqcup\dots\sqcup\sigma_r)$,

$\bullet$ $x,y$ lie in the interior of an embedded $k(r-n)$-disk  contained in
$f'\sigma_1\cap\ldots\cap f'\sigma_n$, and

$\bullet$ $f'$ has the same global $r$-fold points with the same signs as $f$.
\end{lemma}

\begin{proof}
The proof is by induction on $n$.
The base $n=1$ follows by setting $f'=f$.
The required disk is then a small regular neighborhood in $f\sigma_1$ of a path in $f\sigma_1$ joining $x$ to $y$
and avoiding the self-intersection set $\{x\in K\ :\ | f^{-1}fx| \geq 2\}$ of $f$.

In order to prove the inductive step assume that $n\ge2$ and the points $x,y$ lie in the interior of an embedded $k(r-n+1)$-disk $\sigma_-\subset f\sigma_1\cap\ldots\cap f\sigma_{n-1}$.
By general position
$$\dim(\sigma_-\cap f\sigma_n)\le k(r-n+1)+k(r-1)-kr=k(r-n).$$
Since $f$-preimages of $x$ lie in the interiors of $\sigma_1,\ldots,\sigma_r$, the intersections of $f\sigma_i$ and small regular neighborhoods of $x,y$ in $B$ equal to the intersections of affine spaces and the neighborhoods.
Hence the regular neighborhoods of $x,y$ in $\sigma_-\cap f\sigma_n$ are $k(r-n)$-balls.

Take points $x',y'$ in these balls.
Take general position paths $\lambda_+\subset f\Int\sigma_n$ and $\lambda_-\subset \sigma_-$
joining $x'$ to $y'$.
By general position dimension of the self-intersection set of $f$  does not exceed $2k(r-1)-kr<k(r-1)-1$.
So the union $\lambda_+\cup \lambda_-$ is an embedded circle in $\R^{kr}$.
Since $r\ge2$, we have $kr\ge4$.
Hence by general position this circle bounds an embedded $2$-disk $\delta\subset\R^{kr}$.
Since $k(r-1)+2<kr$, by general position $\delta\cap fK=\lambda_+\cup\lambda_-.$
Let $O\delta$ be a small regular neighborhood of $\delta$ in $\R^{kr}$.
Then $O\delta$ is a $kr$-ball and $f^{-1}\beta$ is the union of

$\bullet$ a regular neighborhood $D_n\cong D^{k(r-1)}$ of the arc $\sigma_n\cap f^{-1}\lambda_+$ in $\sigma_n$;

$\bullet$ regular neighborhoods $D_i\cong D^{k(r-1)}$ of the arcs $\sigma_i\cap f^{-1}\lambda_-$ in $\sigma_i$
for each $i=1,\ldots,n-1$.

Then $f|_{D_i}:D_i\to O\delta$ is proper for each $i=1,\ldots,n$, and $O\delta\cap f\sigma_-$ is a proper $k(r-n+1)$-ball in $O\delta$.
Since $n\ge2$, by general position all global $n$-fold points of $f$ in $O\delta$ are contained in
$fD_1\cap\ldots\cap fD_n$.

Since the regular neighborhoods of $x,y$ in $\sigma_-\cap f\sigma_n$ are $k(r-n)$-balls,
the set $O\delta\cap f\sigma_-\cap fD_n$ is a proper $k(r-n)$-submanifold of $O\delta$.
Hence we can apply the Embedded Surgery Lemma \ref{l:surg}.a to the disks $fD_n$ and $O\delta\cap f\sigma_-$ in $O\delta$.
For the obtained map $f':D_n\to O\delta$ the points $x,y\in f'\sigma_1\cap\ldots\cap f'\sigma_r$
are two global $r$-fold points of opposite $r$-intersections signs, lying in the interior of an embedded
$k(r-n)$-disk  contained in $\sigma_-\cap f'D_n$.
Extend $f'$ by $f$ outside $D_n$.
Clearly, the first two bullet points  in the conclusion of Lemma \ref{t:gdi} are fulfilled.
Since all global $n$-fold points of $f$ in $O\delta$ are contained in $fD_1\cap\ldots\cap fD_n$,
 the disk $f'D_n$ may only form global $n$-fold points of $f'$.
Since $n<r$, the third bullet point  in the conclusion of Lemma \ref{t:gdi} is fulfilled.
Thus the map $f'$ is as required.
\end{proof}

\begin{proof}[Proof of the Global Disjunction Theorem \ref{t:elim3}]
By Lemma \ref{t:gdi} for $n=r-1$ we may assume that the points $x,y$ lie in the interior of an embedded $k$-disk $\sigma_-\subset f\sigma_1\cap\ldots\cap f\sigma_{r-1}$.
Choose orientations of $\sigma_1,\ldots,\sigma_{r-1}$.
These orientations define an orientation on $\sigma_-$ (this is analogous to the definition of the $r$-intersection sign given in \S\ref{s:plan}, cf. \cite[\S2.2]{MW15} for a longer formal exposition).
Since $x,y\in f\sigma_1\cap\ldots\cap f\sigma_r$ have opposite $r$-intersections signs,
$x,y\in f\sigma_-\cap f\sigma_r$ have opposite $2$-intersections signs \cite[Lemma 27.cd]{MW15}.

As in the proof of Lemma \ref{t:gdi} (except that we start from $x,y$ not from $x',y'$) we construct a $kr$-ball $O\delta\subset\R^{kr}$ and $k(r-1)$-disks $D_i\subset\Int\sigma_i$
such that $x,y\in O\delta$ and $f|_{D_i}:D_i\to O\delta$ is proper.
Then $O\delta$ intersects the set of global $(r-1)$-fold points of $f$ in a neighborhood of $\lambda_-$ in $\sigma_-$.

Since $kr\ge\dim\sigma_-+k$, we can apply the Surgery of the intersection Lemma \ref{l:surg}.b to the disks $fD_r$ and
$O\delta\cap f\sigma_-$ in $O\delta$.
For the map $f':D_r\to O\delta$ obtained we have $\sigma_-\cap f'D_r=\emptyset$.
We extend $f'$ as $f$ outside $D_r$.
Clearly, $f=f'$ on $K-\Int(\sigma_1\sqcup\ldots\sqcup\sigma_r)$, and $x,y$ are not global $r$-fold points of $f'$.
Since $O\delta$ intersects the set of global $(r-1)$-fold points of $f$ in a neighborhood of $\lambda_-$
in $\sigma_-$, the map $f'$ is as required.
\end{proof}

\section{Appendix: relatives of the topological Radon theorem}\label{s:cgs}

\subsection{Statements}\label{s:cgst}

Here we exhibit relations between the following van Kampen-Flores theorem, Conway-Gordon-Sachs-Lovasz-Schrijver-Taniyama theorem, and topological Radon theorem, by presenting direct proofs of some implications between them; see Main Remark \ref{r:rela} below.
Thus we obtain alternative proofs of some of these results under the assumptions that others hold.
Direct proofs of the implications $(1^+)$ below have apparently not been published before.

Such proofs are based on interesting properties of the van Kampen and the Conway-Gordon-Sachs numbers, see Lemma \ref{l:pai} below.


Consider the following assertions for an integer $d>0$.
The topological Radon theorem ($TR_d$) is stated in \S\ref{s:intrm}.

\smallskip
($TR_d^+$) {\it `Quantitative' topological Radon theorem.} For any general position PL map $f\colon\Delta_{d+1}\to \R^d$ the number of non-ordered pairs $\{x,y\}$ of points in disjoint $k$ and $(d-k)$-faces (for all $k=0,1,\ldots,d$) such that $f(x)=f(y)$, is odd.\footnote{For $k=0,d$ this can be different from the number of intersection points in $\R^d$ of images of disjoint $k$ and $(d-k)$-faces.}

\smallskip
Disjoint closed polygonal lines $L_1,L_2$ in $\R^3$ (or, more generally, disjoint self-intersecting $k$-sphere and $\ell$-sphere in $\R^{k+\ell+1}$)
are {\it linked modulo 2} if a general position singular cone over $L_1$ intersects $L_2$ at an odd number of points \cite[\S77]{ST80}, \cite[\S4]{Sk}.

\smallskip
($VKF_d$) {\it Van Kampen-Flores theorem.}
Let $f:\Delta_{d+2}\to\R^d$ be a general position PL map.

If $d$ is even, then there are disjoint
$(d/2)$-faces whose images intersect.

If $d$ is odd, then there are disjoint $(d+1)/2$-faces whose boundaries have images linked modulo 2.\footnote{\label{f:weak} This result for $d$ even  is due to van Kampen-Flores, 1932-34; see ($VKF_{2k}$) of \S\ref{s:intrm}.
This result for $d=3$ means the existence of linked modulo 2 triangles for any PL embedding of $K_6$ into the 3-space; it is due to Conway, Gordon, and Sachs \cite{CG83, Sa81}.
This result for $d>3$ odd is due to Lovasz, Schrijver, and Taniyama \cite{LS98, Ta00}.
The references for the `quantitative' version ($VKF_d^+$) below are the same.
\newline
For $d$ odd there are the following weaker versions of ($VKF_d$):
\newline
(A) there are disjoint $(d-1)/2$ and $(d+1)/2$-faces whose images intersect.
\newline
(B) there are disjoint $(d+1)/2$-faces the images of whose boundaries are linked, that is, are not isotopic to $(d-1)/2$-spheres in the disjoint balls.
\newline
Here (B) makes sense because the restriction of $f$ to the boundary of any $(d+1)/2$-face is an embedding by general position.
The weaker version (A) follows easily from ($VKF_{d-1}$) by the link construction.
The versions (A) and (B) are equivalent \cite[end of \S6]{BFZ14} (meaning that there is a simple direct proof of the equivalence, not that any two correct statements are equivalent).
Indeed, the implication (B)$\Rightarrow$(A) is obvious.
F. Frick sketched a
simple direct proof of the implication (A)$\Rightarrow$(B).}


\smallskip
($VKF_d^+$) {\it `Quantitative' Van Kampen-Flores theorem.}
Let $f:\Delta_{d+2}\to\R^d$ be a general position PL map.

If $d$ is even, then the number of intersection points in $\R^d$ of images of disjoint $(d/2)$-faces, is odd.
I.e. the number of points $x\in \R^d$ such that $x\in f(\sigma) \cap f(\tau)$ for some disjoint $(d/2)$-faces
$\sigma, \tau$, is odd.

If $d$ is odd, then the number of unordered pairs of disjoint $(d+1)/2$-faces whose boundaries have images linked modulo 2, is odd.\footnote{Also the number of intersection points in $\R^d$ of images of disjoint $(d-1)/2$ and $(d+1)/2$-faces is even (but non-zero, see footnote \ref{f:weak}).}

\smallskip
These well-known results have very many generalizations, but citing them is beyond the purposes of this paper.

\begin{remark}[Main]\label{r:rela}
(a) There are direct proofs of the following implications
(the implications are correct because all the assertions are true).
$$\xymatrix{VKF_d \ar@{<=}[r]^{(3)\ d\text{ even}} \ar@{=>}@(dr,dl)[r]_{(3)} &
VKF_{d-1} \ar@{<=}[r]^{(2)\ d\text{ odd}} & TR_d  \ar@{=>}[r]^{(4)}
& TR_{d-1}}
\quad\text{and}$$
$$\xymatrix{VKF_{d-1}^+ \ar@{<=>}[r]_{(3^+)} & VKF_d^+ \ar@{<=>}[r]_{(1^+)} & TR_d^+
\ar@{=>}[r]_{(4^+)} & TR_{d-1}^+ }
.$$


(b) A direct proof of the right-arrow in $(1^+)$ is obtained
by extending a given general position PL map $\Delta_{d+1}\to\R^d$ to a general position PL map $\Delta_{d+2}\to\R^d$ and applying Lemma \ref{l:pai} below (this idea for $d=2$ emerged in discussions with E. Kolpakov). For $d$ even a version of Lemma \ref{l:pai} is
given in \cite[Lemma 29 and Theorem 30]{HP09} in a less explicit form.

(c) A direct proof of the left-arrow in $(1^+)$ is obtained by restricting a given general position PL map $\Delta_{d+2}\to\R^d$ to $\Delta_{d+1}$ and applying Lemma \ref{l:pai} below (this idea emerged in discussions with  S. Avvakumov).

(d) A direct proof of $(2)$ is a particular case of the
Constraint Lemma \ref{p:redu}.
For even $d$ one can similarly deduce from $(TR_d)$ the weaker properties (A,B) in footnote \ref{f:weak}.

(e) A direct proof of the left-arrow in $(3)$ for $d=2$ is given in \cite[Example 1]{Sk03}.
A direct proof of the right-arrow in $(3^+)$ for $d=3$ was obtained in \cite{Z13}.
For the general case direct proofs of the left-arrow in $(3)$ and of the right-arrow in $(3^+)$ are analogous
(see survey \cite{Sk14}).

(f) Direct proofs of $(4)$ and of the right-arrow in $(3)$ are easily obtained using the cone construction, see \cite[Lemma 1.2]{RST91} for $d=3$.
For a different proof of a more general statement see \cite[Theorem 6.5]{Me11}, cf. \cite[Example 4.7]{Me06} and \cite{Ho06}.\footnote{Melikhov kindly informed me of the following. The `only if' part of \cite[Theorem 4.2.b]{Me06} is false as stated. Its special case that is applied in the proof of \cite[Example 4.7]{Me06} is correct and is proved by the original arguments.
The statement and proof of the general case of \cite[Theorem 4.2.b]{Me06} were corrected in \cite[Theorem 4.6]{Me11}.}

(g) Following the above ideas, direct proofs of $(3^+)$ and of $(4^+)$ are formally obtained using Lemma \ref{l:rec} below.
\end{remark}

\subsection{The van Kampen and the Conway-Gordon-Sachs numbers}\label{s:cgsnum}

Here we state and prove Lemmas \ref{l:pai} and \ref{l:rec}.

Let $f\colon |K|\to \R^d$ a  general position PL map of a complex $K$.

Define {\it the van Kampen number} $v(f)\in\Z_2$ to be the parity of the number of non-ordered pairs $\{x,y\}$ of points from disjoint $k$ and $(d-k)$-faces (for all $k=0,1,\ldots,d$) such that $f(x)=f(y)$.\footnote{If $\dim K<d$, then $v(f)$ is the parity of the number of points $z\in \R^d$ such that $z\in f(\sigma) \cap f(\tau)$ for some disjoint simplices $\sigma,\tau\in K$ with $\dim\sigma+\dim\tau=d$.}
(For another recent application see \cite{ST17}.)


Define {\it the Conway-Gordon-Sachs number} $c(f)\in\Z_2$ to be the number modulo 2 of unordered pairs of  $[(d+1)/2]$- and $(d+1-[(d+1)/2])$-simplices whose boundaries have images linked modulo 2.


Denote by $\Delta_N^k$ the $k$-skeleton of the $N$-dimensional simplex.
Then

for even $d$, $(VKF_d^+)$  states that $v(f)=1$ for any general position PL map
$f\colon\Delta_{d+2}^{d/2}\to\R^d$,

for odd $d$, $(VKF_d^+)$ states that $c(f)=1$ for any general position PL map $f\colon\Delta_{d+2}\to\R^d$,

$(TR_d^+)$ states that $v(f)=1$ for each general position PL map $f\colon\Delta_{d+1}\to \R^d$.

\begin{lemma}\label{l:pai}
If $f\colon\Delta_{d+2}\to \R^d$ is a general position PL map, then
$$v(f|_{\Delta_{d+1}})=\begin{cases} v(f|_{\Delta_{d+2}^{d/2}}) & d \text{ is even}\\
c(f) & d \text{ is odd}\end{cases}.$$
\end{lemma}

\begin{proof} (A reader can first consider the cases $d=1,2,3$.)
For a face $\sigma$ of $\Delta_{d+1}$ denote the `complementary' face by $\overline\sigma$.
Then $\dim\sigma+\dim\overline\sigma=d$.
Denote by $\rho_2:\Z\to\Z_2$ the reduction modulo 2.
For subcomplexes $A,B\subset\Delta_{d+2}$ the sum of whose dimensions is $d$ denote\footnote{If $\dim A,\dim B>0$ then $A\wedge B:=\rho_2|f(A)\cap f(B)|$.
This notation reveals that both Lemmas \ref{l:pai} and \ref{l:rec} have more abstract reformulations in terms of cycles in the deleted product, without mentioning the map $f$.
The statements presented here are then obtained by evaluating the intersection cochain of $f$ on these cycles.}
$$A\wedge B:=\rho_2|\{x\in B\ :\ f(x)\in f(A)\}|\in\Z_2.$$

Denote by $*$ the vertex of $\Delta_{d+2}-\Delta_{d+1}$.
Denote by $*A\subset \Delta_{d+2}$ the cone over $A\subset \Delta_{d+1}$ with vertex $*$.
We omit parenthesis assuming that $+$ is the last-to-do operation and $\wedge$ is the next-to-last operation.
Below the summation is over (non-empty) faces of $\Delta_{d+1}$, or over their ordered pairs, satisfying the assumptions shown under the summation sign.
We have
$$
v(f|_{\Delta_{d+1}}) = \begin{cases}
\sum\limits_{\substack{\sigma \ :\\ \dim\sigma<d/2}} \sigma\wedge\overline\sigma & d\text{ odd}\\
v(f|_{\Delta_{d+1}^{d/2}}) + \sum\limits_{\substack{\sigma \ :\\ \dim\sigma<d/2}} \sigma\wedge\overline\sigma
& d\text{ even}
\end{cases},
$$
$$
v(f|_{\Delta_{d+2}^{d/2}})=
v(f|_{\Delta_{d+1}^{d/2}}) + \sum\limits_{\substack{\sigma\ :\\ \dim\sigma=d/2-1}} *\sigma\wedge\partial\overline\sigma \quad\text{for $d$ even \quad and}
$$
$$
c(f)=\sum\limits_{\substack{\sigma\ :\\ \dim\sigma=(d-1)/2}} \lk\phantom{}_2(f(\partial(*\sigma)),f(\partial\overline\sigma)) =
\sum\limits_{\substack{\sigma\ :\\ \dim\sigma=(d-1)/2}} *\sigma\wedge\partial\overline\sigma\quad
\text{for $d$ odd}.
$$
Here $\lk_2$ is the linking coefficient modulo 2.
Now the lemma holds because
$$
\left.\begin{matrix} \text{even }d: & v(f|_{\Delta_{d+1}})-v(f|_{\Delta_{d+2}^{d/2}}) \\
\text{odd }d: & v(f|_{\Delta_{d+1}})-c(f) \end{matrix}\right\}=
\sum\limits_{\substack{\sigma \ :\\ \dim\sigma\le[(d-1)/2]}} \sigma\wedge\overline\sigma+
\sum\limits_{\substack{\sigma \ :\\ \dim\sigma=[(d-1)/2]}} *\sigma\wedge\partial\overline\sigma =\ldots=0.
$$
Here the first equality is obvious, and the others are obtained by applying the following equality for $k=[(d-1)/2],\ldots, 2,1,0$:
$$
\sum\limits_{\substack{\sigma\ : \\ \dim\sigma=k}}
(\sigma\wedge \overline\sigma + *\sigma\wedge\partial\overline\sigma) \overset{(1)}=
\sum\limits_{\substack{\sigma\ : \\ \dim\sigma=k}} *\partial\sigma\wedge\overline\sigma =
\sum\limits_{\substack{(\sigma,\tau)\ :\\ \dim\sigma=k\\ \tau\subset\partial\sigma}} *\tau\wedge\overline\sigma \overset{(3)}= \sum\limits_{\substack{(\sigma,\tau)\ :\\ \dim\tau=k-1\\ \overline\sigma\subset\partial\overline\tau}} *\tau\wedge\overline\sigma=
\sum\limits_{\substack{\tau\ : \\ \dim\tau=k-1}}*\tau\wedge\partial\overline\tau.
$$
Here\footnote{For $k=0$ the last two sums run over the empty set and hence are zeroes.
The summands can be non-zero for the $(-1)$-dimensional face $\tau=\emptyset$ of $\Delta_{d+1}$, but this face is not included into the summation.}

$\bullet$ the equality (3) holds because
$\tau\subset\partial\sigma\ \Leftrightarrow\ \overline\sigma\subset\partial\overline\tau$.

$\bullet$ the equality (1) holds because
$$*\sigma\wedge\partial\overline\sigma=
\partial(*\sigma)\wedge\overline\sigma+\rho_2|\partial(f(*\sigma)\cap f(\overline\sigma))| =
\partial(*\sigma)\wedge\overline\sigma=
\sigma\wedge\overline\sigma+*\partial\sigma\wedge\overline\sigma,
$$
where the first equality is Leibniz formula, the second equality holds because by general position $f(*\sigma)\cap f(\overline\sigma)$ is a finite number of non-degenerate arcs, and the second equality holds because $\partial(*\sigma)=\sigma\cup(*\partial\sigma)$;

$\bullet$ the other equalities are obvious.
\end{proof}


\begin{lemma}\label{l:rec}
Let $f\colon\Delta_{d+2}\to \R^d$ be a  general position PL map.

(a) Let $\widehat f\colon\Delta_{d+3}\to \R^{d+1}$ be a  general position shift
of the cone over $f$.
Then
$$\begin{cases} v(\widehat f)=v(f) \\
c(\widehat f)=v(f|_{\Delta_{d+2}^{d/2}}) & d \text{ is even} \\
v(\widehat f\ |_{\Delta_{d+3}^{(d+1)/2}})=c(f) & d \text{ is odd}
\end{cases}.$$

$$(b)\qquad \begin{cases}
c(f)=v((\lk_Af)|_{\Delta_{d+1}^{(d+1)/2}}) & d \text{ is odd} \\
v(f|_{\Delta_{d+2}^{d/2}})=c(\lk_Af) & d \text{ is even}
\end{cases}.$$
Here $A$ is any vertex of $\Delta_{d+2}$ and $\lk_Af:\Delta_{d+1}\to S^{d-1}$ is any `link of $f$ at $A$'.
\end{lemma}

The statement of (b) for $d$ even is a simplified form of \cite[Lemma 18]{PT19} (this refines the statement in an earlier version \cite{Sk18u} of the current paper).
The proof of Lemma \ref{l:rec} is analogous to the particular cases $d\le4$ proved implicitly  in the references
on $(3)$ и $(3^+)$ in Remark \ref{r:rela}.e.
For a complete written proof of (b) for $d$ even see \cite[proof Lemma 18]{PT19}.
For an interesting application of (b) see \cite[Theorem 5]{PT19} and \cite{KS21}.


\begin{problem}\label{l:rela}
(a) Find a direct proof that $TR_d^+\Rightarrow VKF_{d-1}^+$, at least for odd $d$, and
that $TR_d\Rightarrow VKF_{d-1}$ for even $d$.

(b) Find direct proofs of the implications converse  to $(4)$ and $(4^+)$.
See proofs of the linear analogues in \cite{Ko18, RRS}.

(c) Recall that the topological Tverberg conjecture and the $r$-fold van Kampen-Flores conjecture are true for $r$ a prime power and are false otherwise (\S\ref{s:intrm}, \S\ref{s:plan}).
It would be interesting to obtain their `quantitative' versions ($VKF_{k,r}^+$) and ($TT_{d,r}^+$).
For example, is the following candidate to ($VKF_{k,r}^+$) true for primes or prime powers $r$?

{\it For each general position PL map $\Delta_{(kr+2)(r-1)}^{k(r-1)}\to\R^{kr}$ the sum of $r$-intersection signs of intersection points in $\R^d$ of images of $r$ pairwise disjoint $k(r-1)$-faces,
is not divisible by $r$.}

Assertion ($TT_{d,r}^+$) can presumably be dug out of \cite{BMZ09, MTW10}, cf. \cite[\S2]{Sk18}; its statement should involve only `rainbow' intersection points.


It would then be interesting to prove that $VKF_{k,r}^+\ \Rightarrow\ TT_{kr,r}^+$
(can Lemma \ref{l:pai} be generalized?).
Cf. \cite{Ko}.


(d) Find
$r$-fold analogue of the Conway-Gordon-Sachs-Lovasz-Schrijver-Taniyama theorem, that is, of ($VKF_{2k+1}$) and ($VKF_{2k+1}^+$).
For $r$-fold analogue of the Conway-Gordon-Sachs theorem, presumably not including optimal bounds, see \cite{Ne91} and an expository paper \cite{PS05}.
In particular, is it correct that {\it if no 4 of 11 points in 3-space lie in the same plane, then there are three Borromean triangles with the vertices at these points?} Cf. \cite{Ko19}.

(e) It would be interesting to obtain an analogous `quantitative' version for non-realizability of $K_5\times K_3$ in $\R^3$ and of $K_5\times K_5$ in $\R^4$, cf. \cite{Sk03}, \cite[Problem 3.9]{Sk14}.
\end{problem}

\section{Appendix: some principles of scientific discussion}\label{s:app}

\subsection{Introduction}\label{ss:app-intr}

The purpose of this section is to reassert the following principles of scientific discussion, by exhibiting violation of these principles and making some suggestions.
See also \cite[Remarks 1.1.e, 1.2.f, 2.1.de]{Sk21d}.
Most of the remarks are accessible to non-specialists.

\begin{remark}[Some principles of scientific discussion]\label{r:prin}
(a) {\it Scientific truth is established by justification and reasoning,} not by authority, majority or administrative pressure.
This is what makes science different from other respectable realms of thought (like religion, intelligence service, etc.).
Thus criticism and responsibility for criticism is vital for scientific discussion.
Without them no reliability standards could be kept.
Suppression of criticism or irresponsible criticism contribute to degeneration of a realm to non-science.
Therefore suppression of criticism or irresponsible criticism must be identified as such,
and their negative influence on science must be overcome.

In particular, the following {\it impartiality principle} should be used:
decisions ought to be based on objective criteria, rather than on bias, prejudice, or preferring to benefit one person over another for improper reasons.


This `justification' principle requires (b) below.


(b) {\it Confidential decisions should match public discussions.}
Every substantial scientific argument affecting a confidential decision should be publicly available.\footnote{Here a `decision' is a decision by an editorial board, dissertaion/diploma committee, etc.
\newline
An example on dissertation defence in Russia: if there are no public objections to a dissertation, but a dissertation committee votes against the dissertation, then the committee is dissolved.
I hope there are analogous rules or traditions in other countries.
\newline
The {\it quantity} of persons publicly supporting some point of view need not be matched by a confidential decision.}
In particular,

$\bullet$ criticism and different opinions should not be suppressed by administrative means
(e.g. by misusing the anonymous review system and writing biased negative reports on papers whose authors propose  criticism and opinions different from reviewer's);

$\bullet$ references to criticism and to different opinions should be cited not suppressed
(even if one disagrees with them).

Observe that the anonymous review system gives a convenient frame for suppression of criticism or irresponsible criticism.
It is valuable that some editors do take care of that problem.
See \cite[\S2]{Sk21d}.

Examples of violation: Remarks \ref{r:hist}, \ref{r:bbz}.f, \ref{r:bs}.c, \ref{r:bz16}.6 (see also Remark \ref{r:tribute}.a), \cite{Skr}.

(c) If original references to a result are presented (not replaced by references to a book or a survey),
then the earliest known to the author published paper should be cited (perhaps together with other references,
but no later than other references).

Examples of violation: Remarks \ref{r:bbz}.bce, \ref{r:bs}.b, \ref{r:bz16}.2.

(d) A distribution of credits for a result should be given only with (a reference to) a detailed description of contribution of different authors to the result (see e.g. footnote \ref{f:discl}; this does not apply to usage `the NN theorem' as a commonly accepted term, because such usage does not mean distribution of credits).

Examples of violation: Remarks \ref{r:bbz}.be, \ref{r:bs}.b , \ref{r:sh}.abc , \ref{r:bz16}.4,5.
\end{remark}

\begin{remark}[What is a known result?]\label{r:known}
(a) A result is known if it is published in a refereed journal, and no criticism involving specific remarks is publicly available, even if the result is unnoticed by a group of scientists.
(E.g. the Constraint Lemma \ref{p:redu} is known since \cite{Gr10}, because the only criticism \cite[\S4]{Sh18} does not justify the incompleteness, see Remark \ref{r:sh}.f.)

(b) A result is not known if it is not published in a refereed journal, or if there is a public criticism  involving specific remarks, and that criticism remains unanswered for some period.
However, an unpublished paper deserves certain proper credit when its result becomes known.
(E.g. the Mabillard-Wagner Theorem \ref{t:mawa} is known only since the survey \cite{Sk18u}, containing the first published proof, and attributing this result entirely to Mabillard-Wagner; see also \cite[\S5]{Sk17}.)


(c) {\it If $A$ is known and $A\Rightarrow B$ is known, then $B$ is known.}
If $B$ is an important statement while groups of mathematicians knowing $A$ and knowing $A\Rightarrow B$ are distant (and perhaps they use different notation for $A$ and $A\Rightarrow B$), a paper establishing $B$ should be published.
However, it would be misleading to write that `$B$ is solved by the paper'.
\footnote{The phrase `the paper observes that $B$ holds by results of X and Y',  where X and Y are papers proving $A$ and  $A\Rightarrow B$, respectively, is also misleading (at least without writing that the observation involved deduction of $B$ from $A$ and $A\Rightarrow B$, not something less trivial).
\newline
Same holds even if sophisticated exposition hides that the only novel part of the paper is deduction of $B$ from $A$ and $A\Rightarrow B$.
Same holds for other `Aristotle' deductions (like  $C$ from $A$, $A\Rightarrow B$ and $B\Rightarrow C$).
\newline
Examples of violation: Remarks \ref{r:sh}.b and \ref{r:bz16}.4.}
\end{remark}

\begin{remark}[Editorial decisions contradict public discussions]\label{r:hist}
Preliminary version of the current paper was sent to the authors of \cite{Oz, Gr10, Fr15, MW15, BFZ14, Ka, JVZ} in April 2016 (before arXiv sumbission, cf. \cite[Remark 1.4.a]{Sk21d}).
Besides approvals, I received letters containing unjustified criticism of my {\it description of references on the counterexample to the topological Tverberg conjecture}.
I asked the authors of these letters to justify their criticism by explaining which sentences
in the description are not proper and why.
I also asked to state their opinion, if different from mine, for citation in the subsequent versions of this paper.
Since then, I received neither explanations nor statements to be cited.
Below I show how this unjustified criticism reappeared through misuse of the anonymous peer review system.

The description of references to the counterexample presented in \cite{BZ16, BBZ, BS17, Sh18, BFZ} violates
important scientific principles (see Remarks \ref{r:prin} and \ref{r:known}; this is justified in Remarks \ref{r:bz16}, \ref{r:bbz}, \ref{r:bs}, \ref{r:sh}, \ref{r:mrbfz}).
Such a criticism for \cite{BZ16} is publicly available since \cite[\S4]{Sk16-2} (see \S\ref{ss:app-bz16}); this criticism also applies to \cite{BBZ, BS17, BFZ}; no reply to this criticism is publicly available.
For the description of references presented in
this survey (and in \cite[\S1.1]{JVZ}, \cite[\S1.1]{AMSW}, version 2 or higher), no criticism was publicly available before \cite[\S4]{Sh18} (for explanation why \cite[\S4]{Sh18} is incompetent see below).
However,

$\bullet$ the paper \cite{BZ16} was accepted, in spite of a referee report repeating public criticism;

$\bullet$ the paper \cite{JVZ} was rejected from the same publication, and among important reasons for rejection
the anonymous referees named that description of references in \cite{JVZ} is inconsistent with that of \cite{BZ16}.\footnote{I sent to the Editors in Summer 2016 critical remarks justifying that the corresponding parts of the reports contain important miscitations, unjustified or plainly wrong statements.
I suggested that a public clarification would improve the standards of `A Tribute to Jiri Matousek', whether
the referees would withdraw the corresponding parts of the reports, or the Editors would publicly state that
the rejection decision was made not taking into account the corresponding parts of the reports, or whatever.
I received no public clarification. See details in Remark \ref{r:tribute}.}

I am grateful to S. Shlosman for criticizing in \cite[\S4]{Sh18} the description of references to the counterexample presented in this survey.
Such publicity allows us to see that the criticism exists and is incompetent.
Indeed, a reader of \cite[\S4]{Sh18} is asked to take the proposed conclusions for granted, they are not justified by detailed references to the papers \cite{Oz, Gr10, BFZ14, MW14, Fr15, MW15, Fr15o, Sk18u}.
(In fact, the conclusions are misleading or even plainly wrong, see Remark \ref{r:sh}.)

Thus public justified criticism of the description of references to the counterexample presented in this survey
is not publicly available.
The misleading description of references is promoted by misuse of the anonymous peer review system (or by incompetent referee reports, which is the same from a reader's point of view).
The papers \cite{BBZ, BS17, BFZ} are published and left without errata.
Math Review on \cite{BFZ} is suppressed and replaced by an incompetent review, see Remarks \ref{r:mrbfz1} and  \ref{r:mrbfz}.\footnote{Additional examples are planned to be documented in an update of \cite[\S2]{Sk21d}; currently the incompetent referee reports and letters justifying incompetence are available upon request.}

This shows that violation of important scientific principles (see Remarks \ref{r:prin} and \ref{r:known}) is
a consistent policy of an influential group in the community of topological combinatorics.\footnote{I was warned about this in 2016 by a mathematician working in this field. He is in no position to confirm this publicly because of the career risk.}
{\bf In order to reassert these principles one needs to move the current discussion of the description of references to the counterexample from the (misused) anonymous peer review system to public space.
Namely, one needs

(a) to refer to (arXiv updates) of \cite{BZ16, BBZ, BS17, BFZ, Sh18} and of this survey
for such an open discussion;

(b) to publish in the journals, where papers \cite{BZ16, BBZ, BS17, BFZ, Sh18} were published, references
to the criticism of these papers in this survey, and to further discussions presumably available as in (a);

(c) to publish the suppressed Math Review of Remark \ref{r:mrbfz}, perhaps with references as in (a);

(d) to publicly exhibit and denounce such a violation of important scientific principles;

(e) to consider using (perhaps introduced gradually and partly) transparent anonymous peer review policy \cite{Skt}.

Observe that by supporting any of the suggestions (a)-(e) one only supports scientific principles, not
the particular description of references in this survey.}
I expect that an open discussion may result in errata to the publications \cite{BZ16, BBZ, BS17, BFZ, Sh18}.
I admit that it could also lead to corrections of the description of the references in this survey
(in the justifications or even in the conclusions).
\end{remark}



In remarks below references are given to numbered sections/statements of the current paper.
The notation for references is updated comparatively to the original sources (public letters and papers).
In spite of that, the references are given to the very version of a paper cited in the source.
Otherwise sources are not changed.
Sentences from the sources are given in italics and in quotation marks.

\subsection{Critical remarks on some papers}\label{ss:crit}

Before arxiv / journal publication earlier versions of the current paper were sent to I. I.~B{\'a}r{\'a}ny,
P. Blagojevi{\'c},
G. Ziegler, and S. Shlosman for their potential comments.
Since the papers criticized below were not sent to me before arxiv / journal publication,
the authors have chosen my criticism
to be public not private.

\begin{remark}[Remarks on \cite{BBZ}]\label{r:bbz}
(a) \cite[Abstract]{BBZ} `{\it ...how the same method also was instrumental in the recent spectacular construction of counterexamples for the general case.}'

This is misleading because for disproof of the topological Tverberg conjecture, of the constraint method
one needs only an elementary lemma (the Constraint Lemma \ref{p:redu} = \cite[Theorem in p. 736]{BBZ}) which is  one of the simplest and earliest known steps towards the counterexample.
See more details in Remark \ref{r:hystor}.

(b) \cite[p. 733]{BBZ} `{\it ... Florian Frick in Berlin noticed that combined with the constraint method this yields counterexamples for all $r\ge6$ that are not prime powers.}'

This is misleading for the same reason as in (a), and additionally because
F. Frick's contribution was rediscovery of the Constraint Lemma \ref{p:redu} (after \cite[2.9.c]{Gr10}).
See more details in Remark \ref{r:hystor}.\footnote{I am glad that in spite of this criticism (and analogous criticism of Remarks \ref{r:bs}.b, \ref{r:sh}.abc and \ref{r:bz16}.4,5) we have useful discussions with F. Frick, see e.g. \cite[p. 2]{AKS} and \cite[footnote 5]{KS20}. In my opinion, blind support of F. Frick's priority (for the counterexamples) contradicting earlier publications is harmful for his reputation.}

(c) \cite[p. 735]{BBZ} `{\it Indeed, the next step was taken by Murad \"Ozaydin in 1987.
In an important and influential paper [12] that
was never published, he proved the topological Tverberg
conjecture for the case when $r=p^k$ is a prime power.}'

This is misleading because reference to the first published proof \cite{vo96} is suppressed.

(d) \cite[p. 735]{BBZ} `{\it A detailed proof can be found in [3].}'

This is plainly wrong because [3]=\cite{BFZ} does not contain a proof (detailed or not) of the \"Ozaydin-Volovikov result discussed (no wonder the above citation of [3] does not mention a specific place where the detailed proof is supposed to be presented).

(e) \cite[p. 737]{BBZ} `{\it Florian Frick... realized that the theorem of Mabillard and Wagner, in
combination with the work of \"Ozaydin and the constraint method, yields all that is needed for obtaining
counterexamples.
More precisely, counterexamples to the generalized Van Kampen-Flores conjecture for nonprime
powers, which one could get from the Mabillard-Wagner theorem, would imply counterexamples to the topological
Tverberg conjecture for nonprime powers [6], [3].}'

This is misleading for the same reason as in (a,b), and additionally because

$\bullet$ the second of the above sentences essentially repeats the Constraint Lemma \ref{p:redu} as it is stated in \cite[2.9.c]{Gr10} = [7, 2.9.c], but   the first published paper \cite{Gr10} is suppressed from the list `[6], [3]' of citations (cf. Remark \ref{r:hystor}.ac);

$\bullet$ the failure of the $r$-fold Van Kampen-Flores conjecture is known since \cite{MW14}
(although it was first explicitly stated by Frick \cite{Fr15}), as explained in Remark \ref{r:hystor}.b.

(f) The paper \cite{BBZ} does not cite other descriptions of references on the counterexample \cite[\S1.1]{AMSW}, \cite[\S1.1]{JVZ}, \cite{Sk16-2}, although these references were publicly available before publication of \cite{BBZ}, and the papers \cite{AMSW, Sk16-2} were even sent to the authors of \cite{BBZ} before their arXiv submissions.
\end{remark}

\begin{remark}[Remarks on (arXiv v2 version of) \cite{BS17}]\label{r:bs}
(a) \cite[p. 10]{BS17}
`{\it In 2010, in a groundbreaking paper, Gromov \cite{Gr10} states that `The topological Tverberg theorem, whenever available, implies the (generalized) Van Kampen-Flores theorem'. This implication holds for any $r$, prime or not. Gromov also gives a proof (or rather a sketch of proof) in three lines.
A detailed proof can be found in \cite{BFZ14}.}'

This is misleading because

$\bullet$ the paper \cite{BFZ14} (as opposed to \cite{Gr10}) does not contain {\it a statement} of the implication (although it does contain its implicit proof);

$\bullet$ it is not mentioned that the detailed proof of the implication (see \cite{BFZ14} or proof of the Constraint Lemma \ref{p:redu} in \S\ref{s:state}) contains about 10 lines.

(b) \cite[pp. 10-11]{BS17}
`{\it Mabillard and Wagner almost succeeded in finding a counterexample to the topological Tverberg conjecture: what was missing was an example where the generalized Van Kampen-Flores theorem fails.
It was Florian Frick \cite{Fr15o} who realized that the above theorem and Ozaydin's example combined with the constraint method (or Gromov's remark), gives a counterexample for every non-prime power $r$.}'

This is misleading because the failure of the $r$-fold Van Kampen-Flores conjecture is known since \cite{MW14}
(although it was first explicitly stated by Frick \cite{Fr15}), as explained in Remark \ref{r:hystor}.b.


(c) The paper \cite{BS17} does not cite other descriptions of references on the counterexample \cite[\S1.1]{AMSW}, \cite[\S1.1]{JVZ}, \cite{Sk16-2}, \cite{Sk18u}, although these references were publicly available at the time of arXiv submission of \cite{BS17}, and the papers \cite{AMSW, Sk16-2} were even sent to I. B{\'a}r{\'a}ny before their arXiv submissions.
(The paper \cite{BS17} does cite \cite{AMSW}, but not as another description of references on the counterexample.)

(d) \cite[p. 8, first lines]{BS17}
`{\it Recently it has been shown by Blagojevi\'c, Frick and Ziegler \cite{BFZ14} that the topological Radon theorem implies Van Kampen-Flores.
The same implication is mentioned (somewhat implicitly) in Gromov \cite{Gr10} as well.}'

The `implicitly' is misleading because \cite{Gr10} {\it explicitly} states a more general result cited in
\cite[p. 10]{BS17} as follows: `{\it The topological Tverberg theorem, whenever available, implies the (generalized) Van Kampen-Flores theorem}'.
\end{remark}


\begin{remark}[Remarks on \cite{Sh18}]\label{r:sh}
(a) \cite[\S3]{Sh18} `{\it Florian Frick came with counterexamples to TTT, in \cite{Fr15o}.}'

This judgement is misleading because F. Frick's contribution was rediscovery of the Constraint Lemma \ref{p:redu}
(after \cite[2.9.c]{Gr10}), which is  one of the simplest and earliest known steps towards counterexample.
As opposed to the exposition of \cite{Sk18u}, this judgment is not justified by references to the papers on counterexamples to the topological Tverberg conjecture.
This is no wonder because those references do not lead to this judgement, as explained in
Remark \ref{r:hystor}, at the beginning of \S\ref{s:plan}, and in (b) below.

\smallskip
(b) \cite[\S4]{Sh18} `{\it In the report \cite{Sk18u} one reads:

"For the counterexample, papers ... by M. Ozaydin, M. Gromov, P. Blagojevic, F. Frick, G. Ziegler, I. Mabillard and U. Wagner are important."

This is a correct statement, since the paper \cite{Fr15o} appears in this list. But
this statement is an understatement; it is the truth, but not the whole truth:
– without \cite{Fr15o} there would be no counterexample. In general, it never happens
that a proof of a theorem is contained in several papers by several authors.
Usually, there exists a pivotal paper, such that the proof in question does
not exists before it, and does exist after it.}

First, the above citation of \cite{Sk18u} in \cite[\S4]{Sh18} is incomplete because it misses \cite[footnote 1]{Sk18u} (similar to footnote \ref{f:discl} of this paper) containing a disclaimer, references, and an invitation to a reader to form his/her own opinion.

Second, juxtaposition of `usually' and `never happens' (also present in the Russian version of \cite{Sh18})
is a logical fallacy.
A pivotal paper indeed {\it usually} exists,
but {\it sometimes} a proof is contained in several papers by several authors.
The latter happens when paper X proves $A$, paper Y proves $A\Rightarrow B$, while the authors of X and Y do not know of each other's achievement (and perhaps use different notation).
Then I would say that $B$ is proved in X and Y.
See more in Remark \ref{r:known}.c.\footnote{S. Shlosman's opinion can well prevail over the opinion I share and expose above.
For this it is necessary that he {\it explicitly} writes that in his opinion
{\it the deduction of $B$ from known $A$ and known $A\Rightarrow B$ is a research achievement}.
Currently this negation of my opinion is hidden by lack of detailed references to the papers
\cite{Oz, Gr10, BFZ14, MW14, Fr15, MW15, Fr15o, Sk18u}.}


Third, a proof grows through extended abstracts (like \cite{MW14}), non-refereed arXiv preprints (like \cite{MW15}), etc.
`A proof exists' after it is published in a refereed journal, see more in Remarks \ref{r:known}.a,b.

\smallskip
(c) \cite[\S4]{Sh18} `{\it In the case of counterexample to TTT such a paper was written by Florian Frick,
building on earlier results of Mabillard and Wagner.}'

Same remark as to (a).
The following shows that (neither \cite{Fr15} nor) \cite{Fr15o} is a pivotal paper such that the counterexample in question does not exists before it, and does exist after it.

$\bullet$ Complete proof of the Mabillard-Wagner Theorem \cite[Theorem 3]{MW14} (explained in Remark \ref{r:hystor}.b) did not exist at the time of writing of \cite{Fr15}, see footnote \ref{f:mw15}.
So counterexample to the topological Tverberg conjecture did not exist after \cite{Fr15}.
Thus the paper \cite{Fr15} proves the same implication as the Constraint Lemma \ref{p:redu} (first proved in \cite[2.9.c]{Gr10}), now restated as `the Mabillard-Wagner Theorem implies counterexamples to the topological Tverberg conjecture'.

$\bullet$ The paper \cite{Fr15o} does not refer to a published proof of the Mabillard-Wagner Theorem
(it only refers to \cite{MW14, MW15}); such a published proof did not exist at the time of writing of \cite{Fr15o} (see footnote \ref{f:mw15}).
So counterexample to the topological Tverberg conjecture did not exist after \cite{Fr15o}.

$\bullet$ If someone is willing to recognize \cite{MW15} as a reliable paper (and so to publicly answer public questions on that paper), then one has to recognize \cite{MW15} as a pivotal paper such that the counterexample  in question does not exists before it, and does exist after it.
Indeed, the paper \cite{MW15} presents two counterexamples, one based on the earlier paper \cite{Fr15} (so implicitly on \cite[2.9.c]{Gr10}), and the other independent of \cite{Fr15}.
The paper \cite{MW15} appeared earlier than \cite{Fr15o}.

\smallskip
(d) \cite[\S4]{Sh18} `{\it The result of Frick is based on the constraint method, developed in 2015
in \cite{BFZ14}. Again, the paper \cite{Sk18u} credits M. Gromov for the discovery of this
method earlier, in his 2010 paper \cite{Gr10}.}'

This is plainly wrong (and is not justified by reference to a specific place in \cite{Sk18u}).
The paper \cite{Sk18u} credits M. Gromov not for the constraint  method, but for the Constraint Lemma \ref{p:redu}.
This elementary lemma is the only part of the constraint method required for disproof of the topological Tverberg conjecture.
See more details in Remark \ref{r:hystor}.ac.

\smallskip
(e) \cite[\S4]{Sh18} `{\it Yet, crediting this (110-pages-) paper for the constraint method (referring to a sketch in the discussion section there) is rather detrimental to this paper than otherwise.}

As (d) explains, \cite{Sk18u} does not credit \cite{Gr10} for the constraint  method.
For the Constraint Lemma \ref{p:redu}, \cite{Sk18u} refers to the short subsection \cite[2.9.c]{Gr10}, not to the 110-pages-paper.

\smallskip
(f) \cite[\S4]{Sh18} `{\it This halfpage sketch of a proof is incomplete...}

No specific critical remarks justifying incompleteness of the proof are presented.
In itself, short length of an argument cannot contribute to calling it incomplete.
Proof of the Constraint Lemma  \ref{p:redu} (\S\ref{s:state}) has about ten lines.


\smallskip
(g) \cite[\S4]{Sh18} `{\it ... and contains typos, so much so that an attempt of correcting them is made in \cite{Sk18u}. Out of respect to Mikhail Gromov it would have been better not to initiate this discussion.}

I share opinion that correction of typos shows no disrespect to the author of the corrected paper.\footnote{A reader inclined to think in moral terms can judge whether calling a proof in a published paper incomplete, without indicating any specific problem, shows respect to the author of the paper (and to principles of scientific discussion).}

\smallskip
(h) \cite[Footnote in p. 1]{Sh18} `{\it This paper clarifies some bibliographical issues in the review paper \cite{Sk18u} in this volume, which clarification I was entrusted to make by the Editorial Board of Russian Math.
Surveys.}'

S. Shlosman misused his position of the Editorial member of Russian Math. Surveys to break the impartiality principle (see Remark \ref{r:prin}.a).
Indeed, \cite[\S4]{Sh18} is polemics on description of references, while such polemics was deleted from \cite{Sk18u} upon request of an anonymous referee supported by the Editors.\footnote{Polemics was requested to be deleted as such.
No criticism of (the submitted version of) \cite{Sk18u} was presented in the report, so the polemics was not shown to be incompetent. On the other hand, the criticism of \cite[\S4]{Sh18} cannot be judged to be competent, because it was not sent before its publication to the author of the criticized paper for his response.
In fact, the criticism of \cite[\S4]{Sh18} is incompetent, which is shown by (b)-(f).}
Although the paper \cite{Sh18} is presumably written upon the invitation of the Editorial Board, the incompetent opinion of \cite[\S4]{Sh18} is not necessarily supported by the Editorial Board.

\smallskip
(i) Orientations on configuration spaces related to the topological Tverberg Conjecture are discussed in some papers before \cite[\S2]{Sh18}, see e.g. \cite{MTW10} and the references therein.
These papers are not cited, so novelty of results \cite[\S2]{Sh18} is not checked.
\end{remark}

\comment


This is in tune with Poincar\'e's quote:
. . .  il n'y a plus des probl\` emes r\'esolus et d'autres qui ne le sont pas,
il y a seulement des probl\`mes plus ou moins r\'esolus.

Let me show that even if we assume that the `pivotal paper principle' of S. Shlosman has no exclusions, we obtain that the counterexample has to be attributed to  \"Ozaydin and Mabillard-Wagner (not to Frick).
Since we are in an argument involving distribution of credits,
we have to be accurate and consider only proofs published in refereed journals.

The implication {\it a counterexample to the generalized Van Kampen-Flores conjecture implies a counterexample to the topological Tverberg conjecture} is first published in \cite{Gr10}.\footnote{In \cite[\S4]{Sh18} it is written that the proof is incomplete without any specific critical remarks, so the judgement `incomplete' has to be dismissed.
In any case, this is published in the survey \cite{S} with credits to Gromov, and to a rediscovery by  Blagojevi\'c-Frick-Ziegler.}
{\it A counterexample to the generalized Van Kampen-Flores conjecture (VKF)} is first published in the survey \cite{S}
with credits to \"Ozaydin and Mabillard-Wagner (whose papers are not published in refereed journals before the survey \cite{Sk18u}; as far as I know, they remain unpublished).
The counterexample to (VKF) immediately follows from results by \"Ozaydin and Mabillard-Wagner (as explained in Remark \ref{r:hystor} and at the beginning of \S\ref{s:plan}).
The survey \cite{S} contains complete proofs (slightly different than the original) of those results.
So counterexample to the topological Tverberg conjecture does not exist before the survey \cite{S}, and exists after it.
[delete the last two phrases?]

The above formal point of view leads to the same conclusion as substantial point of view.
Indeed,

$\bullet$ the above implication proved in \cite{Gr10} is simple; no wonder it was rediscovered.

$\bullet$ the results of \"Ozaydin and especially of Mabillard-Wagner are complicated;
no wonder they required years to be written, refereed, rewritten, rerefereed,  etc before publication.

---------- Forwarded message ---------
From: <umn@mi.ras.ru>
Date: Fri, Jun 9, 2017 at 9:45 PM
Subject: Re: Fwd: submission
To: arkadiy skopenkov <askopenkov@gmail.com>
Cc: <buchstab@mi.ras.ru>

 Аркадий, посылаем вам ответ второго рецензента.

 Работа состоит из двух частей, математической и приоритетной, если
 можно так выразиться.
 В математической части автор излагает довольно подробно и грамотно
 этапы доказательства топологической теоремы Тверберга для степеней
 простых, и объясняет процедуру построения контрпримера для других
 натуральных чисел. Эта часть выглядит разумно и может быть опубликована
 в Успехах.
 Кроме того, статья содержит описание того, кто, какие теоремы, и в
 каком хронолигоческом порядке доказал; этот вопрос несколько раз
 обсуждается также в подстрочных комментариях. Мне представляется что это
 обсуждение и полемика по этому поводу совершенно неуместна на страницах
 Успехов, и все вопросы приритета должны быть полностью удалены из
 текста, если статья будет принята Успехами.

\endcomment

\subsection{Publication of \cite{BZ16} and rejection of \cite{JVZ}}\label{ss:app-bz16}



\begin{remark}[Remarks on \cite{BZ16}]\label{r:bz16}\footnote{These remarks appeared along what I consider a proper way: from private criticism to public criticism to criticism in a referee report.
Remarks (1)-(6) (in a bit different form \cite[\S5]{Sk16-2}) were sent to the authors in May 2016.
These remarks were published in July 2016 \cite[\S5]{Sk16-2} (for reasons explained in \cite[footnote 7]{Sk16-2}).
These remarks, together with most part of remarks (7), (8) constituted my referee report to \cite{BZ16-1}.
Remarks (1)-(11) almost coincide with remarks from my report (October 2016) on a revision of \cite{BZ16-1}.
Although the remarks concern earlier versions of \cite{BZ16}, they are not taken into account, so they are pertinent to \cite{BZ16}.
\newline
My reports were prepared on request of the Editors M. Loebl, J. Ne\v set\v ril and R. Thomas of
{\it `A Journey through Discrete Mathematics. A Tribute to Ji\v ri Matou\v sek'.}
They advised referees to apply high standards for survey papers in terms of clarity of exposition, appropriate for a top combinatorial journal.
Hence, based on Remarks (1)-(11), I recommended rejection or major revision
(and warned that remarks (4), (5) and (6) prevent publication in any scientific journal).
See Remarks \ref{r:hist} and \ref{r:tribute}.}
(1) I suggest to explicitly state the following lemma, and present its proof.

A map  $f\colon K\to \R^d$ of a simplicial complex is an {\it almost $r$-embedding} if $f(\sigma_1)\cap \ldots \cap f(\sigma_r)=\emptyset$ whenever $\sigma_1,\ldots,\sigma_r$ are pairwise disjoint simplices of $K$.

{\bf Constraint Lemma.}
(Gromov, Blagojevi\'c, Frick, Ziegler, \cite[Theorem in p. 736]{BBZ}, cf. the Constraint Lemma \ref{p:redu})
{\it If $k,r$ are integers and there is an almost $r$-embedding of the $k(r-1)$-skeleton of the $(kr+2)(r-1)$-simplex in $\R^{kr}$, then there is an almost $r$-embedding of the $(kr+2)(r-1)$-simplex in $\R^{kr+1}$.}

Currently the lemma is not stated, and is proved separately for $r$ a prime power \cite[\S4.1]{BZ16} or not \cite[\S5.2]{BZ16} in the proof of other results; neither case of the lemma uses the fact that $r$ is a prime power or not.
This is misleading, cf. remark (4) below.

Although explicit statement and proof not repeated twice will make exposition only a little shorter,
this will make exposition substantially  clearer.

\smallskip
(2)  \cite[\S1, 2nd paragraph]{BZ16}:
`{\it The topological Tverberg conjecture was extended to the case when $r$ is a prime power
by Murad \"Ozaydin in an unpublished paper from 1987 [38].}'

I suggest to add a reference [47]=\cite{vo96} to the {\it first published} proof, otherwise the sentence is misleading.\footnote{The authors responded to this suggestion of \cite[\S5]{Sk16-2} on \cite{BZ16-1}
by keeping the above statement in the same form and citing [47]=\cite{vo96} in \cite[\S3.3]{BZ16} not in \cite[\S1]{BZ16}.
Thus a reference to the first published proof was hidden.}

\smallskip
(3) \cite[bottom of p. 1]{BZ16}:
`{\it In a spectacular recent development, Isaac Mabillard and Uli Wagner [31, 32] have developed an $r$-fold
version of the classical `Whitney trick'...}'

I suggest to add something like `see \cite[footnote 2]{Sk16-2} [Added in 2022: this is Remark \ref{r:mw}]
for relation of Mabillard-Wagner idea to earlier references'.
Although this explanation justifies that Mabillard-Wagner's work is spectacular,
the above sentence is misleading without indicating relations to previous publications.
\footnote{The authors responded to this suggestion \cite[\S5]{Sk16-2} on \cite{BZ16-1} by citing the Whitney 1944 paper and thus ignoring many later references cited in \cite[footnote 2]{Sk16-2}.}

\smallskip
(4) \cite[pp. 1-2]{BZ16}: `{\it ...which yields the failure of the generalized van Kampen-Flores theorem when $r\ge6$ is not a prime power. Then Florian Frick observed that this indeed implies the existence of
counter-examples to the topological Tverberg conjecture [26] [13], by a lemma of Gromov [27, p. 445]...}'


This is misleading because the lemma of Gromov (i.e., the Constraint Lemma above) states that {\it the topological Tverberg conjecture implies the generalized van Kampen-Flores conjecture} \cite[2.9.c]{Gr10}, \cite[Theorem in p. 736]{BBZ}.
So it is hidden that F. Frick is being credited for `observing that $\overline B\Rightarrow \overline A$ by a lemma asserting that $A\Rightarrow B$'.
This is unfair to F. Frick: I suppose no mathematician would like being credited for this.
I suggest to write that Frick's important achievement was rediscovery of the Constraint Lemma, of which the wide community of topological combinatorics was not aware, and thus the codimension 3 restriction of Mabillard-Wagner was considered a serious obstacle before Frick.

\smallskip
(5)  \cite[\S5, 1st paragraph]{BZ16}. This is misleading and unfair to F. Frick, see (4).
I suggest to modify this paragraph consistently with suggestion to (4).

\smallskip
(6) I suggest to refer in \S1 to other expositions of references concerning
disproof of the topological Tverberg conjecture, including \cite{Sk16-2} and \cite[\S1.1]{JVZ}.
If  the authors disagree with those expositions, it would be very valuable if they explain in their survey
which sentences exactly they disagree with, and why (cf. Remarks (2)-(4) above).
It is misleading to ignore those references.\footnote{The authors responded to this suggestion of \cite[\S5]{Sk16-2} on \cite{BZ16-1} by citing [42]=\cite{Sk16-2} only in p. 24, in connection with the Mabillard-Wagner Theorem not with the history of disproof of the topological Tverberg conjecture.
Thus a reference to a different point of view was hidden.}

\smallskip
(7) The most complicated (and thus the most interesting!) part of disproof of the topological Tverberg conjecture is the Mabillard-Wagner Theorem in \cite[pp. 24-25]{BZ16}.
Lack of idea and plan of the proof of this result is an annoying gap in a survey intended to cover the
topological Tverberg conjecture.
I suggest to add idea and plan of the proof, at least up to the extent given in \cite[\S4]{Sk16-2}.
In fact, an alternative complete proof in \cite{AMSW} [Added in 2022: or in \S\ref{s:mw}, \cite[\S3.5]{Sk16-4}]
is shorter than the included proof of \cite[Theorem 3.11]{BZ16}.

\smallskip
(8) A survey on the `holy grail' topological Tverberg conjecture will be interesting to a wide audience of mathematicians.
Most of them would have only a vague idea of what spectral sequences or equivariant obstructions are.
So proofs of most important results of the survey (the topological Tverberg conjecture for prime powers $r$ and an important step towards a counterexample, \cite[Theorems 3.11 and 5.1]{BZ16} of \"Ozaydin-Volovikov and \"Ozaydin) are completely unaccessible to them.
It is annoying that the exposition in \cite[\S3]{BZ16} of this {\it survey}  is harder to follow than the exposition of the same results in original {\it research} papers [8, 38, 47]=\cite{BSS, Oz, vo96}.

Theorems 3.11 and 5.1 of \cite{BZ16}  are striking results worth an attempt to write `proof from the book' by
such masters of clear exposition as Blagoevi\'c and Ziegler.
In order to understand their proofs {\it in full details} one needs of course
spectral sequences and equivariant obstructions.
However, the proofs can be structured so as to make {\it the plan} and some parts of the proof accessible to non-specialists, see \cite[\S3]{Sk16-2}, \cite[\S2.3, \S5.2]{Sk16-3}.
(Then the complete proofs will be easier to read for a specialist.)
This can be achieved by {\it stating explicitly} a result proved by a theory in terms not involving this theory.

The proof of Theorem 3.11 in \cite{BZ16} starts with spectral sequences \cite[p. 11]{BZ16}.
This antagonizes non-specialists (in algebraic topology) because they do not know this notion, and neither a reference nor dictionary explains it.
This also antagonizes specialists because they can easily do exercises involving spectral sequences, once the statements \cite[Lemma 2.8]{Sk16-3} are explicitly given, so what they want is proof modulo these statements.

Even for a specialist it is hard to read the 3-page \cite[proof of Theorem 3.11]{BZ16} because it is not well-structured.
E.g. in \cite[p. 12]{BZ16} before the localization theorem it is proved that the map in cohomology is a monomorphism, and after the localization theorem it is observed that it is not, so it seems that we already have a contradiction and \cite[(4)]{BZ16} is not required.
The problem is that there is a statement (`$a\ne0$') which is proved from the contrary,
but neither this statement nor its proof is split from the remaining text.

I also suggest avoiding unnecessary sophisticated language.
E.g. although $S(W_r^{\oplus d})$ is a standard notation for specialists, I suggest to change it
to the simple-minded notation $S^{d(r-1)-1}_{\Sigma_r}$ defined as the set of real $d\times r$-matrices with sum in each row equal to zero, and the sum of squares of the matrix elements equal to 1 \cite{Sk16-2}.
Then the phrases like \cite[second bullet points in p. 10, p. 21 and p. 23]{BZ16} will be unnecessary,
the reader would not have to consult dictionary for `$(N-r)$-simple',
and in the Dold Theorem in \cite[p. 10]{BZ16} one can replace $Y$ by $S^n$.
All this would make the matter more accessible.

\smallskip
(9) \cite[p. 11-12, Proof of Theorem 3.11, p. 23, Proof of Theorem 5.1]{BZ16}.

The notation is unreadable, see display formulas.
It is easy to keep the notation short \cite[\S2.3, \S5.2]{Sk16-3}.

\smallskip
(10) \cite[\S1, 2nd paragraph]{BZ16}: {\it `...with complete proofs for all of the results...'}.

This is misleading because many theorems (of Dold, Leray, most importantly of Mabillard-Wagner) are used without proofs.

\smallskip
(11) The dictionary at the end of the paper \cite{BZ16} does not contain some definitions used in the paper and unknown to  non-specialists. E.g. `spectral sequences' from \cite[p. 11]{BZ16} on; the symbols $\mathcal H$ in
\cite[pp. 11, 23, 34]{BZ16} (for cohomology with twisted coefficients) are not defined; this makes definitions of `primary obstruction', `exact obstruction sequence', `obstruction element' in \cite[p. 34]{BZ16} unclear.


Those definitions the dictionary does contain, are given in more generality than required for this survey,
which makes it harder to follow.


Terms defined in the dictionary are sometimes used without references to the dictionary
(e.g. \cite[\S1, line 4]{BZ16}: `$n$-connected').

So the dictionary is almost of no use to non-specialist.
\end{remark}

\begin{remark}[Letters of 2017; cf. Remark \ref{r:hist}]\label{r:tribute}
(a)  Dear Professors Martin Loebl, Jarik Nesetril and Robin Thomas,
Editors of ``Journey Through Discrete Mathematics. A Tribute to Jiri Matousek''

Please note that the paper \cite{Sk16-4} [Added in 2022: and other arXiv versions] contains (\S4) critical remarks to the paper \cite{BZ16}.

It would be very valuable if you could help to improve the following situation.

* There are two points of view on references concerning disproof of Topological Tverberg Conjecture.
One of them is exposed in \cite{Sk16-4} [Added in 2022: and in other arXiv versions], the other in \cite{BZ16}.

* The survey \cite{Sk16-4} refers to the other point of view (end of \S3.2).
It also explains why the other point of view is inappropriate (\S4).
This explanation was never publicly questioned.
Preliminary version of \cite{Sk16-4} was sent in April 2016 to M. Ozaydin, M. Gromov, F. Frick, I. Mabillard, U. Wagner, P. Blagojevic, G. Ziegler, G. Kalai, D. Jojic, S. Vrecica, R. Zivaljevic.
Besides approvals, I received letters implying that the description of references there is not proper.
I asked the authors of criticism to explain which sentences in the description are not proper and why.
I also asked to state their opinion, if different from mine, for citation in my paper.
I received neither explanations nor statements to be cited.

* The survey \cite{BZ16} does not refer to the other point of view.
Cf. \cite[\S4, (6) and footnote 9]{Sk16-4} [Added in 2022: this is Remark \ref{r:bz16}.6.].

* The paper \cite{JVZ} was rejected from `A Tribute to Ji\v ri Matou\v sek', and among important reasons for rejection anonymous referees named that description of references in \cite{JVZ} is inconsistent with that of \cite{BZ16}.
(This point will automatically disappear after your clarification which I suggest in the other letter.)

* The survey \cite{BZ16} was accepted to `A Tribute to Ji\v ri Matou\v sek' in spite of a referee report indicating misleading presentation of references (and also poor exposition of some parts).
This report was made public and never publicly questioned.
The Editors might have received a number of positive reports from other referees.
However, it is the Editors' own decision whether to require, from a survey paper, to correct misleading presentation of references, or at least to indicate the existence of other point of view, even if only one of the many referees explained that the presentation misleading.
Just as it could be Editors' own decision whether to require, from a research paper,
to correct a mathematical mistake, even if only one of the many referees found it.

In my opinion, suppression of references to other point of view is counter-productive for research.
It makes me personally very sad that this suppression is made a part of `A Tribute to Ji\v ri Matou\v sek'.
I would be glad to correct the corresponding part of \cite[Remark 4.1]{Sk16-4} if this unlucky situation is improved.
Adding Editors' or authors' `added in proof' note referring to the other point of view would be both a simple and an  effective thing to do, if the Editors choose not to consider critical remarks of \cite[\S4]{Sk16-4} in detail.

Sincerely Yours, Arkadiy Skopenkov.

\smallskip
(b) Dear Professor Robin Thomas,

Thank you for your quick reply.
It is nice that Rade's bringing those reports to light provided clarification.
As far as I understand, the Editors' rejection decision was made disregarding those parts of the reports that were criticized in my letter (the criticism being ignored by one of the referees and not denied by the other).
It would be nice if the Editors could state this explicitly in a letter addressed to all the recipients of Rade's and mine June 2016 letters.
Otherwise those recipients would have the impression that the misleading parts of the reports (approx. 4/5 of specific remarks in Report 1 and 1/4 of specific remarks in Report 2) was among important reasons for rejection.
In my opinion, this would be counter-productive for development of mathematics.
It would make me personally very sad if such a thing would happen during editorial process for `A Tribute to Ji\v ri Matou\v sek'.

Sincerely Yours, Arkadiy Skopenkov

P.S. Since the paper \cite{Sk16-4} was updated before I received your clarification, in Remark 4.1 there it is written `The paper \cite{JVZ} was rejected from the same publication, and among important reasons for rejection the anonymous referees named that description of references in \cite{JVZ} is inconsistent with that of \cite{BZ16}.'
I would be glad to change this phrase (or to delete it upon your wish) after your clarification would be sent to all the recipients of Rade's and mine June 2016 letters.
A.S.
\end{remark}

\comment

On 3/2/17, Robin Thomas <thomas@math.gatech.edu> wrote:
> Dear Professor Skopenkov,
>
> Only one of the referees responded. His/her response is copied below.
> Sincerely,
>
> Robin Thomas
>
> ------
>
> as far as I can see, the decision not to publish the paper is not
> challenged.
> This is in accord to my recommendation
>  > You could reject it based on the fact that
>  > the topics treated are very special,
>  > the results are not substantial,
>  > and thus the paper is not good enough.
>
> The discussions about the referee reports, wordings of history, etc.,
> are clearly not productive.
> So in my reviewer role I don’t want to comment further.

Recall that a map  $f\colon K\to \R^d$ of a simplicial complex is called an {\it almost $r$-embedding} if
$f(\sigma_1)\cap \ldots \cap f(\sigma_r)=\emptyset$ whenever $\sigma_1,\ldots,\sigma_r$ are pairwise disjoint
simplices of $K$.

\smallskip
(1) I suggest to explicitly state the following proposition, and present its proof.

{\bf Proposition.} (Gromov, Blagojevi\'c, Frick, Ziegler; this is essentially the implication
`$(A_r)\Rightarrow(T_r)$' of \S\ref{s:plan}.)
{\it If $k,r$ are integers and there is an almost $r$-embedding of the $k(r-1)$-skeleton of
the $(kr+2)(r-1)$-simplex in $\R^{kr}$, then there is an almost
$r$-embedding of the $(kr+2)(r-1)$-simplex in $\R^{kr+1}$.}

Currently the Proposition is not stated in \cite{BZ}, and is proved separately for $r$ a prime power \cite[\S4.1]{BZ} or not \cite[\S5]{BZ} in the proof of other results; neither case of the Proposition uses the fact that
$r$ is a prime power or not.
This is misleading, cf. remark (4) below.

Although explicit statement and proof not repeated twice will make exposition only a a little shorter,
this will make exposition substantially  clearer.

\smallskip
(2) P. 1, 2nd paragraph. `{\it The topological Tverberg conjecture was extended to the case when $r$ is a prime power
by Murad \"Ozaydin in an unpublished paper from 1987 [36].}'

I suggest to add a reference [43] to the {\it first published} proof, otherwise the sentence is misleading.

\smallskip
(3) P. 2, 1st paragraph.
`{\it In a spectacular recent development, Isaac Mabillard and Uli Wagner [32, 33] have developed an $r$-fold
version of the classical `Whitney trick' and which yields the failure of the generalized van Kampen--Flores
theorem when $r\ge6$ is not a prime power.}'

I suggest to add something like `see [S, footnote 2] for relation of Mabillard-Wagner idea to earlier references'.
Although this explanation justifies that Mabillard-Wagner's work is spectacular,
the above sentence is misleading without indicating relations to previous publications.

\smallskip
(4) P. 2, 1st paragraph. `{\it ...and Florian Frick observed that this indeed implies the existence of
counter-examples to the topological Tverberg conjecture [25] [12].}'

I suggest to add that this observation was earlier done by Gromov [26, 2.9.c] and Blagojevi\' c-Frick-Ziegler
[13, Lemma 4.1.iii and 4.2] (before the result of Mabillard-Wagner), of which the wide community of topological combinatorics was not aware,
and that the codimension 3 restriction of Mabillard-Wagner was considered a serious obstacle before Frick.

Although this addition justifies that Frick's observation is important, the above sentence is misleading without indicating relations to previous publications.
Indeed, `the failure of the generalized van Kampen--Flores theorem when $r\ge6$ is not a prime power'
means that

{\it there is an almost $r$-embedding of the $k(r-1)$-skeleton of the $(kr+2)(r-1)$-simplex in $\R^{kr}$}.

`The existence of counter-examples to the topological Tverberg conjecture'
means that

{\it there is an almost $r$-embedding of the $(kr+2)(r-1)$-simplex in $\R^{kr+1}$.}

Thus  Frick's observation is the above Proposition.
The Proposition was explicitly proved in [26, 2.9.c] and implicitly in [13], see subsection `on references' of [S,\S1].

\smallskip
(5) P. 21, last paragraph. `{\it ... We present counterexamples
to the topological Tverberg conjecture for non prime powers derived by Frick [25] [11] from the
remarkable works of \"Ozaydin [36] and of Mabillard and Wagner [32] [33].}'

I suggest to modify this phrase in a way consistent with suggestion (4).

\smallskip
(6) I suggest to add references to other expositions of references concerning
disproof of the topological Tverberg conjecture, including [S] and \cite{JVZ}.
If
the authors
disagree with those expositions, it would be very valuable if
they explain
which sentences exactly
they
disagree with, and why (cf. remarks (2)-(4) above).
However, it is misleading to ignore those references.

\endcomment

\subsection{An incompetent published and a suppressed reviews on \cite{BFZ}}\label{ss:app-bfz}

In Remark \ref{r:mrbfz1} I show that Math Review \cite{Ba21} to \cite{BFZ} is incompetent.
In Remark \ref{r:mrbfz} I present a review which a reader will hopefully find competent.
I submitted this review to Mathematical Reviews in 2020 upon invitation,
but neither it was published nor I received a rejection letter.
Remark \ref{r:mrbfzed} is my letter to the Editors of Mathematical Reviews sent together with the review of Remark \ref{r:mrbfz}.

\begin{remark}[A report on \cite{Ba21}]\label{r:mrbfz1}
(a) {\it `the first counterexamples were provided by the second author'}

This is misleading because the counterexamples were found in a series of papers by M.~\"Ozaydin, M. Gromov, P. Blagojevi\'c, F. Frick, G. Ziegler, I. Mabillard and U. Wagner, see detailed justification and references
in Remark \ref{r:hystor}; see also the beginning of \S\ref{s:plan}, \S\ref{ss:app-bz16}, and Remarks \ref{r:bbz}.be, \ref{r:bs}.b, \ref{r:sh}.abc, \ref{r:mrbfz}, \ref{r:mrbfzed}.

(b) `{\it The counterexamples were proved using the "constraint method" of Mabillard and Wagner.}'

This is wrong / misleading because

$\bullet$ the `constraint method' is due to Blagojevi{\'c}-Frick-Ziegler, not to Mabillard-Wagner;

$\bullet$ the counterexamples require not the `constraint method' in any essential generality, but a simple lemma first proved in \cite{Gr10}:
{\it the topological Tverberg theorem, whenever available, implies the van Kampen-Flores theorem}.

(c) {\it `In this paper, the authors use this method to demonstrate the failure of the topological
Tverberg conjecture for all $r\ge6$ which are not prime powers.'}

This is misleading because of the second bullet point in (b), and because \cite{BFZ} relies upon
unpublished results (with improper statements), see Remark \ref{r:mrbfz}.
\end{remark}

\begin{remark}[Suppressed Math Review on \cite{BFZ}]\label{r:mrbfz}
This paper unfairly claims disproof of the topological Tverberg conjecture (Theorem 4.2).
The claim is unfair because the argument only involves deduction of $(A)$ from $(C)$ and implications $(C)\Rightarrow(B)$ and $(B)\Rightarrow(A)$, all of them earlier known (see below).
A disproof can be found in the survey \cite{Sk18u} (which of course claims no priority for the disproof and contains careful description of relevant references).

The claim is also unfair because reference [14]=\cite{MW15} for Theorem 3.3 is unpublished
(same holds for reference [17]=\cite{Oz} for Theorem 3.4; here and below all references except
\cite{Sk18u} are from the paper under review).
Instead of something like `we checked the argument of that unpublished paper and are ready to bear responsibility for the argument', we see arxiv preprint [14]=\cite{MW15} being called `{\it the full journal version}' in \S1.
[Added in 2021: the authors also demonstrate lack of understanding of those unpublished results by stating
one of them with a superfluous dimension restriction, see below, and by not presenting even the idea of proof of the other result in their survey \cite{BZ16}.]
A published (and hopefully clear) proof of Theorems 3.3 and 3.4 is presented in \cite{Sk18u} (where of course the theorems are attributed to Mabillard-Wagner and \"Ozaydin).


Let $K$ be a (finite simplicial) complex, i.e. a union of some faces of a simplex.
A continuous map  $f\colon K\to \mathbb R^d$ is an {\it almost $r$-embedding} if $f(\sigma_1)\cap \ldots \cap f(\sigma_r)=\emptyset$ whenever $\sigma_1,\ldots,\sigma_r$ are pairwise disjoint faces of $K$.

\smallskip
{\it Theorem 4.2.} {\it If $r$ is not a prime power and $k\ge3$, then

$(A)$ there is an almost $r$-embedding of $(kr+2)(r-1)$-dimensional simplex to $\mathbb R^{kr+1}$.}

\smallskip
{\it Theorem 4.1.} {\it If $r$ is not a prime power and $k\ge3$, then

$(B)$ there is an almost $r$-embedding of any $k(r-1)$-dimensional complex to $\mathbb R^{kr}$.}


\smallskip
{\it Theorem 3.4.} {\it If $r$ is not a prime power, then

$(C)$ there is an equivariant map from the $r$-fold deleted product of any $k(r-1)$-dimensional complex
to  the $r$-fold deleted product of $\mathbb R^{kr}$.}


\smallskip
Theorem 3.4 is `{\it a simple corollary of the work of \"Ozaydin}' (see the paragraph before Theorem 3.4; this was earlier noticed in [15, \S1, Motivation \& Future Work, 2nd paragraph]=\cite{MW14} [Added in 2021: and also in \cite[p. 173, the paragraph before Theorem 3]{MW14}]; the restriction $k\ge 3$, i.e. `{\it the dimension gap $d-\dim K$ is sufficiently large}', is not used in the proof and is superfluous).

The implication $(C)\Rightarrow(B)$ for $k\ge3$ is (a particular case of) Theorem 3.3 and is `{\it a highly nontrivial recent result of Mabillard and Wagner}' (see the paragraph before Theorem 3.3).

The implication $(B)\Rightarrow(A)$ was proved in [9, 2.9.c]=\cite[2.9.c]{Gr10} and implicitly rediscovered in
[5, 8]=\cite{BFZ14, Fr15o} (for more detailed references see \cite[Remark 1.11 on the Constraint Lemma]{Sk18u}
[Added in 2021: and Remark \ref{r:hystor}]).

The exposition justifying the unfair claim is misleading.
E.g. in \S3, the paragraph before `Ingredient 1', it is not explained which ingredient is missing and why.
The paragraph before Theorem 3.2
hides the fact that the only part of the constraint method required for disproof of the topological Tverberg conjecture is the simple implication $(B)\Rightarrow(A)$ proved in [9, 2.9.c]=\cite[2.9.c]{Gr10}.
See more explanations in \cite[\S5]{Sk16-2} [Added in 2021: and in Remark \ref{r:hystor}, \S\ref{s:app}]).


Hopefully this report could stimulate appearance of an update of arXiv:1510.07984 and of errata in JEMS.
\end{remark}

\begin{remark}\label{r:mrbfzed}
Dear Editors of MR,

Please ignore this letter if you think that my review MR3959859 is acceptable in its current form.

Since the review contains a serious criticism, I had to write a thorough justification (although I do understand that opinions of the author and the Editors need not coincide, so by accepting the review, if you choose to accept it, you only confirm that the review is qualified, not that you agree with the conclusions).
Since the review has to be short, I extinguished some less important justifications to this letter.
I am willing to move some justifications from the letter to the report or vice versa, if you ask me to revise the review.

Below I use the notation introduced in the review.

\smallskip
{\bf (i). Does assembling the previously obtained results and previously known implications between them form
a new research result?}

Before I put arXiv:1605.05141v1
to arxiv, I sent it for remarks to a number of mathematicians including the authors of the paper \cite{BFZ} under review.
I received criticism similar to positive answer to the above question implying that the description of references in arXiv:1605.05141v1
is not proper.
I asked the authors of these letters to explain which exactly sentences in the description are not proper and why.
I also asked to state those authors' opinion for citation in  arXiv:1605.05141v1
(because one has to present a different opinion in one's paper even if one disagrees with it [Added in 2022: see Remark \ref{r:prin}.b]).
Since then, I received neither explanations nor statements to be cited.

The answer to the above question depends on what `assembling' is and what `previously known' is.
Hence a clear description of that is important, and the paper \cite{BFZ} does not give such a description.

I would not call a new research result deduction of $X$ from known $Y$ and known $(Y)\Rightarrow(X)$, or any other kind of such an Aristotlean deduction (even if the deduction involves using a statement on any $k$-dimensional complex for a particular complex like $k$-skeleton of a simplex).
I would call this an important expository result if $X$ is a well-known problem.
The paper \cite{BFZ} does not explicitly state that the new research result only involves deduction of $(A)$ from $(C)$ and implications $(C)\Rightarrow(B)$ and $(B)\Rightarrow(A)$, all of them earlier known.
The paper hides this by artificially sophisticated exposition and by suppressing references \cite{JVZ, Sk18u} to another exposition.
So we cannot say that referees and editors approve that this kind of deduction could be considered a new research result.

The Mabillard-Wagner Theorem 3.3 was not shown to be `previously known' in \cite{BFZ}, see the review and (ii) below.

So the paper \cite{BFZ} constitutes an unfair attempt of shifting the assembling to itself from the would-be published version of the Mabillard-Wagner paper [14]=\cite{MW15} (where besides assembling the most non-trivial part would be proved).

A rationale behind this attempt is that for the community of topological combinatorics F. Frick's announcement [8]=\cite{Fr15o} represented a reasonable way to $(A)$, i.e. to a disproof of the topological Tverberg conjecture.
The community generally believed that the strategy of [14, 15]=\cite{MW14, MW15} should work to prove $(B)$.
Some part of the community neglected a distinction between `this should work' and `this is a reliable result' (possibly because of the wish to announce a disproof of a well-known conjecture, to attribute this disproof to a particular person, and to blindly follow the judgement of a nice person one knows).
F. Frick's announcement [8]=\cite{Fr15o} obtaining $(A)$ from $(B)$ and the implication $(B)\Rightarrow(A)$ came in the situation when it was not known to the community that the implication is proved in
[9, 2.9.c]=\cite[2.9.c]{Gr10} (although this implication was implicitly rediscovered in [5]=\cite{BFZ14}, and although [9, 2.9e]=\cite[2.9.e]{Gr10} was discussed during the problem session at 2012 Oberwolfach Workshop on Triangulations).
Moreover, proof of $(A)$ after the Mabillard-Wagner announcement of $(B)$ was considered a serious problem by the community.
All this shows that F. Frick's contribution to the disproof of the topological Tverberg conjecture should be praised, but not in misleading form presented in \cite{BFZ}.
E.g. in \cite{Sk18u} Frick's name was listed in the phrase

`{\it For the counterexample papers \cite{Oz, Gr10, BFZ14, Fr15, MW15} by M.~\"Ozaydin,  M. Gromov, P. Blagojevi\'c, F. Frick, G. Ziegler, I. Mabillard and U. Wagner are important}',

although

$\bullet$ $(C)$ is attributed to \"Ozaydin (both in \cite{BFZ} and in the survey \cite{Sk18u} containing a published proof);

$\bullet$ $(C)\Rightarrow(B)$ is attributed to Mabillard-Wagner (both in \cite{BFZ} and in the survey \cite{Sk18u} containing a published proof);

$\bullet$ $(B)\Rightarrow(A)$ is proved by Gromov in [9, 2.9.c]=\cite[2.9.c]{Gr10}.

\smallskip
{\bf (ii).} Checking the argument of [14]=\cite{MW15} was an important concern for the community of topological combinatorics, and I see no indication that this checking is completed by the time the paper \cite{BFZ} is published (e.g. [14]=\cite{MW15} was still not accepted in 2019).
The paper [15]=\cite{MW14} announcing the results of [14]=\cite{MW15} is an `extended abstract'  according to [8]=\cite{Fr15o}.
As far as I know, the paper [17]=\cite{Oz} was rejected from a journal where it was submitted and not updated since then; for me the paper [17]=\cite{Oz} was useless when I tried to learn the proof of Theorems 3.4 (I learned the proof from a private communication with R. Karasev, who rediscovered the proof himself).

\smallskip
{\bf (iii). It is the more important to reveal the above drawbacks of \cite{BFZ} because support of the paper's unfair claim did not involve the authors' public answer to public criticism} presented in \cite[\S5]{Sk16-2},
(see also \cite[\S1.1]{JVZ} and \cite{Sk18u}), {\bf but did involve suppression of criticism both by not citing \cite{JVZ, Sk16-2, Sk18u} and by misuse of anonymous refereeing system} (see \cite[Remark 4.1]{Sk16-4} [Added in 2022: and Remark  \ref{r:hist}]).

\smallskip
Sincerely Yours, Arkadiy Skopenkov, \url{https://users.mccme.ru/skopenko/}
\end{remark}



{\it In this list books, surveys and expository papers are marked by stars}


\begin{thebibliography}{1}

\UseRawInputEncoding

\newcommand{\abc}{\bibitem[ABC+]{ABC+} * \emph{M. Atiyah, A. Borel, G. J. Chaitin, D. Friedan, J. Glimm, J. J. Gray, M. W. Hirsch, S. MacLane, B. B. Mandelbrot, D. Ruelle, A. Schwarz, K. Uhlenbeck, R. Thom, E. Witten, C.  Zeeman.} Responses to ``Theoretical Mathematics: Toward a cultural synthesis of mathematics and theoretical physics'', by A. Jaffe and F. Quinn. Bull. Am. Math. Soc. 30 (1994) 178--207. arXiv:math/9404229.}

\newcommand{\agles}{\bibitem[AGL]{AGL86} Mathematical Economics,  ed. by A. Ambrosetti, F. Gori, R. Lucchetti,
Lect. Notes Math. 1330, Springer, 1986.}


\newcommand{\akzz}{\bibitem[Ak00]{Ak00} * \emph{П. М. Ахметьев.} Вложения компактов, стабильные
гомотопические группы сфер и теория особенностей, Успехи Мат. Наук.  2000. 55:3. C.~3-62.}

\newcommand{\akoe}{\bibitem[AK19]{AK19} \emph{S. Avvakumov, R. Karasev.} Envy-free division using mapping degree. arXiv:1907.11183.}

\newcommand{\akto}{\bibitem[AK21]{AK21} \emph{G. Arone and V. Krushkal.}
Embedding obstructions in $\R^d$ from the Goodwillie-Weiss calculus and Whitney disks. arXiv:2101.10995. }

\newcommand{\akm}{\bibitem[AKM]{AKM} \emph{M. Abrahamsen, L. Kleist and T. Miltzow.}
Geometric Embeddability of Complexes is $\exists\mathbb R$-complete, arXiv:2108.02585.}

\newcommand{\aksoe}{\bibitem[AKS]{AKS} \emph{S. Avvakumov, R. Karasev and A. Skopenkov.} Stronger counterexamples to the topological Tverberg conjecture, submitted. arxiv:1908.08731.}

\newcommand{\akuoe}{\bibitem[AKu19]{AKu19} \emph{S. Avvakumov, S. Kudrya.}
Vanishing of all equivariant obstructions and the mapping degree. arXiv:1910.12628.}

\newcommand{\alto}{\bibitem[Al21]{Al21} \emph{E. Alkin,}
Hardness of almost embedding simplicial complexes in $\R^d$, II.}

\newcommand{\amsw}{\bibitem[AMS+]{AMSW} \emph{S. Avvakumov, I. Mabillard, A. Skopenkov and U. Wagner.}
Eliminating Higher-Multiplicity Intersections, III. Codimension 2, Israel J. Math. (2021).  arxiv:1511.03501.}


\newcommand{\anzt}{\bibitem[An03]{An03} * \emph{Д. В. Аносов.} Отображения окружности, векторные поля и их применения. М: МЦНМО, 2003.}

\newcommand{\arnf}{\bibitem[Ar95]{Ar95} \emph{V. I. Arnold,}  Topological invariants of plane curves and caustics, University Lecture Series, Vol. 5, Amer. Math. Soc., Providence, RI, 1995.}

\newcommand{\arszo}{\bibitem[ARS01]{ARS01} \emph{P. Akhmetiev, D. Repov\v s and A. Skopenkov},
Embedding products of low-dimensional manifolds in $\R^m$, Topol. Appl. 113 (2001), 7--12.}

\newcommand{\arszt}{\bibitem[ARS02]{ARS02} \emph{P. Akhmetiev, D. Repovs and A. Skopenkov.} Obstructions to approximating maps of $n$-manifolds into $R^{2n}$ by embeddings, Topol. Appl., 123 (2002), 3--14.}

\newcommand{\asoed}{\bibitem[As]{As} \emph{A. Asanau,} \lowercase{A SIMPLE PROOF THAT CONNECTED SUM OF ORDERED
ORIENTED LINKS IS NOT WELL-DEFINED,} Math. Notes, to appear.}

\newcommand{\asoe}{\bibitem[As]{As} \emph{A. Asanau,} On the \lowercase{TRIPLE SELF-INTERSECTION NUMBER FOR GRAPHS IN THE PLANE,} unpublished, 2018.}


\newcommand{\bbsn}{\bibitem[BB79]{BB} \emph{E.~G. Bajm{{\'o}}czy and I.~B{{\'a}}r{{\'a}}ny,}
\newblock On a common generalization of {B}orsuk's and {R}adon's theorem,
\newblock Acta Math.\ Acad.\ Sci.\ Hungar.\ 34:3 (1979), 347-350.}

\newcommand{\bbzos}{\bibitem[BBZ]{BBZ} * \emph{I.~B{{\'a}}r{{\'a}}ny, P.~V.~M. Blagojevi{{\'c}} and G.~M. Ziegler.} Tverberg's Theorem at 50: Extensions and Counterexamples, Notices of the Amer. Math. Soc., 63:7 (2016), 732--739.}


\newcommand{\bcm}{\bibitem[BCM]{BCM} * 13th Hilbert Problem on superpositions of functions, presented by A. Belov, A. Chilikov, I. Mitrofanov, S. Shaposhnikov and A. Skopenkov,
\url{http://www.turgor.ru/lktg/2016/5/index.htm}.}

\newcommand{\beet}{\bibitem[BE82]{BE82} * \emph{V.G. Boltyansky and V.A. Efremovich.} Intuitive Combinatorial Topology. Springer.}

\newcommand{\beetr}{\bibitem[BE82]{BE82} * \emph{В. Г. Болтянский и В. А. Ефремович.} Наглядная топология. М.:  Наука, 1982.}


\newcommand{\bfzof}{\bibitem[BFZ14]{BFZ14} \emph{P. V. M. Blagojevi{\'c}, F. Frick, and G. M. Ziegler,}
Tverberg plus constraints, Bull. Lond. Math. Soc. 46:5 (2014), 953-967, arXiv:1401.0690.}


\newcommand{\bfzos}{\bibitem[BFZ]{BFZ} \emph{P. V. M. Blagojevi{\'c}, F. Frick and G. M. Ziegler,}
Barycenters of Polytope Skeleta and Counterexamples to the Topological Tverberg Conjecture, via Constraints,
J. Eur. Math. Soc., 21:7 (2019) 2107-2116. arXiv:1510.07984.}


\newcommand{\bgos}{\bibitem[BG16]{BG16} \emph{A. Bj\"orner and A. Goodarzi}, On Codimension one Embedding of Simplicial Complexes, in book: A Journey Through Discrete Mathematics, arXiv:1605.01240.}

\newcommand{\biet}{\bibitem[Bi83]{Bi83} * \emph{R. H. Bing.} The Geometric Topology of 3-Manifolds. Providence, R.~I. 1983. (Amer. Math. Soc. Colloq. Publ., 40).}

\newcommand{\bitz}{\bibitem[Bi20]{Bi20} \emph{A. Bikeev.} Realizability of discs with ribbons on the M\"obius strip. arXiv:2010.15833.}


\newcommand{\bito}{\bibitem[Bi21]{Bi21} {\it A. I. Bikeev,}
Criteria for integer and modulo 2 embeddability of graphs to surfaces, arXiv:2012.12070v2.}



\newcommand{\bkkmzof}{\bibitem[BKK]{BKK} \emph{M. Bestvina, M. Kapovich and B. Kleiner,}
Van Kampen's embedding obstruction for discrete groups, Invent. Math. 150 (2002) 219--235. arXiv:math/0010141.}

\newcommand{\bmzf}{\bibitem[BM04]{BM04} \emph{Boyer, J. M. and Myrvold, W. J.} On the cutting edge: simplified $O(n)$ planarity by edge addition,  Journal of Graph Algorithms and Applications, 8:3 (2004) 241--273.}

\newcommand{\bm}{\bibitem[BM15]{BM15} \emph{I. Bogdanov and A. Matushkin.} Algebraic proofs of linear versions of the Conway--Gordon--Sachs theorem and the van Kampen--Flores theorem, arXiv:1508.03185.}


\newcommand{\bmzzn}{\bibitem[BMZ09]{BMZ09} \emph{P. V. M. Blagojevi{\'c}, B. Matschke, G. M. Ziegler,}
Optimal bounds for a colorful Tverberg-Vre\'cica type problem, Advances in Math., 226 (2011), 5198-5215, arXiv:0911.2692.}

\newcommand{\bmzof}{\bibitem[BMZ15]{BMZ15} \emph{P. V. M. Blagojevi{\'c}, B. Matschke, G. M. Ziegler,}
Optimal bounds for the colored Tverberg problem, J. Eur. Math. Soc.,  17:4 (2015) 739--754,
arXiv:0910.4987.}

\newcommand{\bpns}{\bibitem[BP97]{BP97} * \emph{R. Benedetti and C. Petronio.} Branched standard spines of 3-manifolds, Lecture Notes in Math. 1653, Springer-Verlag, Berlin-Heidelberg-New York, 1997.}

\newcommand{\brst}{\bibitem[Br72]{Br72} \emph{J. L. Bryant.} Approximating embeddings of polyhedra in codimension 3, Trans. Amer. Math. Soc., 170 (1972) 85--95.}

\newcommand{\brts}{\bibitem[Br26]{Br26} \emph{P. Bruegel.} The Magpie on the Gallows, 1526, \url{https://en.wikipedia.org/wiki/The_Magpie_on_the_Gallows}.}

\newcommand{\bren}{\bibitem[Br82]{brown1982} * \emph{K.~S. Brown.} \newblock Cohomology of Groups. \newblock Springer-Verlag New York, 1982.}


\newcommand{\bssos}{\bibitem[BS17]{BS17} * \emph{I.~B\'{a}r\'{a}ny and P. Sober\'{o}n,} Tverberg's theorem is 50 years old: a survey, arXiv:1712.06119.}


\newcommand{\brsnn}{\bibitem[BRS99]{BRS99} \emph{D. Repov\v s, N. Brodsky and A. B. Skopenkov.}
A classification of 3-thickenings of 2-polyhedra, Topol. Appl. 1999. 94. P.~307-314.}

\newcommand{\bsseo}{\bibitem[BSS]{BSS} \emph{I.~B\'{a}r\'{a}ny, S.~B. Shlosman, and A.~Sz{\H{u}}cs,}
\newblock On a topological generalization of a theorem of {T}verberg,
\newblock J.\ London Math.\ Soc.\ (II. Ser.) 23 (1981), 158--164.}

\newcommand{\btzs}{\bibitem[BT07]{BT07} \emph{A. Bj\"orner, M. Tancer}, Combinatorial Alexander Duality --- a Short and Elementary Proof, Discr. and Comp. Geom., 42 (2009) 586. arXiv:0710.1172.}

\newcommand{\buse}{\bibitem[Bu68]{Bu68} \emph{A. R. Butz,} Space filling curves and mathematical programming, Information and Control, 12:4 (1968) 314--330.}


\newcommand{\bz}{\bibitem[BZ16]{BZ16} * \emph{P. V. M. Blagojevi\'c and G. M. Ziegler,} Beyond the Borsuk-Ulam theorem: The topological Tverberg story, in: A Journey Through Discrete Mathematics, Eds. M. Loebl,
J. Ne\v set\v ril, R. Thomas, Springer, 2017, 273--341. arXiv:1605.07321v3.}



\newcommand{\carm}{\bibitem[Ca]{Ca} \emph{J. Carmesin.}
Embedding simply connected 2-complexes in 3-space, I-V, arXiv:1709.04642, arXiv:1709.04643, arXiv:1709.04645, arXiv:1709.04652, arXiv:1709.04659.}

\newcommand{\cfsz}{\bibitem[CF60]{CF60} \emph{P. E. Conner and E. E. Floyd}, Fixed points free involutions and equivariant maps, Bull. Amer. Math. Soc., 66 (1960) 416--441.}

\newcommand{\cget}{\bibitem[CG83]{CG83} \emph{J. H. Conway and C. M. A. Gordon},
Knots and links in spatial graphs, J. Graph Theory  7 (1983), 445--453.}

\newcommand{\cten}{\bibitem[Ch]{Ch} \emph{Chuang Tzu,} translated by H. A. Giles, Bernard Quaritch, London, 1889.}

\newcommand{\ctruku}{\bibitem[Ch]{Ch} \emph{Chuang Tzu,} translated to Russian by S. Kuchera, in: Ancient Chinese Philosophy, v. I, Mysl, Moscow, 1972.}


\newcommand{\chnn}{\bibitem[Ch99]{Ch99} * \emph{А. В. Чернавский,} Теорема Жордана.  Мат. Просвещение, 3 (1999), 142--157.}

\newcommand{\hcon}{\bibitem[HC19]{HC19} * \emph{C. Herbert Clemens.} Two-Dimensional Geometries. A Problem-Solving Approach, Amer. Math. Soc., 2019.}

\newcommand{\ckmoo}{\bibitem[CKMS]{CKMS} \emph{M. \v Cadek, M. Kr\v c\'al. J. Matou\v sek, F. Sergeraert,
L. Vok\v r\'inek, U. Wagner.} Computing all maps into a sphere, J. of the ACM, 61:3 (2014). arXiv:1105.6257.}


\newcommand{\ckmvwot}{\bibitem[CKM12+]{CKM12+} \emph{M. \v Cadek, M. Kr\v c\'al. J. Matou\v sek, L. Vok\v r\'inek, U. Wagner.} Polynomial-time computation of homotopy groups and Postnikov systems in fixed dimension, SIAM J. Comput., 43:5 (2014), 1728--1780. arXiv:1211.3093.}

\newcommand{\ckmvw}{\bibitem[CKM+]{CKM+} \emph{M. \v Cadek, M. Kr\v c\'al. J. Matou\v sek, L. Vok\v r\'inek, U. Wagner.} Extendability of continuous maps is undecidable, Discr. and Comp. Geom. 51 (2014) 24--66.
arXiv:1302.2370.}

\newcommand{\ckppt}{\bibitem[CKP+]{CKP+} \emph{E. Colin de Verdi\'ere, V. Kalu\v za, P. Pat\'ak, Z. Pat\'akov\'a and M. Tancer.} A direct proof of the strong Hanani-Tutte theorem on the projective plane. Journal of Graph Algorithms and Applications, 21:5 (2017) 939--981.}

\newcommand{\cksof}{\bibitem[CKS+]{CKS+} * New ways of weaving baskets, presented by G. Chelnokov, Yu. Kudryashov, A.Skopenkov and A. Sossinsky, \url{http://www.turgor.ru/lktg/2004/lines.en/index.htm}.}

\newcommand{\ckv}{\bibitem[CKV]{CKV} \emph{M.~{\v{C}}adek, M.~Kr\v{c}\'{a}l, and L.~Vok\v{r}\'{\i}nek.}
Algorithmic solvability of the lifting-extension problem, Discr. Comp. Geom. 57 (2017), 915--965. arXiv:1307.6444.}


\newcommand{\clr}{\bibitem[CLR]{CLR} * \emph{Т. Кормен, Ч. Лейзерсон, Р. Ривест.} Алгоритмы:
построение и анализ, МЦНМО, Москва, 1999.}

\newcommand{\clreng}{\bibitem[CLR]{CLR} * \emph{T. H. Cormen, C. E.Leiserson, R. L.Rivest, C. Stein.} Introduction to Algorithms, MIT Press, 2009.}

\newcommand{\crzfru}{\bibitem[CR]{CR} * \emph{Р. Курант, Дж. Роббинс,} Что такое математика. М.: МЦНМО, 2004.}

\newcommand{\crzfen}{\bibitem[CR]{CR} * \emph{R. Courant and H. Robbins,} What is Mathematics, Oxford Univ. Press.}

\newcommand{\crsne}{\bibitem[CRS98]{CRS98} * \emph{A. Cavicchioli, D. Repov\v s and A. B. Skopenkov.}
Open problems on graphs, arising from geometric topology, Topol. Appl. 1998. 84. P.~207-226.}

\newcommand{\crsot}{\bibitem[CRS]{CRS} \emph{M. Cencelj, D. Repov\v s and M. Skopenkov,}
Classification of knotted tori in the 2-metastable dimension, Mat. Sbornik, 203:11 (2012), 1654--1681.
arxiv:math/0811.2745.}

\newcommand{\csoo}{\bibitem[CS08]{CS08} \emph{D. Crowley and A. Skopenkov.} A classification of smooth embeddings of 4-manifolds in 7-space, II, Intern. J. Math., 22:6 (2011) 731-757, arxiv:math/0808.1795.}

\newcommand{\csos}{\bibitem[CS16]{CS16} \emph{D. Crowley and A. Skopenkov,} Embeddings of non-simply-connected 4-manifolds in 7-space. I. Classification modulo knots, Moscow Math. J., 21 (2021), 43--98. arXiv:1611.04738.}

\newcommand{\csoso}{\bibitem[CS16o]{CS16o} \emph{D. Crowley and A. Skopenkov,} Embeddings of non-simply-connected 4-manifolds in 7-space. II. On the smooth classification, Proc. A of the Royal Soc. of Edinburgh, to appear. arXiv:1612.04776.}


\newcommand{\crsk}{\bibitem[CS]{CS} \emph{D. Crowley and A. Skopenkov,} Embeddings of non-simply-connected 4-manifolds in 7-space. III. Piecewise-linear classification. draft.}

\newcommand{\cutz}{\bibitem[Cu20]{Cu20} \emph{C. Culter,} Cantor sets are not tangent homogeneous,
Topol. Appl. 271 (2020) 1--9.}


\newcommand{\dies}{\bibitem[Di87]{Di} * \emph{T. tom Dieck,} Transformation groups, Studies in Mathematics, vol. 8, Walter de Gruyter, Berlin, 1987.}

\newcommand{\dent}{\bibitem[De93]{De93}  \emph{T.K. Dey.} On counting triangulations in $d$-dimensions. Comput. Geom.  3:6 (1993) 315--325.}

\newcommand{\denf}{\bibitem[DE94]{DE94}  \emph{T.K. Dey and H. Edelsbrunner.} Counting triangle crossings and halving planes, Discrete Comput. Geom. 12 (1994), 281--289.}


\newcommand{\ers}{\bibitem[ERS]{ERS} * Invariants of graph drawings in the plane, presented by A. Enne, A. Ryabichev, A. Skopenkov and T. Zaitsev, \url{http://www.turgor.ru/lktg/2017/6/index.htm}}


\newcommand{\ffen}{\bibitem[FF89]{FF89} * \emph{А. Т. Фоменко и Д. Б. Фукс.} Курс гомотопической топологии. М.: Наука, 1989.}

\newcommand{\ffene}{\bibitem[FF89]{FF89} * \emph{A.T. Fomenko and D.B. Fuchs.} Homotopical Topology, Springer, 2016.}


\newcommand{\fkosc}{\bibitem[FK17]{FK17} \emph{R. Fulek, J. Kyn{\v{c}}l,} Counterexample to an Extension of the Hanani-Tutte Theorem on the Surface of Genus 4, Combinatorica, 39 (2019) 1267--1279, arXiv:1709.00508.}

\newcommand{\fkos}{\bibitem[FK17]{FK17} \emph{R. Fulek, J. Kyn{\v{c}}l,} Hanani-Tutte for approximating maps of graphs, arXiv:1705.05243.}

\newcommand{\fkon}{\bibitem[FK19]{FK19} \emph{R. Fulek, J. Kyn{\v{c}}l,}
$\Z_2$-genus of graphs and minimum rank of partial symmetric matrices,
35th Intern. Symp. on Comp. Geom. (SoCG 2019), Article No. 39; pp. 39:1–39:16,
\url{https://drops.dagstuhl.de/opus/volltexte/2019/10443/pdf/LIPIcs-SoCG-2019-39.pdf}.
We refer to numbering in arXiv version: arXiv:1903.08637.}

\newcommand{\fktnf}{\bibitem[FKT]{FKT} \emph{M. H. Freedman, V. S. Krushkal and P. Teichner.} Van Kampen's
embedding obstruction is incomplete for 2-complexes in~$\R^4$, Math. Res. Letters. 1994. 1. P.~167-176.}

\newcommand{\fltf}{\bibitem[Fl34]{Fl34} \emph{A. Flores}, \"Uber $n$-dimensionale Komplexe die im $E^{2n+1}$ absolut selbstverschlungen sind, Ergeb. Math. Koll. 6 (1934) 4--7.}

\newcommand{\fo}{\bibitem[Fo]{Fo} * \emph{L. Fortnow.} Time for Computer Science to Grow Up,  \url{https://people.cs.uchicago.edu/~fortnow/papers/growup.pdf}.}

\newcommand{\fozf}{\bibitem[Fo04]{Fo04} * \emph{R. Fokkink.} A forgotten mathematician, Eur. Math. Soc. Newsletter 52 (2004) 9--14.}


\newcommand{\fpstz}{\bibitem[FPS]{FPS} \emph{R. Fulek, M.J. Pelsmajer and M. Schaefer.}
Strong Hanani-Tutte for the Torus, arXiv:2009.01683.}

\newcommand{\frse}{\bibitem[Fr78]{Fr78} \emph{M. Freedman,} Quadruple points of 3-manifolds in $S^4$, Comment. Math. Helv. 53 (1978), 385-394.}

\newcommand{\fres}{\bibitem[FR86]{FR86} \emph{R. Fenn, D. Rolfsen.}
Spheres may link homotopically in 4-space, J. London Math. Soc. 34 (1986) 177-184.}

\newcommand{\frof}{\bibitem[Fr15]{Fr15} \emph{F. Frick}, Counterexamples to the topological Tverberg conjecture,
Oberwolfach reports, 12:1 (2015), 318--321. arXiv:1502.00947.}


\newcommand{\fros}{\bibitem[Fr17]{Fr17} \emph{F. Frick}, O\lowercase{N AFFINE TVERBERG-TYPE RESULTS WITHOUT CONTINUOUS GENERALIZATION}, arXiv:1702.05466}


\newcommand{\fstz}{\bibitem[FS20]{FS20} \emph{F. Frick and P. Sober\'on}, The topological Tverberg problem beyond prime powers, arXiv:2005.05251.}

\newcommand{\fvto}{\bibitem[FV21]{FV21} \emph{M. Filakovsk\'y, L. Vok\v r\'inek.} Computing homotopy classes for diagrams, 	arXiv:2104.10152.}

\newcommand{\fwz}{\bibitem[FWZ]{FWZ} \emph{M. Filakovsk\'y, U. Wagner, S. Zhechev.} Embeddability of simplicial complexes is undecidable. Oberwolfach reports, to appear.}

\newcommand{\fwztz}{\bibitem[FWZ]{FWZ} \emph{M. Filakovsk\'y, U. Wagner, S. Zhechev.} Embeddability of simplicial complexes is undecidable.
Proceedings of the 2020 ACM-SIAM Symposium on Discrete Algorithms, \url{https://epubs.siam.org/doi/pdf/10.1137/1.9781611975994.47}}


\newcommand{\ga}{\bibitem[GA]{GA} * \url{https://en.wikipedia.org/wiki/Galactic_algorithm}}

\newcommand{\gdikrse}{\bibitem[GDI]{GDI} * {\it A. Chernov, A. Daynyak, A. Glibichuk, M. Ilyinskiy, A. Kupavskiy, A. Raigorodskiy and A. Skopenkov,} Elements of Discrete Mathematics As a Sequence of Problems (in Russian),
MCCME, Moscow, 2016. Update: \url{http://www.mccme.ru/circles/oim/discrbook.pdf} .}

\newcommand{\gdikrs}{\bibitem[GDI]{GDI} * {\it А.А. Глибичук, А.Б. Дайняк, Д.Г. Ильинский, А.Б. Купавский, А.М. Райгородский, А.Б. Скопенков, А.А. Чернов,} Элементы дискретной математики в задачах, М, МЦНМО, 2016.
\url{http://www.mccme.ru/circles/oim/discrbook.pdf} .}

\newcommand{\giso}{\bibitem[Gi71]{Gi71} * {\it S. Gitler,} Immersion and Embedding of Manifolds,
Proc. Symp. Pure Math. 22, 87-96 (1971).}

\newcommand{\gkp}{\bibitem[GKP]{GKP} * {\it R. Graham, D. Knuth, and O. Patashnik,} Concrete Mathematics: A Foundation for Computer Science, Addison–Wesley, first published in 1989, \url{https://www.csie.ntu.edu.tw/~r97002/temp/Concrete\%20Mathematics\%202e.pdf}.}

\newcommand{\gmpptw}{\bibitem[GMP+]{GMP+} \emph{X. Goaoc, I. Mabillard, P. Pat\'ak, Z. Pat\'akov\'a, M. Tancer, U. Wagner}, On Generalized Heawood Inequalities for Manifolds: a van Kampen--Flores-type Nonembeddability Result,
arXiv:1610.09063.}

\newcommand{\grsz}{\bibitem[Gr69]{Gr69} \emph{B. Gr\"unbaum.} Imbeddings of simplicial complexes. Comment. Math. Helv., 44:1, 502--513, 1969.}


\newcommand{\gres}{\bibitem[Gr86]{Gr86} * \emph{M. Gromov}, Partial Differential Relations,
Ergebnisse der Mathematik und ihrer Grenzgebiete (3), Springer Verlag, Berlin-New York, 1986.}

\newcommand{\groz}{\bibitem[Gr10]{Gr10} \emph{M. Gromov,}
\newblock Singularities, expanders and topology of maps. Part 2: From combinatorics to topology via algebraic isoperimetry, \newblock Geometric and Functional Analysis 20 (2010), no.~2, 416--526.}

\newcommand{\grsn}{\bibitem[GR79]{GR79} \emph{J. L. Gross	and R. H. Rosen}, A linear time planarity algorithm for 2-complexes, Journal of the ACM, 26:4 (1979), 611--617.}

\newcommand{\gs}{\bibitem[GS]{GS} \emph{М. Гортинский и О. Скрябин.} Критерий вложимости графов в плоскость вдоль прямой, препринт.}

\newcommand{\gssn}{\bibitem[GS79]{GS} \emph{P.~M. Gruber and R.~Schneider,} Problems in geometric convexity. In {\em Contributions to geometry ({P}roc. {G}eom. {S}ympos., {S}iegen, 1978)}, 255--278. Birkh{\"a}user, Basel-Boston, Mass., 1979.}

\newcommand{\gsnn}{\bibitem[GS99]{GS99} \emph{R. Gompf and A. Stipsicz,}
4-manifolds and Kirby calculus, GSM20, AMS, Providence, RI, 1999.}


\newcommand{\gszs}{\bibitem[GS06]{GS06} \emph{D. Goncalves and A. Skopenkov,} Embeddings of homology equivalent manifolds with boundary, Topol. Appl., 153:12 (2006) 2026-2034. arxiv:1207.1326.}

\newcommand{\gssoe}{\bibitem[GSS+]{GSS+} * Projections of skew lines, presented by A. Gaifullin, A. Shapovalov, A. Skopenkov and M. Skopenkov, \url{http://www.turgor.ru/lktg/2001/index.php}.}

\newcommand{\gtes}{\bibitem[GT87]{GT87} * \emph{J. L. Gross and T. W. Tucker.}
Topological graph theory. New York: Wiley-Interscience, 1987.}

\newcommand{\guzn}{\bibitem[Gu09]{Gu09} \emph{A. Gundert.} On the complexity of embeddable simplicial complexes. Diplomarbeit, Freie Universit\"at Berlin, 2009. 	arXiv:1812.08447.}


\newcommand{\ha}{\bibitem[Ha]{Ha} * \emph{F. Harary.} Graph theory.
Рус. пер.: Ф. Харари. Теория графов. М., Мир, 1973.}

\newcommand{\hats}{\bibitem[Ha37]{Ha37} \emph{W. Hantzsche,} Einlagerung von Mannigfaltigkeiten in euklidische R\" aume, Math. Zeitschrift, 43:1 (1937) 38--58.}

\newcommand{\hastk}{\bibitem[Ha62k]{Ha62k} {\em A.~Haefliger,}  Knotted $(4k-1)$-spheres in $6k$-space, Ann. of Math. 75 (1962) 452--466.}

\newcommand{\hastl}{\bibitem[Ha62l]{Ha62l} \emph{A. Haefliger,} Differentiable links, Topology, 1 (1962) 241--244.}

\newcommand{\hast}{\bibitem[Ha63]{Ha63} \emph{A.~Haefliger,} Plongements differentiables dans le domain stable, Comment. Math. Helv. 36 (1962-63) 155--176.}

\newcommand{\hassa}{\bibitem[Ha66A]{Ha66A} \textit{A. Haefliger}. Differential embeddings of~$S^n$ in $S^{n+q}$ for $q>2$. Ann. Math. (2), 83 (1966), 402--~436.}

\newcommand{\hass}{\bibitem[Ha66C]{Ha66C} \emph{A.~Haefliger,}  Enlacements de spheres en codimension superiure a 2, Comment. Math. Helv. 41 (1966-67) 51--72.}

\newcommand{\hase}{\bibitem[Ha68]{Ha68} \emph{A. Haefliger,} Knotted Spheres and Related Geometric Topic,
in Proc. Int. Congr. Math., Moscow, 1966 (Mir, Moscow, 1968), 437--445.}

\newcommand{\hasn}{\bibitem[Ha69]{Ha69} \emph{L.~S.~Harris,} Intersections and embeddings of polyhedra, Topology 8 (1969) 1--26.}

\newcommand{\hasf}{\bibitem[Ha74]{Ha74} * \emph{P. Halmos,} How to talk mathematics. Notices of the Amer. Math. Soc., 21 (1974) 155--158.}

\newcommand{\haef}{\bibitem[Ha84]{Ha84} \emph{N. Habegger,} Obstruction to embedding disks II: a proof of a conjecture by Hudson, Topol. Appl. 17 (1984).}

\newcommand{\haes}{\bibitem[Ha86]{Ha86} \emph{N. Habegger,} Knots and links in codimension greater than 2, Topology, 25:3 (1986) 253--260.}

\newcommand{\hifn}{\bibitem[Hi59]{Hi59} \emph{M. W. Hirsch.} Immersions of manifolds, Trans. Amer. Math. Soc. 93 (1959) 242--276.}

\newcommand{\hjsf}{\bibitem[HJ64]{HJ64} \emph{R. Halin and H. A. Jung.}
Karakterisierung der Komplexe der Ebene und der 2-Sph\"are, Arch. Math. 1964. 15. P.~466-469.}

\newcommand{\hkne}{\bibitem[HK98]{HK98} \emph{N. Habegger and U. Kaiser,} Link homotopy in 2--metastable range, Topology 37:1 (1998) 75--94.}

\newcommand{\hmsnt}{\bibitem[HMS]{HMS93} * \emph{C. Hog-Angeloni, W. Metzler and A. J. Sieradski.}
Two-dimensional homotopy and combinatorial group theory. Cambridge: Cambridge Univ. Press, 1993. (London Math. Soc. Lecture Notes, 197).}

\newcommand{\ho}{\bibitem[Ho]{Ho} * The Hopf fibration, \url{https://www.youtube.com/watch?v=AKotMPGFJYk}}

\newcommand{\hozs}{\bibitem[Ho06]{Ho06} \emph{H. van der Holst,} Graphs and obstructions in four dimensions, J. Combin. Theory Ser. B 96:3 (2006), 388--404.}


\newcommand{\hpzn}{\bibitem[HP09]{HP09} \emph{H. van der Holst and R. Pendavingh,} On a graph property generalizing planarity and flatness, Combinatorica, 29 (2009) 337--361.}

\newcommand{\htsf}{\bibitem[HT74]{HT74} \emph{J. Hopcroft and R. E. Tarjan,} Efficient planarity testing, J. of the Association for Computing Machinery, 21:4 (1974) 549--568.}

\newcommand{\hufn}{\bibitem[Hu59]{hu59} * \emph{S. T. Hu,} Homotopy Theory, Academic Press, New York, 1959.}

\newcommand{\husn}{\bibitem[Hu69]{Hu69} * \emph{J. F. P. Hudson.} Piecewise linear topology, W. A. Benjamin, Inc., New York-Amsterdam, 1969.}


\newcommand{\io}{\bibitem[Io]{Io} * \url{https://en.wikipedia.org/wiki/Category:Impossible_objects}}

\newcommand{\info}{\bibitem[IF]{IF} * \url{http://www.map.mpim-bonn.mpg.de/Intersection_form}}

\newcommand{\irsf}{\bibitem[Ir65]{Ir65} \emph{M.~C.~Irwin,} Embeddings of polyhedral manifolds, Ann. of Math. (2)
82 (1965) 1--14.}

\newcommand{\isot}{\bibitem[Is]{Is} * \url{http://www.map.mpim-bonn.mpg.de/Isotopy}}


\newcommand{\jqnt}{\bibitem[JQ93]{JQ93} * \emph{A. Jaffe, F. Quinn,} ``Theoretical mathematics'': Toward a cultural synthesis of mathematics and theoretical physics. Bull.Am.Math.Soc. 29 (1993) 1-13. arXiv:math/9307227.}

\newcommand{\jozt}{\bibitem[Jo02]{Jo02} \emph{C. M. Johnson.} An obstruction to embedding a simplicial $n$-complex into a $2n$-manifold, Topology Appl. 122:3 (2002) 581--591.}

\newcommand{\jvz}{\bibitem[JVZ]{JVZ} D. Joji\'c, S. T. Vre\'cica, R. T. \v Zivaljevi\' c,
Topology and combinatorics of 'unavoidable complexes', arXiv:1603.08472v1.}


\newcommand{\kalai}{\bibitem[Ka]{Ka} G. Kalai, From Oberwolfach: The Topological Tverberg Conjecture is False, `Combinatorics and more' blog post, February 6, 2015, \url{gilkalai.wordpress.com}}

\newcommand{\kh}{\bibitem[Kh]{Kh} \emph{А.И. Храбров.} Руководство по чтению лекций
\url{http://vm.tstu.tver.ru/topics/pdf_tests/lection.pdf}}

\newcommand{\kho}{\bibitem[Kho]{Kho} \emph{N. Khoroshavkina.} A simple characterization of graphs of cutwidth 2, arXiv:1811.06716.}

\newcommand{\kkrot}{\bibitem[KKR]{KKR} \emph{K. Kawarabayashi, Y. Kobayashi and B. Reed.} The disjoint paths problem in quadratic time, J. of Comb. Theory, Ser. B, 102:2 (2012), 424--435.}

\newcommand{\kmsth}{\bibitem[KM63]{KM63} \emph{M. A. Kervaire and J. W. Milnor,} Groups of homotopy spheres. I,  Ann. of Math. (2) 77 (1963), 504-537.}

\newcommand{\kozeru}{\bibitem[Ko18]{Ko18} * \emph{Е. Колпаков.}
Доказательство теоремы Радона при помощи понижения размерности, Мат. Просвещение, 23 (2018), arXiv:1903.11055.}

\newcommand{\koze}{\bibitem[Ko18]{Ko18} * \emph{E. Kolpakov.}
A proof of Radon Theorem via lowering of dimension, Mat. Prosveschenie, 23 (2018), arXiv:1903.11055.}

\newcommand{\ko}{\bibitem[Ko]{Ko} \emph{E. Kolpakov.} A `converse' to the Constraint Lemma, arXiv:1903.08910.}

\newcommand{\koon}{\bibitem[Ko19]{Ko19} \emph{E. Kogan.} Linking of three triangles in 3-space, arXiv:1908.03865.}

\newcommand{\koto}{\bibitem[Ko21]{Ko21} \emph{E. Kogan.} On the rank of $\Z_2$-matrices with free entries on the diagonal, arXiv:2104.10668.}

\newcommand{\koee}{\bibitem[Ko88]{Ko88} \emph{U. Koschorke.} Link maps and the geometry of their invariants,
Manuscripta Math. 61:4 (1988) 383--415.}

\newcommand{\kps}{\bibitem[KPS]{KPS} * \emph{A. Kaibkhanov, D. Permyakov and A. Skopenkov.}
Realization of graphs with rotation, \url{http://www.turgor.ru/lktg/2005/3/index.htm}.}

\newcommand{\krzz}{\bibitem[Kr00]{Kr00} \emph{V. S. Krushkal.} Embedding obstructions and 4-dimensional thickenings of 2-complexes, Proc. Amer. Math. Soc. 128:12 (2000) 3683--3691. arXiv:math/0004058. }

\newcommand{\ksnn}{\bibitem[KS99]{KS99} * \emph{П. Кожевников и А. Скопенков.} Узкие деревья на плоскости, Мат. Образование. 1999. 2-3. С.~126-131.}

\newcommand{\kstz}{\bibitem[KS20]{KS20} \emph{R. Karasev and A. Skopenkov.}
Some `converses' to intrinsic linking theorems, arXiv:2008.02523.}

\newcommand{\ksto}{\bibitem[KS21]{KS21} \emph{E. Kogan and A. Skopenkov.} A short exposition of the Patak-Tancer theorem on non-embeddability of $k$-complexes in $2k$-manifolds,  arXiv:2106.14010.}

\newcommand{\kstoe}{\bibitem[KS21e]{KS21e} \emph{E. Kogan and A. Skopenkov.}
Embeddings of $k$-complexes in $2k$-manifolds and minimum rank of partial symmetric matrices, arXiv:2112.06636.}

\newcommand{\kuse}{\bibitem[Ku68]{Ku68} * \emph{К. Куратовский.} Топология. Т.~1,~2. М.: Мир, 1969.}


\newcommand{\lazz}{\bibitem[La00]{La00} \emph{F. Lasheras.} An obstruction to 3-dimensional thickening,
Proc. Amer. Math. Soc. 2000. 128. P.~893-902.}

\newcommand{\lfma}{\bibitem[LF]{LF} \url{http://www.map.mpim-bonn.mpg.de/Linking_form}}

\newcommand{\lloe}{\bibitem[LL18]{LL18} \emph{A.S. Levine and T. Lidman.} Simply connected, spineless 4-manifolds, Forum of Math., Sigma, 7 (2019) e14, 1--11, arxiv:1803.01765.}

\newcommand{\lo}{\bibitem[Lo]{Lo} M.~de~Longueville. Notes on the topological Tverberg theorem.
Discrete Math.  247 (2002), no.~1--3, 271--297.
(The paper first appeared in
Discrete Math. 241 (2001) 207--233, but the original version suffered from serious publisher's typesetting errors.)}

\newcommand{\loot}{\bibitem[Lo13]{Lo13} \emph{M. de Longueville.} A course in topological combinatorics. Universitext. Springer, New York (2013).}

\newcommand{\lssn}{\bibitem[LS69]{LS69} \emph{W. B. R. Lickorish and L. C. Siebenmann.}
Regular neighborhoods and the stable range,  Trans. Amer. Math. Soc.. 1969. 139. P.~207-230.}

\newcommand{\lsne}{\bibitem[LS98]{LS98} \emph{L. Lovasz and A. Schrijver,}
A Borsuk theorem for antipodal links and a spectral characterization of linklessly embeddable graphs, Proc. Amer. Math. Soc. 126:5 (1998), 1275-1285.}

\newcommand{\ltof}{\bibitem[LT14]{LT14} \emph{E. Lindenstrauss and M. Tsukamoto,} Mean dimension and an embedding problem: an example, Israel J. Math. 199 (2014).}


\newcommand{\lyzf}{\bibitem[LY04]{LY04} * \emph{Y. Lin and A. Yang,} On 3-cutwidth critical graphs, Discrete Mathematics, 275 (2004), 339--346.}

\newcommand{\lz}{\bibitem[LZ]{LZ} * \emph{S. Lando and A. Zvonkin.} Embedded Graphs. Springer.}



\newcommand{\mast}{\bibitem[Ma73]{Ma73} \emph{С. В. Матвеев.} Специальные остовы кусочно-линейных многообразий, Мат. Сборник. 1973. 92. С.~282-293.}

\newcommand{\maste}{\bibitem[Ma73]{Ma73} \emph{S. V. Matveev.} Special skeletons of PL manifolds (in Russian), Mat. Sbornik. 1973. 92. P.~282-293.}

\newcommand{\mans}{\bibitem[Ma97]{Ma97} \emph{Yu. Makarychev.} A short proof of Kuratowski's graph planarity criterion, J. of Graph Theory, 25 (1997), 129--131.}

\newcommand{\mazt}{\bibitem[Ma03]{Ma03} * \emph{J.~Matou{\v{s}}ek.} Using the {B}orsuk-{U}lam theorem:
Lectures on topological methods in combinatorics and geometry. Springer Verlag, 2008.}


\newcommand{\mazf}{\bibitem[Ma05]{Ma05} \emph{V. Manturov.} A proof of the Vasiliev conjecture on the planarity of singular links, Izv. RAN 2005.}

\newcommand{\metn}{\bibitem[Me29]{Me29} \emph{K. Menger.} \"Uber pl\"attbare Dreiergraphen und Potenzen nicht pl\"attbarer Graphen, Ergebnisse Math. Kolloq., 2 (1929) 30--31.}

\newcommand{\mezf}{\bibitem[Me04]{Me04} \emph{S. Melikhov.} Sphere eversions and realization of mappings, Trudy MIAN 247 (2004) 159-181 (in Russian) arXiv:math.GT/0305158.}

\newcommand{\mezs}{\bibitem[Me06]{Me06} \emph{S. A. Melikhov}, The van Kampen obstruction and its relatives, 	
Proc. Steklov Inst. Math 266 (2009), 142-176 (= Trudy MIAN 266 (2009), 149-183), arXiv:math/0612082.}

\newcommand{\meoo}{\bibitem[Me11]{Me11} \emph{S. A. Melikhov}, Combinatorics of embeddings, arXiv:1103.5457.}

\newcommand{\meos}{\bibitem[Me17]{Me17} \emph{S. Melikhov,} Gauss type formulas for link map invariants, arXiv:1711.03530.}

\newcommand{\meoe}{\bibitem[Me18]{Me18} \emph{S. A. Melikhov,} A triple-point Whitney trick, J. Topol. Anal., 2018, 1--6.}


\newcommand{\miso}{\bibitem[Mi61]{Mi61} \emph{J. Milnor,} A procedure for killing homotopy groups of differentiable manifolds, Proc. Sympos. Pure Math, Vol. III (1961), 39--55.}

\newcommand{\mins}{\bibitem[Mi97]{Mi97} \emph{P. Minc.} Embedding simplicial arcs into the plane, Topol. Proc. 1997. 22. 305--340.}


\newcommand{\moss}{\bibitem[Mo77]{Mo77} * \emph{E. E. Moise.} Geometric Topology in Dimensions 2 and 3 (GTM), Springer-Verlag, 1977.}

\newcommand{\moen}{\bibitem[Mo89]{Mo89} \textit{B. Mohar}. An obstruction to embedding graphs in
surfaces. Discrete Math. 78 (1989) 135--142.}

\newcommand{\mrst}{\bibitem[MRS+]{MRS+} \emph{A. de Mesmay, Y. Rieck, E. Sedgwick, M. Tancer,}
Embeddability in $\R^3$ is NP-hard. arXiv:1708.07734.}

\newcommand{\mesczs}{\bibitem[MS06]{MS06} \emph{S.A. Melikhov, E.V. Shchepin,} The telescope approach to embeddability of compacta. arXiv:math.GT/0612085.}

\newcommand{\mstwof}{\bibitem[MST+]{MST+} \emph{J. Matou\v sek, E. Sedgwick, M. Tancer, U. Wagner}, Embeddability in the 3-sphere is decidable, Journal of the ACM 65:1 (2018) 1--49, arXiv:1402.0815.}


\newcommand{\mtzo}{\bibitem[MT01]{MT01} * \emph{B. Mohar and C. Thomassen.} Graphs on Surfaces.
The John Hopkins University Press, 2001.}

\newcommand{\mtwoz}{\bibitem[MTW10]{MTW10} \emph{J. Matou\v sek, M. Tancer, U. Wagner.} A geometric proof of
the colored Tverberg theorem, Discr. and Comp. Geometry, 47:2 (2012), 245--265. arXiv:1008.5275.}


\newcommand{\mtwoo}{\bibitem[MTW]{MTW} \emph{J. Matou\v sek, M. Tancer, U. Wagner.}
Hardness of embedding simplicial complexes in $\R^d$, J. Eur. Math. Soc. 13:2 (2011), 259--295. arXiv:0807.0336.}


\newcommand{\mwofo}{\bibitem[MW14]{MW14} \emph{I. Mabillard and U. Wagner.} Eliminating Tverberg Points, I. An Analogue of the Whitney Trick, Proc. of the 30th Annual Symp. on Comp. Geom. (SoCG'14), ACM, New York, 2014, pp. 171--180.}

\newcommand{\mwof}{\bibitem[MW15]{MW15} \emph{I. Mabillard and U. Wagner.}
Eliminating Higher-Multiplicity Intersections, I. A Whitney Trick for Tverberg-Type Problems. arXiv:1508.02349.}


\newcommand{\mwos}{\bibitem[MW16]{MW16} \emph{I. Mabillard and U. Wagner.} Eliminating Higher-Multiplicity Intersections, II. The Deleted Product Criterion in the $r$-Metastable Range. arXiv:1601.00876v2.}

\newcommand{\mwosd}{\bibitem[MW16']{MW16'} \emph{I. Mabillard and U. Wagner.} Eliminating Higher-Multiplicity Intersections, II. The Deleted Product Criterion in the r-Metastable Range,
Proceedings of the 32nd Annual Symposium on Computational Geometry (SoCG'16).}


\newcommand{\neno}{\bibitem[Ne91]{Ne91} \emph{S. Negami.} Ramsey theorems for knots, links and spatial graphs,
Trans. Amer. Math. Soc., 324 (1991), 527--541.}



\newcommand{\nkon}{\bibitem[NKS]{NKS} * \emph{L. T. Nguyen, J. Kim, B. Shim.}
Low-Rank Matrix Completion: A Contemporary Survey. arXiv:1907.11705.}

\newcommand{\noss}{\bibitem[No76]{No76} * \emph{С. П. Новиков.} Топология-1. М.: Наука, 1976. (Итоги науки и техники. ВИНИТИ. Современные проблемы математики. Основные направления, 12).}

\newcommand{\nwns}{\bibitem[NW97]{NW97} \emph{A. Nabutovsky, S. Weinberger}. Algorithmic aspects of homeomorphism problems. arXiv:math/9707232.}


\newcommand{\omoe}{\bibitem[Om18]{Om18} * \emph{А. Омельченко,} Теория графов. М.: МЦНМО, 2018.}

\newcommand{\ossf}{\bibitem[OS74]{OS74} \emph{R. P. Osborne and R. S. Stevens.} Group presentations
corresponding to spines of 3-manifolds, I, Amer. J.~Math. 1974. 96. P.~454-471; II, Amer. J.~Math. 1977. 234.
P.~213-243; III, Amer. J.~Math. 1977. 234 P.~245-251.}


\newcommand{\oz}{\bibitem[Oz]{Oz} \emph{M. \"Ozaydin,} Equivariant maps for the symmetric group, unpublished,
\url{http://minds.wisconsin.edu/handle/1793/63829}.}

\newcommand{\panof}{\bibitem[Pan15]{Pan15} \emph{K. Panagiotis.} A note on the topology of irreducible $SO(3)$-manifolds, 	arXiv:1508.06150.}

\newcommand{\paof}{\bibitem[Pa15]{Pa15} \emph{S. Parsa,} On links of vertices in simplicial $d$-complexes embeddable in the euclidean $2d$-space, Discrete Comput. Geom. 59:3 (2018), 663--679.
This is arXiv:1512.05164v4 up to numbering of sections, theorems etc; we refer to numbering in arxiv version.}

\newcommand{\paoe}{\bibitem[Pa18]{Pa18} \emph{S. Parsa,} On links of vertices in simplicial $d$-complexes
embeddable in the euclidean $2d$-space, arXiv:1512.05164v6.}

\newcommand{\patz}{\bibitem[Pa20]{Pa20} \emph{S. Parsa,} On links of vertices in simplicial $d$-complexes
embeddable in the euclidean $2d$-space, arXiv:1512.05164v8.}


\newcommand{\patzl}{\bibitem[Pa20]{Pa20} \emph{S. Parsa,}
Correction to: On the Links of Vertices in Simplicial $d$-Complexes Embeddable in the Euclidean $2d$-Space,
Discrete Comput. Geom. (2020).}

\newcommand{\patza}{\bibitem[Pa20a]{Pa20a} \emph{S. Parsa,} On the Smith classes, the van Kampen obstruction and embeddability of $[3]*K$, arXiv:2001.06478.}

\newcommand{\patzb}{\bibitem[Pa20b]{Pa20b} \emph{S. Parsa,} On the embeddability of $[3]*K$, arXiv:2001.06506.}

\newcommand{\pak}{\bibitem[Pa]{Pa} * \emph{I. Pak}, Lectures on Discrete and Polyhedral Geometry, \url{http://www.math.ucla.edu/~pak/geompol8.pdf}.}

\newcommand{\peze}{\bibitem[Pe08]{Pe08} \emph{Д. Пермяков.} Классификация погружений графов в плоскость,
Вестник МГУ, сер.1, 2008, N5, 55-56.}

\newcommand{\peos}{\bibitem[Pe16]{Pe16} \emph{Д. Пермяков.} Матем. сб., 207:6 (2016),  93--112.}

\newcommand{\pest}{\bibitem[Pe72]{Pe72} * \emph{B. B. Peterson.} The Geometry of Radon's Theorem, Amer. Math. Monthly 79 (1972), 949-963.}


\newcommand{\prnf}{\bibitem[Pr95]{Pr95} * \emph{V. V. Prasolov.} Intuitive topology. Amer. Math. Soc., Providence, R.I., 1995.}

\newcommand{\prnfr}{\bibitem[Pr95]{Pr95} * \emph{В. В. Прасолов.} Наглядная топология. М.: МЦНМО, 1995.}


\newcommand{\przs}{\bibitem[Pr06]{Pr06} * \emph{V. V. Prasolov.}
Elements of Combinatorial and Differential Topology, 2006, GSM 74, Amer. Math. Soc., Providence, RI.}

\newcommand{\przsru}{\bibitem[Pr04]{Pr04} * \emph{В. В. Прасолов.}
Элементы комбинаторной и дифференциальной топологии. М.: МЦНМО, 2004. \url{http://www.mccme.ru/prasolov}.}

\newcommand{\przse}{\bibitem[Pr07]{Pr07} * \emph{V. V. Prasolov.} Elements of homology theory. 2007, GSM 74, Amer. Math. Soc., Providence, RI.}


\newcommand{\przseru}{\bibitem[Pr06]{Pr06} * \emph{В. В. Прасолов.} Элементы теории гомологий. М.: МЦНМО, 2006.}


\newcommand{\psns}{\bibitem[PS96]{PS96} * \emph{V. V. Prasolov, A. B. Sossinsky } Knots, Links, Braids, and 3-manifolds. Amer. Math. Soc. Publ., Providence, R.I., 1996.}


\newcommand{\pszf}{\bibitem[PS05]{PS05} * \emph{В. В. Прасолов и М. Б. Скопенков.}
Рамсеевская теория зацеплений, Мат. Просвещение. 2005. 9. С.~108--115.}

\newcommand{\pszfen}{\bibitem[PS05]{PS05} * \emph{V. V. Prasolov and M.B. Skopenkov.}
Ramsey link theory, Mat, Prosvescheniye, 9 (2005), 108--115.}

\newcommand{\psoo}{\bibitem[PS11]{PS11} \emph{Y. Ponty and C. Saule.} A combinatorial framework for designing (pseudoknotted) RNA algorithms, Proc. of the 11th Intern. Workshop on Algorithms in Bioinformatics, WABI'11, 250--269.}


\newcommand{\pstz}{\bibitem[PS20]{PS20} \emph{S. Parsa and A. Skopenkov.} On embeddability of joins and their `factors', arXiv:2003.12285.}


\newcommand{\psszn}{\bibitem[PSS]{PSS} \emph{M. J. Pelsmajer, M. Schaefer and D. Stasi.} Strong Hanani-Tutte on the projective plane. SIAM J. Discrete Math., 23:3 (2009) 1317--1323.}

\newcommand{\pton}{\bibitem[PT19]{PT19} \emph{P. Pat\'ak and M. Tancer.} Embeddings of $k$-complexes into $2k$-manifolds. arXiv:1904.02404.}

\newcommand{\pw}{\bibitem[PW]{PW} \emph{I. Pak, S. Wilson}, G\lowercase{EOMETRIC REALIZATIONS OF POLYHEDRAL COMPLEXES}, \url{http://www.math.ucla.edu/~pak/papers/Fary-full31.pdf}.}


\newcommand{\razf}{\bibitem[RA05]{RA05} * \emph{J. L. Ram\'irez Alfons\'in.} Knots and links in spatial graphs: a survey. Discrete Math., 302 (2005), 225--242.}

\newcommand{\rep}{\bibitem[Rep]{Rep} Referee's report on the paper ``Some `converses' to intrinsic linking theorems', \url{https://www.mccme.ru/circles/oim/materials/ksreport.pdf}}

\newcommand{\rnoo}{\bibitem[RN11]{RN11} * \emph{R. L. Ricca, B. Nipoti.} Gauss' linking number revisited.
J. of Knot Theory and Its Ramif. 20:10 (2011) 1325--1343. \url{https://www.maths.ed.ac.uk/~v1ranick/papers/ricca.pdf} .}

\newcommand{\rrstz}{\bibitem[RRS]{RRS} * \emph{V. Retinskiy, A. Ryabichev and A. Skopenkov.}
Motivated exposition of the proof of the Tverberg Theorem (in Russian).
Mat. Prosveschenie, 27 (2021), 166--169. arXiv:2008.08361.}

\newcommand{\rssec}{\bibitem[RS68]{RS68} \emph{C. P. Rourke and B. J. Sanderson,} Block bundles II, Ann. of Math. (2), 87 (1968) 431--483.}

\newcommand{\rsst}{\bibitem[RS72]{RS72} * \emph{C. P. Rourke and B. J. Sanderson,}
\newblock Introduction to Piecewise-Linear Topology,
\newblock \emph{Ergebn.\ der Math.} 69, Springer-Verlag, Berlin, 1972.}

\newcommand{\rsstr}{\bibitem[RS72]{RS72} * \emph{К. П. Рурк и Б. Дж. Сандерсон.} Введение в кусочно-линейную топологию, Москва. Мир. 1974.}

\newcommand{\rsns}{\bibitem[RS96]{RS96} * \emph{D. Repov\v s and A. B. Skopenkov.}
Embeddability and isotopy of polyhedra in Euclidean spaces,
Proc. of the Steklov Inst. Math. 1996. 212. P.~173-188.}

\newcommand{\rsne}{\bibitem[RS98]{RS98} \emph{D. Repov\v s and A. B. Skopenkov.}
A deleted product criterion for approximability of a map by embeddings, Topol. Appl. 1998. 87 P.~1-19.}

\newcommand{\rsnn}{\bibitem[RS99]{RS99} * \emph{D. Repov\v s and A. B. Skopenkov.} New results on embeddings of polyhedra and manifolds into Euclidean spaces,
Russ. Math. Surv. 54:6 (1999), 1149--1196.}


\newcommand{\rsnnd}{\bibitem[RS99']{RS99'} * \emph{Д. Реповш и А. Скопенков.}
Кольца Борромео и препятствия к вложимости, Труды МИРАН. 1999. 225. С.~331-338.}

\newcommand{\rszz}{\bibitem[RS00]{RS00} \emph{D. Repov\v s and A. Skopenkov.} Cell-like resolutions of polyhedra by special ones,  Colloq. Math. 2000. 86:2. P. 231--237.}

\newcommand{\rszzd}{\bibitem[RS00']{RS00'} * \emph{Д. Реповш и А. Скопенков.} Характеристические классы для начинающих, Мат. Просвещение. 2000. 4. С.~151-176.}

\newcommand{\rszo}{\bibitem[RS01]{RS01} \emph{D. Repovs and A. Skopenkov.} On contractible $n$-dimensional compacta, non-embeddable into $\R^{2n}$, Proc. Amer. Math. Soc. 129 (2001) 627--628.}

\newcommand{\rszt}{\bibitem[RS02]{RS02} * \emph{Д. Реповш и А. Скопенков.} Теория препятствий для начинающих,
Мат. Просвещение. 2002. 6. C.~60-77.}

\newcommand{\rszf}{\bibitem[RS04]{RS04} \emph{N. Robertson and P. Seymour.} Graph Minors. XX. Wagner's conjecture, J. of Comb. Theory, B, 92:2 (2004) 325--357.}

\newcommand{\rssnf}{\bibitem[RSS]{RSS95} \emph{D. Repov\v s, A. B. Skopenkov  and E. V. \v S\v cepin.}
On uncountable collections of continua and their span, Colloq. Math. 1995. 69:2. P.~289-296.}

\newcommand{\rssnfd}{\bibitem[RSS']{RSS95'} \emph{D. Repov\v s, A. B. Skopenkov and E. V \v S\v cepin.}
On embeddability of $X\times I$ into Euclidean space, Houston J.~Math. 1995. 21. P.~199-204.}

\newcommand{\rssz}{\bibitem[RSS+]{RSSZ} * \emph{A. Rukhovich, A. Skopenkov, M. Skopenkov, A. Zimin},
Realizability of hypergraphs, \url{http://www.turgor.ru/lktg/2013/1/index.htm}.}

\newcommand{\rstnt}{\bibitem[RST']{RST93} \emph{N. Robertson, P. Seymour and R. Thomas}, Linkless embeddings of graphs in 3-space, Bull. of the Amer. Math. Soc., 21 (1993) 84--89.}

\newcommand{\rstno}{\bibitem[RST]{RST91} \emph{N. Robertson, P. Seymour and R. Thomas}, A survey of
linkless embeddings, Graph Structure Theory (Seattle, WA, 1991), Contemp. Math. 147, (1993) 125--136.}


\newcommand{\rwzl}{\bibitem[RWZ+]{RWZ+} \emph{Y. Ren, C. Wen, S. Zhen, N. Lei, F. Luo, D.X. Gu},
Characteristic class of isotopy for surfaces, J. Syst. Sci. Complex. 33 (2020) 2139--2156.}


\newcommand{\saeo}{\bibitem[Sa81]{Sa81} \emph{H. Sachs.} On spatial representation of finite graphs,
in: Finite and infinite sets, Colloq. Math. Soc. Janos Bolyai, North Holland, Amsterdam (37) 1981.}

\newcommand{\sano}{\bibitem[Sa91]{Sa91} \emph{K. S. Sarkaria.}
A one-dimensional Whitney trick and Kuratowski's graph planarity criterion, Israel J.~Math. 73 (1991), 79--89.
\url{http://kssarkaria.org/docs/One-dimensional.pdf}.}

\newcommand{\sanov}{\bibitem[Sa91g]{Sa91g} \emph{K. S. Sarkaria.} A generalized Van Kampen-Flores theorem, Proc. Amer. Math. Soc. 111 (1991), 559--565.}

\newcommand{\sant}{\bibitem[Sa92]{Sa92} \emph{K. S. Sarkaria.} Tverberg’s theorem via number fields. Israel J. Math., 79:317–320, 1992.}

\newcommand{\sazz}{\bibitem[Sa00]{Sa00} \emph{K. S. Sarkaria.} Tverberg partitions and Borsuk-Ulam theorems. Pacific J. Math., 196:1 (2000) 231--241.}

\newcommand{\sczf}{\bibitem[Sc04]{Sc04} \emph{T. Sch\"oneborn.} On the Topological Tverberg Theorem, arXiv:math/0405393.}


\newcommand{\scot}{\bibitem[Sc13]{Sc13} * \emph{M. Schaefer.} Hanani-Tutte and related results. In Geometry --- intuitive, discrete, and convex, Bolyai Soc. Math. Stud., 24 (2013), 259--299.
\url{http://ovid.cs.depaul.edu/documents/htsurvey.pdf} }


\newcommand{\sctz}{\bibitem[Sc20]{Sc20} \emph{M. Schaefer.} The Graph Crossing Number and
its Variants: A Survey. The Electr. J. of Comb. (2020), DS21, \url{https://www.combinatorics.org/files/Surveys/ds21/ds21v5-2020.pdf}}


\newcommand{\scef}{\bibitem[Sc84]{Sc84} \emph{E.~V.~\v S\v cepin.} Soft mappings of manifolds, Russian Math. Surveys, 39:5 (1984).}

\newcommand{\shfs}{\bibitem[Sh57]{Sh57} \emph{A. Shapiro,} Obstructions to the embedding of a complex in a Euclidean space, I, The first obstruction, Ann. Math. 66 (1957), 256--269.}


\newcommand{\shen}{\bibitem[Sh89]{Sh89} * \emph{Ю. А. Шашкин,} Неподвижные точки, М., Наука, 1989.}

\newcommand{\shoe}{\bibitem[Sh18]{Sh18} * \emph{S. Shlosman},  Topological Tverberg Theorem: the proofs and the counterexamples, Russian Math. Surveys, 73:2 (2018), 175–182. arXiv:1804.03120}

\newcommand{\sisn}{\bibitem[Si69]{Si69} \emph{K. Sieklucki.} Realization of mappings, Fund. Math. 1969. 65. P.~325-343.}

\newcommand{\sios}{\bibitem[Si16]{Si16} \emph{S. Simon,} Average-Value Tverberg Partitions via Finite Fourier Analysis, Israel J. Math., 216 (2016), 891-904, arXiv:1501.04612.}



\newcommand{\sknf}{\bibitem[Sk94]{Sk94} \emph{А. Скопенков.} Геометрическое доказательство теоремы
Нойвирта об утолщаемости 2-мерных полиэдров, Math. Notes. 1995. 58:5. P.~1244-1247.}


\newcommand{\skne}{\bibitem[Sk98]{Sk98} \emph{A. B. Skopenkov.} On the deleted product criterion for embeddability in $\R^m$, Proc. Amer. Math. Soc. 1998. 126:8. P.~2467-2476.}

\newcommand{\skzz}{\bibitem[Sk00]{Sk00} \emph{A. Skopenkov,} On the generalized Massey--Rolfsen invariant for link maps, Fund. Math. 165 (2000), 1--15.}

\newcommand{\skzt}{\bibitem[Sk02]{Sk02} \emph{A. Skopenkov,} On the Haefliger-Hirsch-Wu invariants for embeddings and immersions, Comment. Math. Helv. 77 (2002), 78--124.}

\newcommand{\skzth}{\bibitem[Sk03]{Sk03} \emph{M. Skopenkov,} Embedding products of graphs into Euclidean spaces,
Fund. Math. 179 (2003),~191--198, arXiv:0808.1199.}

\newcommand{\skzthd}{\bibitem[Sk03']{Sk03'} \emph{M. Skopenkov,} On approximability by embeddings of cycles in the plane, Topol. Appl. 134 (2003),~1--22, arXiv:0808.1187.}

\newcommand{\skzf}{\bibitem[Sk05]{Sk05} * \emph{A. Skopenkov,}
On the Kuratowski graph planarity criterion, Mat. Prosveschenie, 9 (2005), 116-128. arXiv:0802.3820.}


\newcommand{\skzs}{\bibitem[Sk05i]{Sk05i} \emph{A. Skopenkov,} A new invariant and parametric connected sum of embeddings, Fund. Math. 197 (2007) 253--269. arxiv:math/0509621.}

\newcommand{\skzei}{\bibitem[Sk05]{Sk05} \emph{A.  Skopenkov,} A classification of smooth embeddings of
4-manifolds in 7-space, I, Topol. Appl., 157 (2010) 2094--2110. arXiv:math/0512594.}

\newcommand{\skze}{\bibitem[Sk06]{Sk06} * \emph{A. Skopenkov,} Embedding and knotting of manifolds in Euclidean spaces, London Math. Soc. Lect. Notes, 347 (2008) 248--342. arXiv:math/0604045.}

\newcommand{\skzsi}{\bibitem[Sk06']{Sk06'} \emph{A. Skopenkov,} A classification of smooth embeddings of 3-manifolds in 6-space, Math. Zeitschrift, 260:3 (2008) 647--672. arxiv:math/0603429.}

\newcommand{\skozp}{\bibitem[Sk08]{Sk08} \emph{A.  Skopenkov,} Embeddings of $k$-connected $n$-manifolds into
$\R^{2n-k-1}$. arxiv:math/0812.0263; earlier version published in Proc. Amer. Math. Soc., 138 (2010) 3377--3389.}

\newcommand{\skoz}{\bibitem[Sk10]{Sk10} * \emph{А. Скопенков,} Вложения в плоскость графов с вершинами степени 4,
Мат. Просвещение, 21 (2017), arXiv:1008.4940.}

\newcommand{\skoo}{\bibitem[Sk11]{Sk11} \emph{M. Skopenkov,} When is the set of embeddings finite up to isotopy? Intern. J. Math. 26:7 (2015), 28 pp. arXiv:1106.1878.}

\newcommand{\sks}{\bibitem[Sk14]{Sk14} * \emph{A. Skopenkov,} Realizability of hypergraphs and Ramsey link theory, arXiv:1402.0658.}


\newcommand{\skof}{\bibitem[Sk15]{Sk15} * \emph{А. Скопенков,} Алгебраическая топология с геометрической точки зрения, Москва, МЦНМО, 2015 (1е издание).}

\newcommand{\skofe}{\bibitem[Sk15]{Sk15} * \emph{A. Skopenkov,} Algebraic Topology From Geometric Viewpoint (in Russian), MCCME, Moscow, 2015 (1st edition). }

\newcommand{\skofel}{\bibitem[Sk15e]{Sk15e} * \emph{А. Скопенков,} Алгебраическая топология
с геометрической точки зрения, эл. версия, \url{http://www.mccme.ru/circles/oim/home/combtop13.htm\#photo}}


\newcommand{\skotzr}{\bibitem[Sk20]{Sk20} * \emph{А. Скопенков,} Алгебраическая топология с геометрической точки зрения, Москва, МЦНМО, 2020 (2е издание).
Часть книги: \url{http://www.mccme.ru/circles/oim/obstruct.pdf}}

\newcommand{\skotz}{\bibitem[Sk20]{Sk20} * \emph{A. Skopenkov,} Algebraic Topology From Geometric Viewpoint (in Russian), MCCME, Moscow, 2020 (2nd edition).
Part of the book: \url{http://www.mccme.ru/circles/oim/obstruct.pdf} . Accepted for the English translation by Springer. Preprint of a part: \url{https://www.mccme.ru/circles/oim/obstructeng.pdf}. }

\newcommand{\skofp}{\bibitem[Sk15]{Sk15} \emph{A. Skopenkov,} Classification of knotted tori,
Proc. A of the Royal Soc. of Edinburgh, 150:2 (2020), 549-567. Full version: arXiv:1502.04470.}


\newcommand{\skos}{\bibitem[Sk16]{Sk16} * \emph{A. Skopenkov,} A user's guide to the topological Tverberg Conjecture, Russian Math. Surveys, 73:2 (2018), 323--353.  arXiv:1605.05141.
Section 4 of the published version is available as {\it A. Skopenkov,} On van Kampen-Flores, Conway-Gordon-Sachs and Radon theorems, arXiv:1704.00300.}



\newcommand{\skosd}{\bibitem[Sk16']{Sk16'} * \emph{A. Skopenkov,} Stability of intersections of graphs in the plane and the van Kampen obstruction, Topol. Appl. 240(2018) 259--269, arXiv:1609.03727.}


\newcommand{\skosc}{\bibitem[Sk16c]{Sk16c} * \emph{A. Skopenkov,}  Embeddings in Euclidean space: an introduction to their classification, to appear in Boll. Man. Atl. http://www.map.mpim-bonn.mpg.de/Embeddings\_in\_Euclidean\_space:\_an\_introduction\_to\_their\_classification}

\newcommand{\skosie}{\bibitem[Sk16e]{Sk16e} * \emph{A. Skopenkov,} Embeddings just below the stable range: classification, to appear in Boll. Man. Atl.
http://www.map.mpim-bonn.mpg.de/Embeddings\_just\_below\_the\_stable\_range:\_classification}

\newcommand{\skost}{\bibitem[Sk16t]{Sk16t} * \emph{A. Skopenkov,} 3-manifolds in 6-space, to appear in Boll. Man. Atl.
http://www.map.mpim-bonn.mpg.de/3-manifolds\_in\_6-space}

\newcommand{\skosf}{\bibitem[Sk16f]{Sk16f} * \emph{A. Skopenkov,} 4-manifolds in 7-space, to appear in Boll. Man. Atl. http://www.map.mpim-bonn.mpg.de/4-manifolds\_in\_7-space}

\newcommand{\skosh}{\bibitem[Sk16h]{Sk16h} * \emph{A. Skopenkov,} High codimension links, to appear in Boll. Man. Atl. \url{http://www.map.mpim-bonn.mpg.de/High_codimension_links}.}

\newcommand{\skosi}{\bibitem[Sk16i]{Sk16i} * \emph{A. Skopenkov,} Isotopy, submitted to Boll. Man. Atl.
\url{http://www.map.mpim-bonn.mpg.de/Isotopy}.}

\newcommand{\skose}{\bibitem[Sk17]{Sk17} \emph{A. Skopenkov,}
Eliminating higher-multiplicity intersections in the metastable dimension range, submitted. arxiv:1704.00143.}

\newcommand{\skosed}{\bibitem[Sk17v]{Sk17v} * \emph{A. Skopenkov,}
On van Kampen-Flores, Conway-Gordon-Sachs and Radon theorems,  arxiv:1704.00300.}

\newcommand{\sk}{\bibitem[Sk17o]{Sk17o} \emph{A. Skopenkov,} On the metastable Mabillard-Wagner conjecture.  arxiv:1702.04259.}

\newcommand{\skmos}{\bibitem[Sk17d]{Sk17d} \emph{M. Skopenkov}. Discrete field theory: symmetries and conservation laws, arxiv:1709.04788.}

\newcommand{\skoe}{\bibitem[Sk18]{Sk18} * \emph{A. Skopenkov.} Invariants of graph drawings in the plane.
Arnold Math. J., 6 (2020) 21--55; full version: arXiv:1805.10237.}


\newcommand{\skoeo}{\bibitem[Sk18o]{Sk18o} * \emph{A. Skopenkov.} A short exposition of S. Parsa's theorems on intrinsic linking and non-realizability. Discr. Comp. Geom. 65:2 (2021), 584--585; full version:  arXiv:1808.08363.}


\newcommand{\skona}{\bibitem[Sk19]{Sk19} * \emph{A. Skopenkov,} A short exposition of the Levine-Lidman example of spineless 4-manifolds, arXiv:1911.07330.}

\newcommand{\sktze}{\bibitem[Sk21m]{Sk21m} * \emph{A. Skopenkov.} Mathematics Through Problems: from olympiades and math circles to a profession. Part I. Algebra. Amer. Math. Soc., Providence, 2021. Preliminary version: \url{https://www.mccme.ru/circles/oim/algebra_eng.pdf}}

\newcommand{\sktz}{\bibitem[Sk20u]{Sk20u} * \emph{A. Skopenkov.} A user's guide to basic knot and link theory,
in Topology, Geometry, and Dynamics, Contemporary Mathematics, vol. 772, Amer. Math. Soc., Providence, RI, 2021, pp. 281--309.
Russian version: Mat. Prosveschenie 27 (2021), 128--165. arXiv:2001.01472.}

\newcommand{\sktzo}{\bibitem[Sk20o]{Sk20o} \emph{A. Skopenkov.} On some results of S. Abramyan and T. Panov, arXiv:2005.11152.}

\newcommand{\sktzr}{\bibitem[Sk20e]{Sk20e} * \emph{A. Skopenkov.}
Extendability of simplicial maps is undecidable, 	arXiv:2008.00492.}



\newcommand{\sktzd}{\bibitem[Sk21d]{Sk21d} * \emph{A. Skopenkov.}
On different reliability standards in current mathematical research, arXiv:2101.03745.
More often updated version: \url{https://www.mccme.ru/circles/oim/rese_inte.pdf}.}

\newcommand{\skd}{\bibitem[Sk]{Sk} * \emph{А. Скопенков.} Алгебраическая топология с алгоритмической точки зрения, 
\url{http://www.mccme.ru/circles/oim/algor.pdf}.}

\newcommand{\skde}{\bibitem[Sk]{Sk} * \emph{A. Skopenkov.} Algebraic Topology From Algorithmic Viewpoint, draft of a book, mostly in Russian,
\url{http://www.mccme.ru/circles/oim/algor.pdf}.}


\newcommand{\skon}{\bibitem[Skw]{Skw} * \emph{A. Skopenkov.} Whitney trick for eliminating multiple intersections, slides for talks at St Petersburg, Brno, Kiev, Moscow,  \url{https://www.mccme.ru/circles/oim/eliminat_talk.pdf}.}

\newcommand{\skl}{\bibitem[EEF]{EEF} * {\it Proposed by D. Eliseev, A. Enne, M. Fedorov, A. Glebov, N. Khoroshavkina, E. Morozov, A. Skopenkov, R. \v Zivaljevi\'c.}
A user's guide to knot and link theory, \url{https://www.turgor.ru/lktg/2019/3} .}

\newcommand{\skr}{\bibitem[Skr]{Skr} * \emph{A. Skopenkov.} Realizability of hypergraphs, slides for talks,  \url{https://www.mccme.ru/circles/oim/algor1_beamer.pdf}.}

\newcommand{\skt}{\bibitem[Skt]{Skt} * \emph{A. Skopenkov.} Transparent anonymous peer review,
\url{https://www.mccme.ru/circles/oim/home/transp_peer_review.htm} .}

\newcommand{\rslktg}{\bibitem[KRR+]{RRSl} * Towards higher-dimensional combinatorial geometry, presented by
E. Kogan, V. Retinskiy, E. Riabov and A. Skopenkov, \url{https://www.turgor.ru/lktg/2020/3/index.html}.}

\newcommand{\sm}{\bibitem[Sm]{Sm} S. Smirnov.}

\newcommand{\sper}{\bibitem[Sp]{Sp} * Sperner's lemma defeats the rental harmony problem, \url{https://www.youtube.com/watch?v=7s-YM-kcKME}.}

\newcommand{\sset}{\bibitem[SS83]{SS83} \emph{Е. В. Щепин, М. А. Штанько.} Спектральный критерий вложимости компактов в евклидовы пространства, Труды Ленинградской Международной Топологической конференции. Л.: Наука, 1983. С.~135-142.}

\newcommand{\ssnt}{\bibitem[SS92]{SS92} \emph{J.~Segal and S.~Spie\.z.} Quasi embeddings and embeddings of polyhedra in $\R^m$,  Topol. Appl., 45 (1992) 275--282.}

\newcommand{\sszt}{\bibitem[SS03]{SS03} \emph{F. W. Simmons and F. E. Su.}
Consensus-halving via theorems of Borsuk-Ulam and Tucker, Math. Social Sciences 45 (2003) 15–25. \url{https://www.math.hmc.edu/~su/papers.dir/tucker.pdf}.}

\newcommand{\ssot}{\bibitem[SS13]{SS13} \emph{M. Schaefer and D. Stefankovi\v c.} Block additivity of $\Z_2$-embeddings. In Graph drawing, volume 8242 of Lecture Notes in Comput. Sci., 185--195.
Springer, Cham, 2013. \url{http://ovid.cs.depaul.edu/documents/genus.pdf}}

\newcommand{\sssne}{\bibitem[SSS]{SSS} \emph{J. Segal, A. Skopenkov and S. Spie\. z.}
Embeddings of polyhedra in $\R^m$ and the deleted product obstruction, Topol. Appl. 1998. 85. P.~225-234.}

\newcommand{\sstnf}{\bibitem[SST95]{SST95} \emph{R. S. Simon, S. Spie\. z and H. Toru\'nczyk.}
T\lowercase{HE EXISTENCE OF EQUILIBRIA IN CERTAIN GAMES, SEPARATION FOR FAMILIES OF CONVEX FUNCTIONS
AND A THEOREM OF BORSUK-ULAM TYPE}, Israel J. Math 92 (1995) 1--21.}

\newcommand{\sstzt}{\bibitem[SST02]{SST02} \emph{R. S. Simon, S. Spie\. z and H. Toru\'nczyk.}
E\lowercase{QUILIBRIUM EXISTENCE AND TOPOLOGY IN SOME REPEATED GAMES WITH INCOMPLETE INFORMATION},
Trans. Amer. Math. Soc. 354:12 (2002) 5005-5026.}

\newcommand{\stez}{\bibitem[ST80]{ST80} * {\it H.~Seifert and W.~Threlfall.}
A textbook of topology, v~89 of {\em Pure and Applied Mathematics}.
Academic Press, New York-London, 1980.}


\newcommand{\stzs}{\bibitem[ST07]{ST07} * \emph{А. Скопенков и А. Телишев.}
И вновь о критерии Куратовского планарности графов, Мат. Просвещение, 11 (2007), 159--160.}

\newcommand{\stzse}{\bibitem[ST07]{ST07} * \emph{A. Skopenkov and A. Telishev}, Once again on the Kuratowski graph planarity criterion, Mat. Prosveschenie, 11 (2007), 159-160. arXiv:0802.3820.}

\newcommand{\stos}{\bibitem[ST17]{ST17} \emph{A. Skopenkov  and M. Tancer,}
Hardness of almost embedding simplicial complexes in $\R^d$, Discr. Comp. Geom., 61:2 (2019), 452--463. arXiv:1703.06305.}

\newcommand{\stwh}{\bibitem[SW]{SW} * \url{http://www.map.mpim-bonn.mpg.de/Stiefel-Whitney_characteristic_classes}}

\newcommand{\sz}{\bibitem[SZ05]{SZ} \emph{T. Sch\"oneborn and G. Ziegler}, The Topological Tverberg Theorem and Winding Numbers, J. Comb. Theory, Ser. A, 112:1 (2005) 82--104, arXiv:math/0409081.}

\newcommand{\szno}{\bibitem[Sz91]{Sz91} \emph{A.~Sz\"ucs,} On the cobordism groups of immersions and embeddings,
Math. Proc. Camb. Phil. Soc., 109 (1991) 343--349.}


\newcommand{\ta}{\bibitem[Ta]{Ta} * Handbook of Graph Drawing and Visualization. ed. by R. Tamassia, CRC Press, 2016.}


\newcommand{\tazz}{\bibitem[Ta00]{Ta00} \emph{K. Taniyama,} Higher dimensional links in a simplicial complex embedded in a sphere, Pacific Jour. of Math. 194:2 (2000), 465-467.}

\newcommand{\theo}{\bibitem[Th81]{Th81} * \emph{C.~Thomassen,} Kuratowski's theorem, J.~Graph. Theory 5 (1981), 225--242.}

\newcommand{\tooo}{\bibitem[To11]{To11} \emph{Tonkonog D.} Embedding 3-manifolds with boundary into closed 3-manifolds, Topol. Appl. 158 (2011), 1157-1162. arXiv:1003.3029.}


\newcommand{\tsbzf}{\bibitem[TSB]{TSB} \emph{D. M. Thilikos, M. Serna and H. L. Bodlaender},
Cutwidth I: A linear time fixed parameter algorithm, J. of Algorithms, 56:1 (2005), 1--24.}


\newcommand{\tsbzfd}{\bibitem[TSB05']{TSB05'} \emph{D. M. Thilikos, M. Serna and H. L. Bodlaender},
Cutwidth II: , J. of Algorithms, 56:1 (2005), 25--49.}



\newcommand{\umse}{\bibitem[Um78]{Um78} \emph{B. Ummel.} The product of nonplanar complexes does not imbed in 4-space, Trans. Amer. Math. Soc., 242 (1978) 319--328.}




\newcommand{\val}{\bibitem[Val]{Val} Valknut, \url{https://en.wikipedia.org/wiki/Valknut}}


\newcommand{\vi}{\bibitem[Vi]{Vi} * \emph{O. Viro.}
Some integral calculus based on Euler characteristic, Lect. Notes in Math. 1346.}

\newcommand{\vizt}{\bibitem[Vi02]{Vi02} * \emph{Э. Б. Винберг.} Курс алгебры. Москва. Факториал Пресс. 2002.}

\newcommand{\vizteng}{\bibitem[Vi02]{Vi02} * \emph{E. B. Vinberg.} A Course in Algebra. Graduate Studies in Mathematics, vol. 56. 2003.}

\newcommand{\vinhzs}{\bibitem[VINH07]{VINH07} * \emph{О. Я. Виро, О. А. Иванов, Н. Ю. Нецветаев и В. М. Харламов.}
Элементарная топология, МЦНМО. 2007.}

\newcommand{\vktt}{\bibitem[vK32]{vK32} \emph{E.~R.~van~Kampen}, Komplexe in euklidischen R\"aumen, Abh. Math. Sem. Hamburg, 9 (1933) 72--78; Berichtigung dazu, 152--153.
English translation by Tu T$\hat a$m Ngu$\tilde{\hat e}$n-Phan:
\url{https://sites.google.com/site/tutamnguyenphan/van_Kampen.pdf}}

\newcommand{\kafo}{\bibitem[vK41]{vK41} \emph{E. R. van Kampen,} Remark on the address of S. S. Cairns,
in Lectures in Topology, 311--313, University of Michigan Press, Ann Arbor, MI, 1941.}

\newcommand{\vo}{\bibitem[Vo96]{vo96} \emph{A. Yu. Volovikov,} On a topological generalization of the Tverberg theorem. Math. Notes 59:3 (1996), 324--326.}

\newcommand{\vopns}{\bibitem[Vo96v]{Vo96v} \emph{A. Yu. Volovikov,} On the van Kampen-Flores Theorem.
Math. Notes 59:5 (1996), 477--481.}

\newcommand{\vznt}{\bibitem[VZ93]{VZ93} \emph{A. Vu\v ci\'c and R. T. \v Zivaljevi\'c}, Note on a conjecture of Sierksma, Discr. Comput. Geom. 9 (1993), 339-349.}

\newcommand{\vzzn}{\bibitem[VZ09]{VZ09} \emph{S. T. Vre\'cica and R. T. \v Zivaljevi\'c},  Chessboard complexes
indomitable, J. of Comb. Theory, Ser. A 118:7 (2011), 2157--2166. arXiv:0911.3512.}


\newcommand{\walst}{\bibitem[Wa62]{Wa62} \emph{C.~T.~C.~Wall}, Classification of $(n-1)$-connected $2n$-manifolds, Ann. of Math., 75 (1962) 163--189.}


\newcommand{\wallss}{\bibitem[Wa67]{Wa67} \emph{C.~T.~C.~Wall.} Classification problems in differential topology, IV, Thickenings, Topology 1966. 5. P. 73--94.}

\newcommand{\waldss}{\bibitem[Wa67m]{Wa67m} \emph{F. Waldhausen.} Eine Klasse von 3-dimensional Mannigfaltigkeiten, I. Invent. Math. 1967. 3. P.~308-333.}

\newcommand{\walsz}{\bibitem[Wa70]{Wa70} \emph{C. T. C. Wall,} Surgery on compact manifolds,
1970, Academic Press, London.}

\newcommand{\wess}{\bibitem[We67]{We67} \emph{C.~Weber.} Plongements de poly\`edres dans le domain metastable, Comment. Math. Helv. 42 (1967), 1--27.}

\newcommand{\whit}{\bibitem[Wl]{Wl} * \url{https://en.wikipedia.org/wiki/Whitehead_link}}

\newcommand{\winum}{\bibitem[Wn]{Wn} * \url{https://en.wikipedia.org/wiki/Winding_number}}

\newcommand{\wrss}{\bibitem[Wr77]{Wr77} \emph{P. Wright.} Covering 2-dimensional polyhedra by 3-manifolds spines.
Topology. 16 (1977), 435--439.}

\newcommand{\wufe}{\bibitem[Wu58]{Wu58} \emph{W. T. Wu.} On the realization of complexes in a euclidean space (in Chinese): I, Sci Sinica, 7 (1958) 251--297; II, Sci Sinica, 7 (1958) 365--387; III, Sci Sinica, 8 (1959) 133--150.}

 \newcommand{\wufn}{\bibitem[Wu59]{Wu59} \emph{W.~T.~Wu.} On the isotopy of a finite complex in Euclidean space, I, II, Science Record, N.S. 3:8 (1959) 342--347, 348--351.}

\newcommand{\wusf}{\bibitem[Wu65]{Wu65} * \emph{W. T. Wu.} A Theory of Embedding, Immersion and Isotopy of Polytopes in an Euclidean Space. Peking: Science Press, 1965.}


\newcommand{\yann}{\bibitem[Ya99]{Ya99} \emph{Z. Yang.} Computing Equilibria and Fixed Points: The Solution of Nonlinear Inequalities, Kluwer, Springer Science + Business Media, 1990.}

\newcommand{\z}{\bibitem[Ze]{Z} * \emph{E. C. Zeeman}, A Brief History of Topology, UC Berkeley, October 27, 1993, On the occasion of Moe Hirsch's 60th birthday, \url{http://zakuski.utsa.edu/~gokhman/ecz/hirsch60.pdf}.}

\newcommand{\zioz}{\bibitem[Zi10]{Zi10} * \emph{D. \v Zivaljevi\'c}, Borromean and Brunnian Rings
\url{http://www.rade-zivaljevic.appspot.com/borromean.html}.}

\newcommand{\zioo}{\bibitem[Zi11]{Zi11} * \emph{G. M. Ziegler}, 3N Colored Points in a Plane, Notices of the Amer. Math. Soc., 58:4 (2011), 550-557.}


\newcommand{\zot}{\bibitem[Zi13]{Z13} \emph{A. Zimin.} Alternative proofs of the Conway-Gordon-Sachs Theorems, arXiv:1311.2882.}


\newcommand{\zss}{\bibitem[ZSS]{ZSS} * Элементы математики в задачах: через олимпиады и кружки к профессии
Сборник под редакцией А. Заславского, А. Скопенкова и М. Скопенкова. Изд-во МЦНМО, 2018.
\url{http://www.mccme.ru/circles/oim/materials/sturm.pdf}.}


\newcommand{\zu}{\bibitem[Zu]{Zu} \emph{J. Zung.} A non-general-position Parity Lemma,
\url{http://www.turgor.ru/lktg/2013/1/parity.pdf}.}







\agles
\aksoe
\amsw

\bibitem[Ba21]{Ba21}  {\it S. Basu,} Math Review MR3959859 to \cite{BFZ}, Math Reviews (2021).

\bbsn
\bbzos
\beet
\bfzof
\bfzos

\bibitem[BLZ]{BLZ}
P. V. M. Blagojevi\' c, W. L\"uck and G. M. Ziegler,
Equivariant topology of configuration spaces, J. of Topology, 8 (2015), 414-456,
Arxiv:1207.2852

\bm
\bmzzn
\bren
\bssos
\bsseo

\bibitem[BZ16-1]{BZ16-1} P. V. M. Blagojevi\v c and G. M. Ziegler,
Beyond the Borsuk-Ulam theorem: The topological Tverberg story, arXiv:1605.07321v1.

\bz

\cget
\ckv
\crsot

\dies

\ers

\fktnf
\ffene
\fozf
\fo


\bibitem[Fr15]{Fr15} \emph{F. Frick}, Counterexamples to the topological Tverberg conjecture,
arXiv:1502.00947v1.

\bibitem[Fr15o]{Fr15o} \emph{F. Frick}, Counterexamples to the topological Tverberg conjecture,
Oberwolfach reports, 12:1 (2015), 318--321.

\fros

\bibitem[Fr78]{Fr78}
{\em M. Freedman}, Quadruple points of 3-manifolds in $S^4$, Comment. Math. Helv. 53 (1978), 385-394.

\fstz
\fvto

\gres
\groz
\gssn

\bibitem[Ha62]{Ha62} A.~Haefliger, {\em Knotted $(4k-1)$-spheres in $6k$-space}, Ann. of Math. 75 (1962) 452--466.

\hast
\hkne
\hozs
\hpzn
\hufn

\jvz


\bibitem[Ka]{Ka} G. Kalai, From Oberwolfach: The Topological Tverberg Conjecture is False, `Combinatorics and more' blog post, February 6, 2015, \url{gilkalai.wordpress.com}

\ko
\koon
\koze
\kstz
\ksto

\bibitem[Lo]{Lo} M.~de~Longueville.
\newblock Notes on the topological {T}verberg theorem.
\newblock Discrete Math.\ 247 (2002), no.~1--3, 271--297.
\newblock (The paper first appeared in
Discrete Math. 241 (2001) 207--233, but the original version suffered from serious publisher's typesetting errors.).

\lsne

\mazt
\mezs
\meoo
\meos
\meoe
\mtwoz
\mwofo

\url{https://research-explorer.app.ist.ac.at/download/2159/4735/IST-2016-534-v1%2B1_Eliminating_Tverberg_points_I._An_analogue_of_the_Whitney_trick.pdf}

\mwof
\mwos
\mwosd

\neno

\oz

\pest
\przs
\pszfen
\pton
\pw

\rrstz
\rsst
\rstno

\saeo
\sanov
\sazz
\shoe
\sios
\skde
\skzth
\skze
\sks
\skofe



\bibitem[Sk16-2]{Sk16-2} * \emph{A. Skopenkov,} A user's guide to disproof of topological Tverberg conjecture, arXiv:1605.05141v2.

\bibitem[Sk16-3]{Sk16-3} * \emph{A. Skopenkov,} A user's guide to topological Tverberg conjecture, arXiv:1605.05141v3.

\bibitem[Sk16-4]{Sk16-4} * \emph{A. Skopenkov,} A user's guide to topological Tverberg conjecture, arXiv:1605.05141v4.

\bibitem[Sk18u]{Sk18u} * \emph{A. Skopenkov,} A user's guide to the topological Tverberg Conjecture, Russian Math. Surveys, 73:2 (2018), 323--353. Full updated version: arXiv:1605.05141.

\skose
\sk
\skoe
\sktzd

\bibitem[Skr]{Skr} * \emph{A. Skopenkov,} On referee independence and personell management at Moscow State University (in Russian), \url{https://www.mccme.ru/circles/oim/home/indepen.htm}

\skt
\sszt
\sstnf
\sstzt
\stez
\stos

\tazz

\vo
\vopns
\vznt

\wess

\yann

\bibitem[Zi98]{Zi98}
R. T. \v Zivaljevi\'c. User's guide to equivariant methods in combinatorics. II. Publ. Inst. Math.
(Beograd) (N.S.), 64(78) (1998) 107-132.

\zioo
\zot

\end{thebibliography}
\end{document}